\documentclass[11pt]{amsart}
\usepackage{amsbsy,amssymb,amscd,amsfonts,latexsym,amstext,delarray,
amsmath,graphicx}  \headsep=15pt
\input xypic

\newtheorem{thm}{Theorem}[section]
\newtheorem{prop}[thm]{Proposition}
\newtheorem{cor}[thm]{Corollary}
\newtheorem{lem}[thm]{Lemma}
\newtheorem{defn}[thm]{Definition}

\newtheorem{ex}[thm]{Example}

\newtheorem{ques}[thm]{Question}

\numberwithin{equation}{section}

\newcommand{\ie}{{\it i.e.\/}\ }
\newcommand{\eg}{{\it e.g.\/}\ }
\newcommand{\cf}{{\it cf.\/}\ }

\def\text{\hbox}

\def\bF{{\mathbb F}}
\def\bG{{\mathbb G}}

\def\bT{{\mathbb T}}

\def\bX{{\mathbb X}}

\def\A{{\mathbb A}}
\def\C{{\mathbb C}}
\def\F{{\mathbb F}}

\def\N{{\mathbb N}}
\renewcommand{\P}{{\mathbb P}}
\def\Q{{\mathbb Q}}
\def\R{{\mathbb R}}
\def\Z{{\mathbb Z}}
\def\K{{\mathbb K}}

\def\fr{{\mathfrak r}}
\def\fR{{\mathfrak R}}
\def\fH{{\mathfrak H}}
\def\n{{\mathfrak n}}

\def\ft{{\mathfrak t}}
\def\fI{{\mathfrak I}}
\def\fs{{\mathfrak s}}

\def\cA{{\mathcal A}}
\def\cB{{\mathcal B}}
\def\cC{{\mathcal C}}

\def\cH{{\mathcal H}}
\def\cI{{\mathcal I}}
\def\cJ{{\mathcal J}}

\def\cL{{\mathcal L}}

\def\cP{{\mathcal P}}

\def\cR{{\mathcal R}}
\def\cS{{\mathcal S}}
\def\cT{{\mathcal T}}
\def\cU{{\mathcal U}}

\def\cW{{\mathcal W}}
\def\cX{{\mathcal X}}

\def\cZ{{\mathcal Z}}

\def\Aut{{\rm Aut}}

\def\End{{\rm End}}

\def\Gal{{\rm Gal}}
\def\GL{{\rm GL}}
 
\def\Hom{{\rm Hom}}

\def\Ker{{\rm Ker}}

\def\sign{{\rm sign}}
\def\SL{{\rm SL}}
\def\Spec{{\rm Spec}}
\def\Sp{{\rm Spec}}

\def\Tr{{\rm Tr}}

\title{Cyclotomy and endomotives}
\author{Matilde Marcolli}
\address{Department of Mathematics, California Institute of Technology \\
Mail Code 253-37, 1200 E. California Blvd. \\
Pasadena, CA 91125, USA} \email{matilde\@@caltech.edu}
\date{}

\begin{document}
\maketitle

\begin{abstract}
We compare two different models of noncommutative geometry
of the cyclotomic tower, both based on an arithmetic algebra
of functions of roots of unity and an action by endomorphisms,
the first based on the Bost-Connes (BC) quantum statistical mechanical 
system and the second on the Habiro ring, where the Habiro functions 
have, in addition to evaluations at roots of unity, also
full Taylor expansions.  Both have compatible endomorphisms
actions of the multiplicative semigroup of positive integers.
As a higher dimensional generalization, we consider a crossed
product ring obtained using Manin's multivariable generalizations  
of the Habiro  functions and an action by endomorphisms of
the semigroup of integer matrices with positive determinant. 
We then construct a corresponding class of multivariable BC 
endomotives,  which are obtained geometrically from self maps 
of higher dimensional algebraic tori, and we discuss some of 
their quantum statistical mechanical properties. 
These multivariable BC endomotives are universal 
for (torsion free) $\Lambda$-rings, compatibly with the Frobenius action. 
Finally, we discuss briefly how Habiro's universal Witten--Reshetikhin--Turaev 
invariant of integral homology 3-spheres may relate invariants of 
3-manifolds to gadgets over $\F_1$ and semigroup actions on 
homology 3-spheres to endomotives. 
\end{abstract}

\section{Introduction}\label{IntroSec}

Several of the different current approaches to algebraic geometry over the 
``field with one element" are based on identifying a class of additional data
to be imposed over a ring, or a scheme over $\Z$, for it to be defined over 
the ``underlying" $\F_1$. Among these approaches, the one developed
by Soul\'e \cite{Soule} is based on finiteness and universality conditions on
the structure of cyclotomic points, whereas another approach recently
developed by Borger \cite{Bor} is based on the structure of $\Lambda$-ring
as a descent condition from $\Z$ to $\F_1$. 
In the recent work \cite{CCM2}, it was shown that, in the Soul\'e 
approach to $\F_1$-geometry, the tower of extensions $\F_{1^n}$
with their Galois action are nicely encoded by a noncommutative
space, the Bost--Connes (BC) endomotive. Endomotives are a class of
noncommutative spaces introduced in \cite{CCM} based on
projective limits of Artin motives with semigroup actions, and
morphisms given by correspondences. The BC algebra of \cite{BC}  
is a natural example in this class. 

In this paper we show that the BC endomotive, and some direct
generalizations, are also natural objects from the point of view
of Borger's approach to geometry over $\F_1$, since the endomorphism
action that defines the structure of endomotive is also a consistent lift of
Frobenius morphisms in the sense of $\Lambda$-rings.

The results obtained in this paper fall under the common
theme of comparing different models of noncommutative
algebras associated to the cyclotomic tower with the
action of the positive integers by endomorphisms. 
As argued in \cite{CMR1}, \cite{CMR2}, and \S 3 of \cite{CMbook}, the
Bost--Connes algebra can be regarded as the
noncommutative ring of functions of a noncommutative space
naturally associated to the zero-dimensional Shimura variety
(the cyclotomic tower), when one considers on it the action
of the semigroup of positive integers by endomorphisms. The
latter corresponds to imposing the equivalence relation
of commensurability on 1-dimensional $\Q$-lattices up to
scaling, in the formulation given in \cite{CM1}. Here we
also consider another noncommutative algebra of functions
associated to the same geometric object, which instead
of using the group algebra of $\Q/\Z$, or equivalently the
algebra of continuous functions on $\hat\Z$, to describe 
functions on roots of unity, is based on the Habiro ring
of ``analytic functions of roots of unity'' and on the 
same family of endomorphisms. 
Both the Habiro ring and the abelian part of the
Bost-Connes algebra
describe different aspects of the same underlying geometry
and this serves as a very useful guideline to obtain 
various results on generalizations of the BC endomotive. 
For instance, the 
multivariable generalizations of the Habiro ring constructed
by Manin in \cite{Man} also give rise to a family of 
noncommutative algebras related to actions of endomorphisms
of higher dimensional algebraic tori. 
The resulting noncommutative spaces also
admit a description in terms of endomotives in the sense of
\cite{CCM} and this leads to the construction of 
multivariable BC endomotives. 
We show that these multivariable BC algebras are 
closely related to the
theory of $\Lambda$-rings, in the formulation given in 
\cite{BoSm}, as another approach to 
the theory of geometry over the ``field with one element''.
We show that the (multivariable) BC endomotives are universal
(torsion free) $\Lambda$-rings, in the sense that any such
$\Lambda$-ring that is of finite rank over $\Z$ and has no
nilpotents can be embedded in a product of 
BC endomotives, compatibly with the Frobenius action.

In order to construct quantum statistical mechanical
systems for the multivariable BC endomotives, one needs to
correct for the fact that the natural time evolution that
generalizes the one-variable case has infinite multiplicities
in the spectrum of the Hamiltonian coming from the presence
of an $\SL_n(\Z)$-symmetry. The problem is completely 
analogous to the one encountered in \cite{CM1} in the case
of 2-dimensional $\Q$-lattices, and one can use essentially
the same method used in \cite{CM1} to circumvent the problem
by using the trace of a type II$_1$ factor associated to the 
$\SL_n(\Z)$-action. We show that multiple zeta values of 
cones define classes of states on the multivariable BC algebras.
 
The Habiro ring of functions was originally introduced (see
\cite{Hab1} and \cite{Hab2}) as the natural receptacle for
a universal invariant of integral homology 3-spheres that
specializes to the Witten--Reshetikhin--Turaev invariant
$\ft_\zeta(M)$ at each root of unity $\zeta$ and whose
Taylor expansion at $\zeta=1$ gives the Ohtsuki series $\ft^O(M)$. 
In view of the relation between the Habiro ring and the 
noncommutative geometry of the
cyclotomic tower, it is then natural to try to relate the
universal Witten--Reshetikhin--Turaev invariant itself to
quantum statistical mechanics and noncommutative geometry, on
one side, and to the geometry of gadgets over the
``field with one element'' in the sense of \cite{Soule}, on
the other. We concentrate on these questions in the last 
part of the paper, where we consider the question of
constructing convolution algebras associated to surgery
presentations of 3-manifolds with actions of $\N$ by
endomorphisms, and of possible relations to endomotives
via the WRT invariants.

\bigskip

{\bf Acknowledgment.}  
I thank Yuri Manin for useful  
conversations and for suggesting the possible 
relation to $\Lambda$-rings. I thank
Jack Morava for several useful discussions and
Peter Teichner for useful comments. 
I also thank James Borger  for reading earlier
drafts of this manuscript.  I especially wish to thank 
Alain Connes for extensive comments and suggestions 
that greatly improved the paper. This work is partially 
supported by NSF grant DMS-0651925. Part of this work 
was done during stays at the MPI and at MSRI, which 
I thank for the hospitality and for support.

\section{The Bost--Connes system and Habiro 
analytic functions}\label{HabiroBCSec}

We recall here briefly the main properties of 
the Bost--Connes (BC) quantum statistical 
mechanical system introduced in 
\cite{BC}, in the formulation given in \cite{CM1} and
in \S 3 of \cite{CMbook}. We then also recall the
main properties of the Habiro ring of analytic functions
of roots of unity, and we show that the same semigroup
of endomorphisms used to for the crossed product BC
algebra acts on the Habiro ring and gives rise to a
similar crossed product construction. We then show 
that using Habiro functions together with evaluations
at roots of unity one can recover the Hilbert space
representations of the BC algebra, as well as versions
where the full Taylor expansions at roots of unity 
is used.

\subsection{The Bost--Connes endomotive, the
cyclotomic tower, and $\bF_1$}\label{BCSec}

The cyclotomic tower is defined by the direct
system of inclusions of rings $\Z[\zeta_n] \hookrightarrow
\Z[\zeta_m]$, for primitive roots of unity with $n|m$,
with covering groups $\GL_1(\hat\Z)\cong \Gal(\Q^{ab}/\Q)$.
Under the identification $\bar\Q^*_{tors}\cong \Q/\Z$,
one can identify the algebra of functions of roots of 
unity with 
\begin{equation}\label{QZhatZ}
C^*(\Q/\Z) \cong C(\hat\Z),
\end{equation}
where the isomorphism is Pontrjagin duality between
$\Q/\Z$ and $\hat\Z$. More precisely, when one views it
as the group ring $C^*(\Q/\Z)$ this is an algebra of
convolution of functions with finite support, and only
after passing to $C(\hat\Z)$ one has an algebra of
functions with pointwise product. In addition to the automorphisms
of the cyclotomic tower given by the action of $\hat\Z^*$,
there are endomorphisms given by the multiplicative
semigroup of positive integers, acting on the generators
$e(r)$ of the group algebra $\Q[\Q/\Z]$ by
\begin{equation}\label{rhoacter}
\rho_n(e(r))= \frac{1}{n} \sum_{ns=r} e(s).
\end{equation}
As shown in \cite{CM1}, when one identifies the
algebra $C(\hat\Z)$ with functions on the set of
1-dimensional $\Q$-lattices up to scaling, the
action \eqref{rhoacter} implements the equivalence
relation of commensurability. 

The noncommutative BC algebra is explicitly given in terms 
of generators and relations as the algebra (over $\Q$)
generated by elements $e(r)$, with $r\in \Q/\Z$,
and $\mu_n$, with $n\in \N$, subject to the relations
\begin{equation}\label{relsmuBC}
\begin{array}{ll}
(\mu_n^*)^* = \mu_n \\
\mu_n\mu_m=\mu_{nm} & \forall n,m\in \N \\
\mu_n \mu_m^* = \mu_m^* \mu_n & \text{ when }  (n,m)=1\\
\mu_n^* \mu_n =1
\end{array}
\end{equation}
and
\begin{equation}\label{relserBC}
e(r+s)=e(r)e(s), \ \ \ e(0)=1
\end{equation}
and the relation
\begin{equation}\label{relermun}
\mu_n e(r) \mu_n^* =\frac{1}{n} 
\sum_{ns=r} e(s),
\end{equation}
which implements the semigroup action \eqref{rhoacter}
and identifies the resulting algebra with the semigroup 
crossed product
\begin{equation}\label{BCalg}
\cA_{\Q,BC} = \Q[\Q/\Z] \rtimes_\rho \N.
\end{equation}
The algebra $\cA_\Q\otimes_\Q \C$ has as $C^*$-algebra
completion the crossed product algebra
\begin{equation}\label{BCstar}
\cA_{BC} = C^*(\Q/\Z)\rtimes \N=C(\hat\Z)\rtimes \N.
\end{equation}

One can also consider ``categorifications" of the BC system
\cite{Mor} by describing the BC algebra as the convolution
algebra associated to a small category. We take up this point
of view in the last section of this paper, when we discuss
algebras and categories associated to surgery presentations
of 3-manifolds, \cf also \cite{MaZa}. 

As argued in \cite{CMR1}, \cite{CMR2}, the BC algebra 
is (up to Morita equivalence) the algebra of functions 
on the noncommutative zero-dimensional Shimura variety
$C_0(\A_f) \rtimes \Q^*_+$,
whose set of classical points is the usual 
zero-dimensional Shimura variety $\A_f^*/\Q^*_+ =\hat\Z^*$
associated to the cyclotomic tower. This classical
space is recovered from the noncommutative space 
as the set of low temperature extremal KMS states 
of the BC quantum statistical mechanical system.
The action of $\Q^*_+$ on the finite adeles $\A_f$
corresponds to a partially defined action on $\hat\Z$,
which in turn corresponds to the semigroup action 
\eqref{rhoacter}.
This point of view based on towers and Shimura 
varieties was generalized to 2-dimensions in \cite{CM1} 
and \cite{CMR1}, \cite{CMR2}, and to more general 
Shimura varieties in \cite{HaPau}.

The quantum statistical mechanical procedure of
associating to a noncommutative space with a 
natural time evolution a classical space given by
the set of low temperature extremal KMS states was
formulated in greater generality in \cite{CCM}
in the context of a theory of {\em endomotives}.
The latter are noncommutative spaces obtained
from towers of Artin motives and semigroups of
endomorphisms. They carry a Galois action and
a time evolution to which quantum statistical
mechanical techniques can be applied.

The original case of the BC system is obtained 
in \cite{CCM} as a particular case of a general 
procedure that constructs endomotives from self maps
of pointed algebraic varieties, via a projective 
system of Artin motives obtained as inverse
images of the fixed point under the given family 
of self-maps. In the case of the BC endomotive, 
the variety is $\bG_m$, the multiplicative group, 
and the self maps are the semigroup $S=\N$ of morphisms
\begin{equation}\label{ttk}
s_k: P(t,t^{-1}) \mapsto P(t^k, t^{-k}), \ \ \ \ k\in \N,
\end{equation}
on $\Q[t,t^{-1}]$. The point $1$ is fixed by all 
these maps and the inverse image under
the map $t\mapsto t^k$ is 
$X_k=\Spec(\Q[t, t^{-1}]/(t^k-1))$. 
The BC endomotive is obtained, at the algebraic 
level, by taking the projective limit 
$X=\varprojlim_k X_k$ or equivalently the direct 
limit of the corresponding algebras 
$A=\varinjlim_k A_k$, with
$A_k=\Q[\Z/k\Z]$ and the crossed product 
$\cA_{\Q,BC}=A\rtimes S=\Q[\Q/\Z]\rtimes \N$.

The endomorphisms \eqref{ttk} correspond (\cf 
\cite{CCM2}) to the endomorphisms
of $\Q[\Q/\Z]$ given by
\begin{equation}\label{sigmanBC}
\sigma_n(e(r)) = e(nr),  \ \ \ \  \forall n\in \N, 
\ \forall r\in \Q/\Z .
\end{equation}
These satisfy 
\begin{equation}\label{sigmarhorels}
\sigma_n \rho_n (x)= x, \ \ \ \text{ and } \ \ \ 
\rho_n \sigma_n (x) = e_n x, \ \ \ \forall x\in \Q[\Q/\Z],
\end{equation}
with the $\rho_n$ as in \eqref{rhoacter}. Here $e_n$ is 
an idempotent in $\Q[\Q/\Z]$ that gives the 
projection onto the range of $\rho_n$. In particular,
one has $\sigma_n(x)=\mu_n^*\, x\, \mu_n$.

It is useful to recall here more in detail
the properties of the endomorphisms $\sigma_n$,
$\rho_n$, and the idempotents $e_n$, because this
will constitute an important difference
between the model of the noncommutative geometry
of the cylotomic tower based on the BC algebra
and the one based on the Habiro ring. The idempotent
$e_n$ is given by
\begin{equation}\label{enmun}
e_n = \mu_n \mu_n^*,
\end{equation} 
and it satisfies the identity
\begin{equation}\label{enes}
e_n = \frac{1}{n} \sum_{nr=0} e(r),
\end{equation}
where one can see easily that the right-hand-side
does not depend on the choice of a solution 
$r\in \Q/\Z$ of $nr=0$. Thus, in particular, 
one has $e_n\in \Q[\Q/\Z]$.

We also recall one more property of the BC algebra,
namely its Hilbert space representations associated
to embeddings in $\C$ of the roots of unity, which we
need in the following. 

As shown in \cite{BC}, the choice of an 
element $\rho\in \hat\Z^*$ determines a 
representation $\pi_\rho$ of the BC algebra 
as bounded operators on the
Hilbert space $\ell^2(\N)$ with
\begin{equation}\label{munOp}
\mu_n \, \epsilon_k = \epsilon_{nk}
\end{equation}
on the canonical basis $\{ \epsilon_k \}$ of
$\ell^2(\N)$ and 
\begin{equation}\label{erOp}
\pi_\rho(e(r))\, \epsilon_k=\zeta_r^k\, \epsilon_k,
\end{equation}
where $\zeta_r =\rho(e(r))$ is the root of unity
obtained using the chosen $\rho\in \hat\Z^*$ to 
embed the abstract roots of unity $\Q/\Z$ in $\C$. 

Finally, we also recall the fact that an integer
model for the BC algebra was constructed in \cite{CCM2},
as the ring $\cA_{\Z,BC}$ generated by elements $e(r)$,
with $r\in \Q/\Z$ and $\tilde\mu_n$ and $\mu_n^*$,
subject to the same relations \eqref{relserBC}, but with
the relations \eqref{relsmuBC} replaced by the analogous
relations
\begin{equation}\label{relstildemu}
\begin{array}{ll}
\tilde\mu_n \tilde\mu_m =\tilde \mu_{nm} & \forall n,m \in \N \\
\mu_n^* \mu_m^* = \mu_{nm}^* & \forall n,m \in \N \\
\mu_n^* \tilde\mu_n = n & \forall n\in\N \\
\tilde\mu_n \mu_m^* = \mu_m^* \tilde\mu_n & \forall (n,m)=1.
\end{array}
\end{equation}
The endomorphisms $\sigma_n$
of \eqref{sigmanBC} are still defined on $\Z[\Q/\Z]$ and they satisfy 
relations replacing \eqref{relermun} and \eqref{sigmarhorels}, in the form
\begin{equation}\label{sigmanZrels}
\mu_n^* x= \sigma_n(x)  \mu_n^* \ \ \ \text{ and } \ \ \  x \tilde\mu_n = \tilde\mu_n \sigma_n(x). 
 \ \ \ \forall x\in \Z[\Q/\Z]. 
\end{equation}
The $\tilde\mu_n$ can be realized as bounded operators on $\ell^2(\N)$ acting by
$$ \tilde\mu_n \epsilon_m = n \epsilon_{nm}. $$
They no longer are adjoints to the $\mu_n^*$ but they satisfy $\tilde\mu_n = n\mu_n$. 
One then sets (\cite{CCM2})
\begin{equation}\label{tilderhoBC}
 \tilde\rho_n(x) = \tilde\mu_n \, x \, \mu_n^*, 
\end{equation} 
to replace the original relation $\rho_n(x) = \mu_n \, x \, \mu_n^*$, which
involves denominators. The $\tilde\rho_n$ are no longer ring
homomorphisms, but only correspondences, so that $\cA_{\Z, BC}$
is not a semigroup crossed product ring.  

It was also shown in \cite{CCM2} that the Bost--Connes endomotive has a model over the
``field with one element" $\bF_1$, in the sense that it is obtained from a projective limit of 
affine varieties $\mu^{(n)}$ over $\bF_1$ (in the sense of Soul\'e \cite{Soule}, as well as 
in the sense of \cite{CoCons}),  with the endomorphisms $\sigma_n$ of
\eqref{sigmanBC} defining morphisms of ``gadgets" over $\bF_1$ (again in the terminology
of \cite{Soule}). In fact, the abelian part of the Bost--Connes endomotive corresponds to
the inductive system of ``extensions" $\bF_{1^n}$ of $\bF_1$, and the endomorphisms
$\sigma_n$ induce the Frobenius correspondence on all the associated positive
characteristic reductions of the abelian part of the Bost--Connes algebra $\Z[\Q/\Z]\otimes_\K \K$,
with ${\rm char}(\K)=p>0$.

\subsection{The Habiro ring, endomorphisms, and Taylor series expansions}\label{HabiroSec}

We now recall the basic properties of the Habiro ring of ``analytic functions on roots of
unity" constructed in \cite{Hab1}. We then show that the same family of endomorphisms
$\sigma_n$ acting on the abelian part of the BC algebra also act on the Habiro ring
in a compatible way, where the compatibility is given by the evaluations of Habiro functions
at roots of unity. We introduce a semigroup crossed product ring obtained using the action
of $\N$ on the Habiro ring, that provides an analog of the BC algebra.
We discuss the effect of the action of the endomorphisms $\sigma_n$
on the Taylor expansions of Habiro functions at roots of unity.

\medskip

As in \cite{Hab1}, \cite{Hab2}, \cite{Man}, we consider the ring
$\widehat{\Z[q]}$ defined 
as the inverse limit
\begin{equation}\label{hatZq}
\widehat{\Z[q]} = \varprojlim_n \Z[q]/ ((q)_n),
\end{equation}
where
\begin{equation}\label{qn}
(q)_n = (1-q)(1-q^2)\cdots (1-q^n).
\end{equation}
These are ordered by divisibility, since $(q)_k | (q)_n$ for $k\leq
n$, with corresponding inclusions 
$((q)_n)\subseteq ((q)_k)$ of ideals and homomorphisms $\Z[q]/((q)_n)
\twoheadrightarrow  
\Z[q]/((q)_k)$.

Thus, one can think of elements in the Habiro ring as series $\sum_n P_n(q)$,
with $P_n(q)$ in $((q)_n)$. Notice that can be seen as functions of roots of
unity, with the pointwise product, instead of the convolution product of the
group ring $\Z[\Q/\Z]$.

\smallskip

Let $\cZ$ denote the set of all roots of unity in $\C$. By the results
of \cite{Hab1}, \cite{Hab2}, 
for any $\zeta \in \cZ$, there exists an evaluation map 
\begin{equation}\label{evzeta}
ev_\zeta : \widehat{\Z[q]} \to \Z[\zeta]
\end{equation}
that is a surjective ring homomorphisms. Moreover, assembling these maps for
all roots of unity one obtains an {\em injective} 
\begin{equation}\label{evtotal}
ev: \widehat{\Z[q]} \to \prod_{\zeta \in \cZ} \Z[\zeta], 
\end{equation}
that is, elements of $\widehat{\Z[q]}$ are determined uniquely by
their evaluations at roots of unity. This is proved in \cite{Hab1}, Theorem 6.3, 
along with a characterization of the rings of coefficients $R$ for which the
corresponding Habiro ring $\widehat{R[q]}$ has the property that
$ev: \widehat{R[q]} \to \prod_\zeta R[\zeta]$ is injective (see also \S 2 of
\cite{Man}).

\subsubsection{Endomorphisms}

We then have the following result on the action of $\N$ by 
endomorphisms of the Habiro ring.

\begin{prop}\label{semiact}
The semigroup $\N$ of positive integers acts by endomorphisms of $\widehat{\Z[q]}$ by
\begin{equation}\label{sigmanZq}
\sigma_n (f) (q) = f(q^n).
\end{equation}
\end{prop}

\proof First notice that the homomorphism $\sigma_n: \Z[q] \to \Z[q]$ given by
$\sigma_n(f)(q) = f(q^n)$ descends to a homomorphism 
$\sigma_n : \widehat{\Z[q]}\to \widehat{\Z[q]}$. In fact, notice that $(q)_m | \sigma_n(q)_m$
for all $n$ and $m$, hence the homomorphism $\sigma_n$ induces compatible
homomorphisms $\sigma_n: \Z[q]/(q)_m \to \Z[q]/(q)_m$. We still denote by $\sigma_n$
the map on the inverse limit. 
\endproof

Suppose then that $\zeta$ is an $m$-th root of unit . The evaluation
map $ev_\zeta$ is in this case defined on $\Z[q]/((q)_m)$ and one sees that the $\sigma_n$
then give compatible maps $\Z[q]/((q)_m) \to \Z[q]/((q)_{m/n})$, for $n|m$, lifting the map 
$\Z[\zeta]\to \Z[\zeta^n]$ given by $\sigma_n(f)(\zeta)=f(\zeta^n)$. See Corollary \ref{notzetan} below.

\subsubsection{Semigroup crossed product}

By analogy with the BC system, we would then like to form a semigroup crossed
product of the Habiro ring by the action of $\N$ by endomorphisms. This requires
introducing an analog of the $\rho_n$ and of the isometries $\mu_n$ implementing
them as additional generators of the crossed product ring. 

\begin{defn}\label{ZhatinftyDef}
In the following we let $\widehat{\Z[q]}_\infty$ denote the direct limit of the system of 
maps $\sigma_n:\widehat{\Z[q]} \to \widehat{\Z[q]}$.
\end{defn}

We can then define a ring generated by the functions of the Habiro ring
$\widehat{\Z[q]}$ and the action of the endomorphisms $\sigma_n$ in the following way.

\begin{defn}\label{AZqinfty} Let $\cA_{\Z,q}$
denote the ring generated by the Habiro ring $\widehat{\Z[q]}$ together with additional
generators $\mu_n$ and $\mu_n^*$, for $n\in \N$, subject to the relations \eqref{relsmuBC}
and the additional relations 
\begin{equation}\label{relqsigmamu}
\mu_n\sigma_n(f) = f \mu_n , \ \ \ \ \  \mu_n^* f = \sigma_n(f) \mu_n^* ,
\end{equation}
for all $f\in \widehat{\Z[q]}$ and all $n\in \N$, with $\sigma_n(f)(q) = f(q^n)$, as above.
We also require that the elements $\mu_n\mu_n^*$ belong to the direct limit $\widehat{\Z[q]}_\infty$.
\end{defn}

We can identify the ring $\cA_{\Z,q}$ with a crossed product ring.
To this purpose, we make the following preliminary observation.

\begin{lem}\label{hatRnprop} 
Let $\hat R_n$ denote the range of $\sigma_n$ acting on
$\widehat{\Z[q]}$, that is,  
$\hat R_n =\{ f\in \widehat{\Z[q]}\,|\, \exists h\in \widehat{\Z[q]},
\, f =\sigma_n(h) \}$. 
There is a ring homomorphism $\eta_n: \hat R_n \to \widehat{\Z[q]}$,
given by $\eta_n (f)=h$, which satisfies
$\sigma_n\circ \eta_n =id|_{\hat R_n}$ and $\eta_n \circ \sigma_n = id
|_{\widehat{\Z[q]}}$. 
\end{lem} 

\proof It suffices to see that $\hat R_n$ is indeed a subring of
$\widehat{\Z[q]}$ and that 
$\eta_n (f)=h$ gives a well defined ring homomorphism $\hat R_n \to
\widehat{\Z[q]}$.  To this purpose 
it suffices to observe that, if we can represent in two different ways
$f=\sigma_n(h_1)=\sigma_n(h_2)$, then the image of $h_1-h_2$ in
$\Z[q]/((q)_N)$ lies in the ideal $\cI_N$ 
generated by $(q)_N$. Thus, the map $\eta_n (f)=h$ is well defined in
$\widehat{\Z[q]}$. 
\endproof

The statement above is just giving the {\em injectivity} of the maps $\sigma_n$.
In fact, the condition $\Ker \sigma_n =0$ follows from the fact that an element
in the kernel would evaluate at zero on all the $\zeta^n$. The injectivity of the
evaluation \eqref{evtotal} together with the divisibility of  the group of roots of
unity give the result.
This should be contrasted with the situation of the BC system, where
the $\sigma_n$ are surjective and not injective.

\smallskip

We cannot directly extend the partial inverse $\eta_n$ of $\sigma_n$
to a homomorphism 
on the whole of $\widehat{\Z[q]}$, since, unlike what happens in
the case of the Bost--Connes algebra, we do not have an idempotent 
in $\widehat{\Z[q]}$ analogous to the idempotent 
$\mu_n\mu_n^* = e_n \in \Q[\Q/\Z]$. However, we can give a more
detailed description of the ring $\cA_{\Z,q}$ and we find that it is indeed
a crossed product ring. We need some preliminary results.

\begin{lem}\label{ZqLem1}
Consider the rings $A_N=\widehat{\Z[q]}_N$ generated by all
the elements of the form $\mu_N\, f\, \mu_N^*$, for $f\in  \widehat{\Z[q]}$, 
with $\widehat{\Z[q]}_1 =\widehat{\Z[q]}$. These subrings 
$A_N=\widehat{\Z[q]}_N\subset A_\infty =\cup_N A_N$  
satisfy the property that $A_N \cdot A_M \subset A_{NM/(N,M)}$.
Moreover, the endomorphisms
$\sigma_n$ of $A_1=\widehat{\Z[q]}$ extend to endomorphisms
$\sigma_n: A_N \to A_N$ when $n\not| N$
and $\sigma_n : A_N \to A_{N/n}$ when $n|N$. 

There are also well defined endomorphisms $\rho_n(a) =
\mu_n \, a\, \mu_n^*$ mapping $\rho_n: A_N \to A_{nN}$.   
These satisfy 
\begin{equation}\label{rhonsigmanen}
\rho_n(\sigma_n(a)) = e_n a e_n, \ \ \ \  \sigma_n(\rho_n(a))=a
\end{equation}
where the idempotent $e_n =\mu_n \mu_n^*$  maps $A_1$ by $f\mapsto e_n f e_n$ 
onto the subring $\mu_n \hat R_n \mu_n^* \subset A_n$, and more generally it maps
$$ A_N \to \mu_{nN/(n,N)} \, \hat R_{n/(n,N)} \, \mu_{nN/(n,N)}^* \subset A_{nN/(n,N)}, 
\ \ \ \   a_N \mapsto  e_n a_N e_n, $$
where the $\hat R_n$ are as in Lemma \ref{hatRnprop}.
\end{lem}

\proof Each $A_N=\widehat{\Z[q]}_N$ is a ring, by the property that 
$\mu_N^* \mu_N =1$. Moreover, the relations \eqref{relstildemu} also
imply the statement about multiplication in $\cA_{\Z,q}$ of two elements 
$a_N =\mu_N f \mu_N^* \in A_N$ and $a_M =\mu_M h \mu_M^* \in A_M$. 
In fact, this gives $a_N \cdot a_M = \mu_N f \mu_N^* \mu_M h \mu_M^*$. 
If $(N,M)=1$ this is 
$$ \mu_N f \mu_N^* \mu_M h \mu_M^* =\mu_N f \mu_M \mu_N^* h \mu_M^*
=\mu_N \mu_M \sigma_M(f) \sigma_N(h) \mu_N^* \mu_M^* \in A_{NM}. $$
If instead $N=a N'$ and $M=a M'$ with $(N',M')=1$, we find, using $\mu_a^* \mu_a=1$,
$$ \mu_N f \mu_N^* \mu_M h \mu_M^* =\mu_N f \mu_{M'} \mu_{N'}^* h \mu_M^*
=\mu_N \mu_{M'} \sigma_{M'}(f) \sigma_{N'}(h) \mu_{N'}^* \mu_M^* $$ $$ =
\mu_{a N' M'}   \sigma_{M'}(f) \sigma_{N'}(h) \mu_{a N' M'}^* \in A_{a N' M'}. $$

Clearly $\rho_n(a) = \mu_n \, a\, \mu_n^*$ defines an endomorphism
since $\mu_n^* \mu_n = 1$.  Consider the restriction of $\rho_n$ to
$\hat R_n \subset A_1$. Elements of $\hat R_n$ are of the form $\sigma_n(h)$ for some
$h\in A_1$. By the relations between $\sigma_n$ and $\mu_n$ and $\mu_n^*$,
we have $\sigma_n(h)=\mu_n^* h \mu_n$, so that the action of $\rho_n$ on
elements of this form is given by $\mu_n \mu_n^* h \mu_n \mu_n^*$. The
relations \eqref{rhonsigmanen} follow from the fact that $\sigma_n(a) =\mu_n^* a \mu_n$
and $\rho_n(a) =\mu_n a \mu_n^*$. 

Moreover, for $a\in A_N$, $a=\mu_N f \mu_N^*$ for some $f\in A_1$, 
we have $$\rho_n(\sigma_n(a))=\mu_{nN/(n,N)} \sigma_{n/(n,N)}(f) \mu_{nN/(n,N)}^* =
e_n a e_n,$$ which shows that the ring $e_n A_N e_n$, obtained 
by compressing $A_N$ with the idempotent $e_n=e_n^2$, is the subring 
$\mu_{nN/(n,N)} \hat R_{n/(n,N)}  \mu_{nN/(n,N)}^*$ of $A_{nN/(n,N)}$. In the case
of elements $f\in \hat R_n \subset A_1$, the identification 
$e_n A_1 e_n = \mu_n \hat R_n \mu_n^*$ is induced by $\rho_n(h)= e_n \eta_n(h) e_n$, with
$\eta_n: \hat R_n \to A_1$ the homomorphism of Lemma \ref{hatRnprop}.
\endproof

We then consider then the ring $A_\infty =\cup_N  
A_N$ and the induced maps $\sigma_n$ and $\rho_n$.
We have the following result.

\begin{lem}\label{ZqLem2}
The ring $A_\infty =\cup_N  A_N$ is the direct limit $\widehat{\Z[q]}_\infty$ of the
inductive system given by the endomorphisms $\sigma_n : \widehat{\Z[q]} \to \widehat{\Z[q]}$.
The induced morphisms $\sigma_n: \widehat{\Z[q]}_\infty \to \widehat{\Z[q]}_\infty$ are 
automorphisms with inverses the $\rho_n$.
\end{lem}

\proof Let $\widehat{\Z[q]}_\infty$ be the direct limit of the system of injective homomorphisms
$\sigma_n : \widehat{\Z[q]} \to \widehat{\Z[q]}$. The rings $A_N$ are arranged by inclusions
$j_{M,N}: A_N \hookrightarrow A_M$, when $N|M$ with the embedding given by identifying
$a=\mu_N f \mu_N^*$ with the same element written as $a= \mu_M \sigma_k(f) \mu_M^*$,
with $k=M/N$. Thus, we can identify $A_\infty$ with the direct limit of the $A_N$
under these inclusions. The ring homomorphisms $\rho_N: A_1 \to A_N$ with 
$f\mapsto \rho_N(f)=\mu_N f \mu_N^*$ define maps of direct systems, since by construction
$\rho_M \sigma_k(f) = j_{M,N} \rho_N(f)$ for $M=kN$.  This induces a morphism on the 
direct limits $\widehat{\Z[q]}_\infty\to A_\infty$. We see that this is in fact an isomorphism,
since one can construct an inverse to this map using the $\sigma_N : A_N \to A_1$ and their
compatibility with the direct systems. In fact, using the assumption that the elements
$e_n =\mu_n\mu_n^*$ belong to the direct limit $\widehat{\Z[q]}_\infty$, upon writing the
idempotent $e_n =1 -p_n$ for some projector $p_n$ and using the relations induced on
$\widehat{\Z[q]}_\infty$ by \eqref{relqsigmamu}, we obtain 
$\sigma_n(p_n)=\mu_n^*(1-\mu_n\mu_n^*)\mu_n=0$, hence by the injectivity of the $\sigma_n$
we have $e_n=1$. This in fact shows that the $\sigma_n$ become 
automorphisms on the direct limit $\widehat{\Z[q]}$ with inverses the $\rho_n$. One can
see this also from the fact that one can write the direct limit $\widehat{\Z[q]}_\infty$ as sequences
of elements $f_\ell \in \widehat{\Z[q]}$ where $f_{k\ell} = \sigma_k(f_\ell)$ and one sees
in this way that the maps induced by the $\sigma_k$ on the direct limit are surjective as well
as injective.
\endproof

The fact that the endomorphisms $\sigma_n : \widehat{\Z[q]}\to \widehat{\Z[q]}$ induce
automorphisms of the ring $\widehat{\Z[q]}_\infty$ can also be seen from the fact that,
dually, one is taking a projective limit of a system of surjective maps of spaces. The map
induced on the projective limit is then injective as well as surjective, as one sees already 
in the case of a single surjective map $s: X\to X$, where the projective limit is the space
of sequences $(x_n)$ with $x_n= s(x_{n+1})$ and the map $(x_n) \mapsto (s(x_n))$ is
then injective since $(s(x_n))=(s(y_n))$ implies $x_{n-1} =y_{n-1}$ for all $n$.

\smallskip

We then obtain a description of the ring $\cA_{\Z,q}$ as a crossed product as follows.

\begin{thm}\label{AinftyZq}
Let $\cA_{\Z,q}$ be, as in Definition \ref{AZqinfty}, the ring
generated by the elements 
of $\widehat{\Z[q]}$ and the elements $\mu_n$, $\mu_n^*$ subject to
the relations \eqref{relsmuBC} and \eqref{relqsigmamu}.
Then $\cA_{\Z,q}$ is the crossed product ring
\begin{equation}\label{hatZqN}
\cA_{\Z,q}=\widehat{\Z[q]}_\infty \rtimes \Q_+^*
\end{equation}
of the action of the group $\Q_+^*$ on the ring $\widehat{\Z[q]}_\infty =\varinjlim_n (\sigma_n:\widehat{\Z[q]} \to   \widehat{\Z[q]})$.
\end{thm}

\proof We have seen in Lemma \ref{ZqLem2} above that the endomorphisms
$\sigma_n$ become invertible in the direct limit $A_\infty=\widehat{\Z[q]}_\infty$ with
inverses given by the $\rho_n(a)=\mu_n a \mu_n^*$.  Thus the crossed product
of the action of the multiplicative semigroup $\N$ on $A_\infty$ extends to a group 
crossed product by the multiplicative group $\Q_+^*$. In other words, the
action of $r=m/n \in \Q_+^*$, with $(m,n)=1$, on $A_\infty$ that gives the crossed
product $A_\infty \rtimes \Q^*_+$ is given by $\rho_r(a)=\mu_m \mu_n^* a \mu_n \mu_m^*$.

To see that there is an isomorphism between $A_\infty \rtimes \Q_+^*$ and $\cA_{\Z,q}$
we first show that they are isomorphic as $\Z$-modules and then that the ring structures 
also agree. As a $\Z$-module $\cA_{\Z,q}$ is spanned by elements of the form $\mu_N f \mu_M^*$,
for $f\in \widehat{\Z[q]}$ and $N,M\in \N$. In fact, using the relations 
\eqref{relsmuBC} and $\sigma_n(f) =\mu_n^* f \mu_n$, which follows from \eqref{relqsigmamu},
one can see that products of monomials of the form $\mu_N f \mu_M^*$ is still of the same form. 
We then show that elements of $A_\infty \rtimes \Q_+^*$ can also be always written in this form. In fact,
such elements are in the span of elements $a \mu_n$ for $a\in A_\infty$ and their adjoints
$\mu_n^* a$. Any element of the form $\mu_N f \mu_N^* \mu_n$ can be written as
$$ \mu_N f \mu_{n/(n,N)} \mu_{N/(n,N)}^* = \mu_{Nn/(n,N)} \sigma_{n/(N,n)}(f) \mu_{N/(n,N)}^* .$$
Conversely, any element of the form $\mu_N f \mu_M^*$ with $(N,M)=1$ can be written
equivalently as $a_N \mu_N \cdot \mu_M^*$, with $a_N=\mu_N f \mu_N^*$, which is
the product of two elements $a_N \mu_N$ and $\mu_M^*=(\mu_M)^*$ in $A_\infty \rtimes \Q_+^*$.
The product in $A_\infty \rtimes \Q_+^*$ is defined by $a \mu_n b \mu_m = a \rho_n(b) \mu_{nm}$,
which, by the relations $\rho_n(b)=\mu_n b \mu_n^*$ and $\mu_n^* \mu_n=1$, is the 
product of $a \mu_n$ and $b\mu_m$ in $\cA_{\Z,q}$.  
\endproof

This again should be compared with the different situation of the BC algebra, where 
one has a semigroup crossed product $\Q[\Q/\Z] \rtimes \N$. We can give a concrete 
realization of the ring $\widehat{\Z[q]}_\infty$ with the action of $\Q^*_+$ in
terms of functions of rational powers of $q$. 

\begin{prop}\label{polyfrac}
Let $\cP_\Z=\Z[q^r; r\in \Q_+^*]$ denote the ring of polynomials in rational powers of the 
variable $q$. Let $\hat\cP_\Z$ denote the completion
\begin{equation}\label{fracpolhat}
\hat\cP_\Z = \varprojlim_N \cP_\Z / \fI_N ,
\end{equation}
where $\fI_N$ is the ideal generated by the $(q^r)_N =(1-q^r)\cdots (1-q^{rN})$ for
$r\in \Q_+^*$. Then sending an element $a=\rho_n(f)=\mu_n f \mu_n^*$, for 
$f\in \widehat{\Z[q]}$, to the function $f(q^{1/n})$ gives an isomorphism 
\begin{equation}\label{isofracpol}
\Phi: \widehat{\Z[q]}_\infty \stackrel{\simeq}{\to}  \hat\cP_\Z .
\end{equation}
The action of $\Q_+^*$ on $\widehat{\Z[q]}_\infty$ corresponds, under this identification,
to the action on $\hat\cP_\Z$ given by $\rho_r(f)(q)=f(q^r)$. 
\end{prop}

\proof First notice that the map that sends $\mu_n f \mu_n^*$ to $f(q^{1/n})$
is well defined as a ring homomorphism $\widehat{\Z[q]}_\infty \to \hat\cP_\Z$.
In fact, if $f_1-f_2 \in ((q)_N)$ we see that $f_1(q^{1/n})-f_2(q^{1/n})$ is in the 
ideal $\fI_N$. In other words, if we write elements in the Habiro ring as series
$f(q)=\sum_m P_m(q)$, with $P_m(q)$ in $((q)_m)$, the map $\Phi$
gives $\Phi(\rho_n(f))(q)= \sum_m P_m(q^{1/n})$ where $P_m(q^{1/n})$
is in $((q^{1/n})_m)\subset \fI_n$. 

By construction, all functions in $\hat\cP_\Z$
can be described as combinations of functions obtained from elements 
in the Habiro ring by changes of variable $q \mapsto q^r$ with $r\in \Q_+^*$.
Thus, the morphism $\Phi$ is surjective. To see that is is also injective, notice
that if $\Phi(\mu_n f \mu_n^*) \in \fI_N$ for some $N$, then there exists a
sufficiently large $M\in \N$, such that $n|M$, and such that $f(q^{M/n}) \in ((q)_m)$
for some $m\in \N$. This implies that, for $k=M/n$, we have $\sigma_k(f) \in ((q)_m)$,
but the injectivity of the $\sigma_k$ implies that $f=0$ in the Habiro ring $\widehat{\Z[q]}$.
Thus, $\Phi$ is an isomorphism. It is then clear by construction that the action of $\Q^*_+$
on $\widehat{\Z[q]}_\infty$ takes the form $\rho_r(f)(q)=f(q^r)$ on $\hat\cP_\Z$.
\endproof

\medskip

\subsubsection{Taylor expansions at roots of unity}

We first recall how the Taylor expansion at roots of unity for Habiro
functions is 
defined. For any choice of a root of unity $\zeta\in \cZ$, one has an
embedding  
of rings $\Z[q] \hookrightarrow \Z[\zeta,q]$. Moreover, suppose given 
$f_N \in \Z[q]/((q)_N)$. Since for $N > i \, {\rm ord}(\zeta)$ we have
$(q-\zeta)^i | (q)_N$, 
we have a corresponding ring homomorphism 
\begin{equation}\label{evTi}
\ft_\zeta^{(i)}: \Z[q]/((q)_N) \to \Z[\zeta,q]/((q-\zeta)^i), \ \ \
f\mapsto \ft_{\zeta}^{(i)}(f) 
\end{equation}
with $\ft_\zeta^{(0)}=ev_\zeta$. These maps are compatible with the
projective system,  
hence they define, for an element $f\in \widehat{\Z[q]}=
\varprojlim_N \Z[q]/((q)_N)$ a Taylor expansion at $\zeta \in \cZ$,
given by the resulting 
ring homomorphism
\begin{equation}\label{Taylor}
\ft_\zeta: \widehat{\Z[q]} \to \Z[\zeta][[q-\zeta]].
\end{equation}
It was shown in \cite{Hab1} that the homomorphism \eqref{Taylor} is
injective. The following result is already contained in \cite{Hab1},
\cite{Hab2}, but we reproduce it here explicitly for later use.

\begin{lem}\label{TaylorPolyn}
For $f\in \Z[q]/((q)_N)$ the class of a polynomial $P\in \Z[q]$, and
for a choice of 
$i$ and $\zeta$ such that $i {\rm ord}(\zeta)< N$, we have
\begin{equation}\label{tauiP}
\ft_{\zeta}^{(i)}(f)(q)=\sum_{k=0}^i  \frac{1}{k!} P^{(k)}(\zeta) \,
(q-\zeta)^k. 
\end{equation}
\end{lem}

\proof For $P(q)=a_0+a_1 q + a_2 q^2 +\cdots + a_\ell q^\ell$ we write 
$$ P(q)= a_0 + a_1 (\zeta + (q-\zeta)) + \cdots + a_\ell (\zeta +
(q-\zeta))^\ell = $$ 
$$ \sum_{j=0}^\ell a_j \sum_{k=0}^j \binom{j}{k} \zeta^{j-k}
(q-\zeta)^k  = $$ 
$$ \begin{array}{l} 
 a_0 + a_1 \zeta + \cdots + a_\ell \zeta^\ell + \\[2mm] 
 (a_1 + 2 a_2 \zeta+\cdots + \ell a_\ell \zeta^{\ell-1}) (q-\zeta) + \\[2mm]
 (a_2 + \binom{3}{2} a_3 \zeta +\cdots
+ \binom{\ell}{2} \zeta^{\ell-2}) (q-\zeta)^2 + \\[2mm]
 \cdots + a_\ell (q-\zeta)^\ell = \end{array} $$
$$  P(\zeta) + P'(\zeta) (q-\zeta) + \frac{1}{2} P''(\zeta) (q-\zeta)^2 +\cdots
+ \frac{1}{\ell !} P^{(\ell)}(\zeta) (q-\zeta)^\ell , $$
so that, for $k<i$, the coefficient of $(q-\zeta)^k$ in
$\ft_\zeta^{(i)}(f)$ is 
given by the expression
$$ \frac{1}{k!} P^{(k)}(\zeta) = \sum_{j=k}^\ell a_j \binom{j}{k}
\zeta^{j-k}. $$  
This gives \eqref{tauiP} by taking the above expression for $P(q)$ modulo
$(q-\zeta)^i$.
\endproof

The endomorphisms $\sigma_n$ of the Habiro ring induce endomorphisms
of the rings $\Z[\zeta][[q-\zeta]]$ via the Taylor expansion at roots
of unity. 

\begin{lem}\label{sigmanTaylor}
The endomorphisms $\sigma_n :\widehat{\Z[q]}\to \widehat{\Z[q]}$ of
the Habiro ring 
given by $\sigma_n(f)(q)=f(q^n)$ induce an action of the
multiplicative semigroup $\N$ 
on the Taylor expansions at roots of unity $\sigma_n
(\ft_\zeta(f))=\ft_\zeta(\sigma_n(f))$, 
where the Taylor expansion $\ft_\zeta(\sigma_n(f))$ is given by
\begin{equation}\label{Tsigman}
\ft_\zeta(\sigma_n(f)) = \sum_k \frac{1}{k!} (P\circ \sigma_n)^{(k)}
(\zeta) \, (q-\zeta)^k, 
\end{equation}
for $f$ represented by a polynomial $P(q)=a_0+a_1 q + a_2 q^2 +\cdots
+ a_\ell q^\ell$ and 
$(P\circ \sigma_n)(q)= a_0 + a_1 q^n +\cdots +a_\ell q^{n\ell}$.
\end{lem}

\proof The Taylor expansion of Habiro functions at a given root of
unity $\zeta \in \cZ$ 
is an injective ring homomorphism $\ft_\zeta: \widehat{\Z[q]} \to
\Z[\zeta][[q-\zeta]]$. 
Thus, the map $\ft_\zeta(f) \mapsto \ft_\zeta(\sigma_n(f))$ is a
well defined ring 
homomorphism from the image of $\widehat{\Z[q]}$ in
$\Z[\zeta][[q-\zeta]]$ to itself. 
One then checks directly that, given $f\in \widehat{\Z[q]}$
represented by a polynomial 
$P(q)=a_0+a_1 q + a_2 q^2 +\cdots + a_\ell q^\ell$, the Taylor
expansion of $\sigma_n(f)$ 
is obtained by considering the polynomial
$$ P(q^n)= a_0 + a_1 (\zeta+(q-\zeta))^n + \cdots + a_\ell (\zeta +
(q-\zeta))^{n\ell}. $$ 
The coefficient of $(q-\zeta)^k$ in this expansion is of the form
$$ \sum_j a_j \binom{jn}{k} \zeta^{jn-k} =  \frac{1}{k!} (P\circ
\sigma_n)^{(k)} (\zeta), $$ 
So that we get 
$$ \ft_\zeta (\sigma_n(f)) = \sum_{k\geq 0} \frac{1}{k!} (P\circ
\sigma_n)^{(k)}  
(\zeta) (q-\zeta)^k. $$
\endproof

Notice, however, the following simple remark on the 
endomorphisms $\sigma_n$ and the Taylor expansion.

\begin{cor}\label{notzetan} 
The evaluations $ev_\zeta$ of Habiro functions satisfy the property
\begin{equation}\label{evzetan}
ev_\zeta(\sigma_n(f)) = ev_{\zeta^n} (f).
\end{equation}
However, on the full Taylor expansion this no longer holds and one finds
\begin{equation}\label{tauzetan}
\ft_\zeta (\sigma_n(f)) \neq \ft_{\zeta^n}(f).
\end{equation}
\end{cor}

\proof The case of the evaluations $ev_\zeta$ follows directly from Proposition 
\ref{semiact}. On the full Taylor expansion we have
$$ \ft_\zeta(\sigma_n(f)) = \sum_{k\geq 0} \frac{1}{k!} (P\circ \sigma_n)^{(k)} 
(\zeta) (q-\zeta)^k, $$
in the notation of Lemma \ref{sigmanTaylor}, while we have
$$ \ft_{\zeta^n}(f) = \sum_{k\geq 0} \frac{1}{k!} P^{(k)}(\zeta^n) (q-\zeta^n)^k, $$
where
$$ \frac{1}{k!} P^{(k)}(\zeta^n)= \sum_j a_j \binom{j}{k} \zeta^{n(j-k)} . $$
\endproof

Notice, however, that while we do not have the simple equality \eqref{tauzetan}, 
one can still compare the Taylor expansion of $f$ at $\zeta^n$ and the one of
$\sigma_n(f)$ at $\zeta$ in the following way.  If $J_\zeta$ denotes the ideal
$\Ker(ev_\zeta)$, then, as one sees from Lemma \ref{TaylorPolyn}, the Taylor
expansion of $f$ at $\zeta$ comes from the filtration in powers of $J_\zeta$.
The identity \eqref{evzetan} gives $\sigma_n^{-1}( J_\zeta )= J_{\zeta^n}$ and
the same holds for powers of these ideals. This means that there are operators
$\cT_{\zeta,n}$, which can be put in the form of triangular matrices, that
relates the power series $\ft_\zeta(\sigma_n(f))$ and $\ft_{\zeta^n}(f)$.

\medskip

\subsubsection{Compatibility of endomorphism actions} 

To show in what sense the action on the Habiro ring of the endomorphisms $\sigma_n$ is
``compatible" with the action \eqref{sigmanBC} of the BC algebra, we show how to obtain
the $C^*$-algebra of operators  \eqref{erOp} on the Hilbert space $\ell^2(\N)$ using Habiro
functions together with all their evaluations at roots of unity. Notice that this construction
does not give a representation of the Habiro ring on $\ell^2(\N)$, because we are going to
associate to a Habiro function all of its evaluations at roots of unity. While through
each evaluation at a given $\zeta\in \cZ$ one obtains a surjective homomorphism of 
the Habiro ring onto a piece of the abelian part of the BC algebra factoring through 
the ring $\Z[\zeta]$, these do not assemble to give a representation of the
product $\prod_\zeta \Z[\zeta]$. However, the construction given here, which is based
on these partial representations, suffices to explain in what sense one has a 
compatibility, through evaluations at roots of unity, between the actions of $\sigma_n$ 
on the Habiro ring and on the abelian part of the BC algebra.

\medskip

\begin{lem}\label{ZzetaRep}
For a given $\zeta\in \cZ$, the map $P\mapsto E_P$ with
\begin{equation}\label{Paction}
E_P \epsilon_n = P(\zeta^n) \epsilon_n 
\end{equation}
defines a representation of the ring $\Z[\zeta]$ as bounded operators
on the Hilbert space 
$\ell^2(\N)$. This defines a norm $\| f \|_{\cZ}$ for elements
$f\in \widehat{\Z[q]}$.
\end{lem}

\proof For $\zeta\in \cZ$ and $P\in \Z[\zeta]$, $P(\zeta)=a_0 + a_1
\zeta + \cdots + a_r \zeta^r$, 
we have $| P(\zeta^k) |\leq \sum_i |a_i|$, for all $k\in \N$. Thus, we have
$$ \| P \|_{\cB(\ell^2(\N))} = \sup_{v\neq 0} \frac{\| P v \|}{\| v
\|} \leq \sum_i |a_i| < \infty . $$ 
\endproof

Thus, we can consider the abelian $C^*$-subalgebra of $\cB(\ell^2(\N))$ 
generated by all the operators $E_{ev_\zeta(f)}$ as in \eqref{Paction}, for 
$f\in \widehat{\Z[q]}$ and $\zeta \in \cZ$.

\begin{defn}\label{abAlgZq}
Let $C^*_{\cZ}(\widehat{\Z[q]})$ denote the abelian $C^*$-algebra generated by all the operators in $\cB(\ell^2(\N))$ of the form
\begin{equation}\label{Ezetaf}
E_{\zeta,f} \, \epsilon_n = ev_{\zeta^n}(f)\, \epsilon_n,
\end{equation}
for $\zeta\in \cZ$ and $f\in \widehat{\Z[q]}$.
\end{defn}

We have the following result relating the algebra constructed 
above to the abelian part of the Bost--Connes algebra.

\begin{prop}\label{BCandZqev}
Any element $\rho\in  \hat\Z^*$ determines an isomorphism between the
$C^*$-algebra $C^*_{\cZ}(\widehat{\Z[q]})$ and the abelian part $C^*(\Q/\Z)$ 
of the Bost--Connes algebra.
\end{prop}

\proof The choice of an element $\rho\in \hat\Z^*$ can be viewed as the choice of 
an embedding of the roots of unity in $\C$. In particular, any $\rho\in \hat\Z^*$
determines a representation of the Bost--Connes algebra on $\ell^2(\N)$, where
the action of the abelian part $C^*(\Q/\Z)$ is given by the action of the generators
$e(r)$, for $r\in \Q/\Z$ as
\begin{equation}\label{rhoer}
\pi_\rho(e(r)) \epsilon_n = \zeta_r^n \, \epsilon_n,
\end{equation}
where $\zeta_r = \rho(r)$ is the root of unity corresponding to the generator $e(r)$.
Moreover, given two such elements $\rho,\rho'\in \hat\Z^*$, one has
$$ \pi_\rho(C^*(\Q/\Z))=\pi_{\rho'}(C^*(\Q/\Z)) \subset \cB(\ell^2(\N)), $$
since, for a given generator $e(r)$ with $r\in \Q/\Z$, one has
$$ \pi_{\rho'}(e(r))= \pi_\rho (e(\gamma(r))), $$
where the automorphism $\gamma$ of $\Q/\Z$ is given by $\gamma =\rho^{-1} \rho' \in \hat\Z^*$.
We use here the identification $\hat\Z =\Hom(\Q/\Z,\Q/\Z)$. 

We see that, for any given choice of $\rho\in \hat\Z^*$, we have an inclusion
of the algebra $\pi_\rho(C^*(\Q/\Z))\subset \cB(\ell^2(\N))$
in $C^*_{\cZ}(\widehat{\Z[q]}) \subset \cB(\ell^2(\N))$. In fact,
it suffices to show that, for all the generators $e(r)$ of $C^*(\Q/\Z)$, with $r\in \Q/\Z$
the operator $\pi_\rho(e(r))$ is in the algebra $C^*_{\cZ}(\widehat{\Z[q]})$. This follows
by writing $\pi_\rho(e(r))$ equivalently as the operator $ev_{\rho(r)}(f)$, for $f(q)=q$.
In fact, this acts like
$$ E_{\rho(r),f} \epsilon_n = \zeta_r^n \, \epsilon_n, $$ 
with $\zeta_r=\rho(r) \in \cZ$.

To see that the reverse inclusion also holds, consider an element $ev_\zeta(f) \in \Z[\zeta]$,
acting on $\ell^2(\N)$ as in \eqref{Ezetaf}. Suppose given an element $\rho\in \hat\Z^*$. Then
there exists an $r\in \Q/\Z$ such that $\zeta =\rho(r)$. If $ev_\zeta(f)=
P(\zeta)=a_0+a_1 \zeta +\cdots+ a_n \zeta^n$, we can then write the operator $E_{\zeta,f}$ as
$$ E_{\zeta,f} = a_0 + a_1 \pi_\rho(e(r)) + \cdots + a_n \pi_\rho (e(nr)), $$
so that we obtain $C^*_{\cZ}(\widehat{\Z[q]}) \subset \pi_\rho(C^*(\Q/\Z))$.

We have shown in this way that a choice of $\rho\in \hat\Z^*$ determines an
isomorphism $$\iota_\rho : C^*(\Q/\Z) \to C^*_{\cZ}(\widehat{\Z[q]})$$ given on generators by 
\begin{equation}\label{iotarho}
\iota_\rho(e(r))= E_{\rho(r),f}, \ \ \ \text{ for } \ \  f(q)=q.
\end{equation}
\endproof

Notice that the $C^*$-algebra $C^*_{\cZ}(\widehat{\Z[q]})$ is not
a representation of the Habiro ring $\widehat{\Z[q]}$, hence in
particular \eqref{iotarho} is {\em not} a ring homomorphism between the
abelian part of the BC algebra and the Habiro ring, though it
gives a convenient parameterization of the operators 
in terms of Habiro functions and roots of unity. 

The action of $\hat\Z^*$ by automorphisms of of $C^*(\Q/\Z)$ is then given  by
\begin{equation}\label{actZq}
\gamma (E_{\zeta,f})=E_{\theta_\gamma(\zeta),f},
\end{equation}
where $\theta: \hat\Z^* \to \Gal(\Q^{ab}/\Q)$ is the class field theory isomorphism.

The ``compatibility" between the action of the endomorphisms $\sigma_n$ of
the Habiro ring and of the abelian part of the BC algebra can then be
expressed as the property that
\begin{equation}\label{compatibility}
 \sigma_n (E_{\zeta,f}) =E_{\zeta^n,f} = E_{\zeta,f^n} =
E_{\zeta,\sigma_n(f)}. 
\end{equation}

\medskip

\subsection{Hopf algebra and endomotive}

As we have seen above, we do not have a direct relation between the Habiro
ring and the crossed product ring $\cA_{\Z,q}$ and the BC algebra in terms
of Hilbert space representations. The main difficulty in relating them directly
is the fact that the maps $\sigma_n$ have very different properties in the two
cases. There is, however, some further structure of the Habiro ring, which 
makes it possible to obtain from the endomorphisms $\sigma_n$ a construction  
similar to the endomotive description of the BC system.

\smallskip

First recall that, as introduced in \cite{CCM}, endomotives are a category
of noncommutative spaces constructed out of projective limits of Artin
motives with semigroup actions, with morphisms that are also constructed 
out of correspondences between Artin motives, with a suitable compatibility
condition with respect to the semigroup actions. We refer the reader to
\cite{CCM} and also to \S 4 of \cite{CMbook} for details. In particular,
the BC endomotive is obtained by considering the projective system
of Artin motives given by the preimages $s_k^{-1}(1)$ under the
system of endomorphisms $s_k : \bG_m \to \bG_m$ of the multiplicative
group that correspond to the endomorphisms of the ring $\Z[t,t^{-1}]$ given
by $s_k(P)(t)=P(t^k)$. 

\smallskip

It is known that the Habiro ring has the structure of a Hopf algebra. The
argument reproduced here below is due to Habiro: I learned of it
from Morava \cite{Mor}. 

\begin{lem}\label{HabiroHopf}
The coproduct $\Delta(q)=q\otimes q$, counit $\varepsilon(q)=1$, and 
antipode $S(q)=q^{-1}$ on the ring $\Z[q,q^{-1}]$ induce a coproduct, counit, 
and antipode on $\widehat{\Z[q]}$ with respect to the completed tensor product
\begin{equation}\label{complotimes}
\widehat{\Z[q]}\hat\otimes \widehat{\Z[q]} =\varprojlim_n (\Z[q]\otimes\Z[q])/
((q)_n \otimes 1, 1\otimes (q)_n).
\end{equation}
This makes the Habiro ring into a Hopf algebra. 
\end{lem}

\proof The compatibility of the coproduct with the quotients defining the Habiro ring
follows from the fact that
$$ \Delta (1-q^N)= (1-q^N) \otimes 1 + q^N \otimes (1-q^N) $$
so that 
$$ \Delta((1-q^N)^{2M})\in ((1-q^N)^M)\otimes (1) + (1) \otimes ((1-q^N)^M). $$
One then uses an equivalent description of  the Habiro ring as projective limit
$$ \widehat{\Z[q]} = \varprojlim_{N,M} \Z[q]/((1-q^{N!})^M), $$
to see that the coproduct is well defined on the projective limit, with the
corresponding completion, obtained using the ideals 
$((1-q^{N!})^M\otimes 1, 1\otimes (1-q^{N!})^M)$. For the antipode the
argument is similar, since $(1-q^{-N!})^M \in (1-q^{N!})^M\,\, \widehat{\Z[q]}$.
\endproof

The endomorphisms $\sigma_n$ of the Habiro ring are induced by the endomorphisms
$s_k(P)(q)=P(q^k)$ of the polynomial ring $\Z[q,q^{-1}]$, which in turn define endomorphisms
of the multiplicative group $\bG_m$. They in turn determine endomorphisms of the
affine group scheme $\hat\bG$ dual to the Hopf algebra $\widehat{\Z[q]}$, with
$\hat\bG(R)=\Hom(\widehat{\Z[q]},R)$, for a given commutative ring $R$. 
Let $\fs_k: \hat\bG\to \hat\bG$ be the morphism induced dually by $\sigma_n$.
We then let $\bX_k = \fs_k^{-1}(1)$ with the projective system of maps
$\bX_n \to \bX_m$ given by $x\mapsto \fs_k(x)$ for $n=mk$. Then the
projective limit $\bX = \varprojlim_n \bX_n$ with the induced action gives
an endomotive associated to the Habiro ring and its dual group $\hat\bG$, 
which parallels the construction of the BC endomotive from the multiplicative 
group $\bG_m$.

\medskip

\subsection{Operators from Taylor expansions}\label{TaySec}

The content of this subsection is a side remark with respect to the main
theme of the paper and will not be needed in what follows, so it
can be skipped by the reader interested in moving on directly to the 
next section.

\smallskip

We construct here Hilbert space operators associated to the data of
Taylor expansions $\ft_\zeta(f)$ of the Habiro functions at all $\zeta\in \cZ$.
These will in general give unbounded operators, except in the case where
only finitely many terms of the expansion are nonzero. We consider a 
time evolution that generalizes the BC dynamics and recovers the extremal 
KMS states of the BC system at low temperature from the values at roots of 
unity, as well as the higher terms in the Taylor expansion from suitable 
limits of positive temperature Gibbs states.

\medskip

Let $\cH=\ell^2(\N)\otimes \ell^2(\N\cup \{0\})$ with the canonical basis 
$\epsilon_{n,m}=\epsilon_n\otimes \epsilon_m$. 
We set
\begin{equation}\label{tauzetafOps}
T^{(i)}_{\zeta,f} \, \epsilon_{n,m} := \sum_{k=0}^i   \frac{1}{k!} (P\circ \sigma_n)^{(k)}(\zeta) \, 
\epsilon_{n,m+k}.
\end{equation}
We also introduce shift operators $\delta_k$ defined as
\begin{equation}\label{deltak}
 \delta_k \epsilon_{n,m}= \epsilon_{n,m+k}
\end{equation}
so that we write
\begin{equation}\label{zetataudelta}
T_{\zeta,f}^{(i)} =  \sum_{k=0}^i   \frac{1}{k!} ev_\zeta((f\circ \sigma_n)^{(k)}) \delta_k,
\end{equation}
where the notation $ev_\zeta((f\circ \sigma_n)^{(k)}) $ is understood in the sense of
\eqref{tauzetafOps}, according to the definition of the Taylor expansions at roots of
unity of the Habiro functions discussed above.

\smallskip

\begin{lem}\label{taufOps}
The operators \eqref{tauzetafOps} are bounded operators on $\cH$ for all $\zeta\in \cZ$
for all $f \in \Z[q]/((q)_N)$ with $N> i {\rm ord}(\zeta)$. These operators determine
operators $T_{\zeta, f}$ on $\cH$ for all $f\in \widehat{\Z[q]}$, satisfying
$T_{\zeta, f h}=T_{\zeta, f} T_{\zeta, h}$, for each $\zeta\in \cZ$ and 
$f,h\in \widehat{\Z[q]}$. The operators $T_{\zeta,f}$ are typically unbounded, except
in the case where only finitely many terms of the series are nonzero.
\end{lem}

\proof Only finitely many terms are involved in the sum \eqref{zetataudelta}, hence
they are bounded operators. Consider the the operators of the form
\begin{equation}\label{Tzetaf}
T_{\zeta,f}\, \epsilon_{n,m} =\sum_{k\geq 0} \frac{1}{k!} ev_\zeta((f\circ\sigma_n)^{(k)}) 
\epsilon_{n,m+k}  .
\end{equation}
To see the multiplicative property for fixed $\zeta\in \cZ$, notice that 
$\ft_\zeta(fh)=\ft_\zeta(f)\ft_\zeta(h)$,
since $\ft_\zeta: \widehat{\Z[q]} \to \Z[\zeta][[q-\zeta]]$ is a ring homomorphism.
The property then follows using 
$\ft_\zeta(f) = \sum_{k\geq 0} \frac{1}{k!} f^{(k)}(\zeta) (q-\zeta)^k$, together 
with the fact that the homomorphisms $\sigma_n$ of $\widehat{\Z[q]}$ induce 
homomorphisms of the Taylor series expansions as in Lemma \ref{sigmanTaylor}.
To see that the operators \eqref{Tzetaf} are typically unbounded, consider elements
in the Habiro ring defined by series $f(q)=\sum_m P_m(q)$ with $P_m \in ((q)_m)$
of the form $P_m(q) = n_m (q)_m$, for an arbitrary sequence $(n_m)$ of non-negative 
integers. Taking $\zeta=1$ one sees that the Taylor expansion of $f$ has coefficients
that are all positive integers, since $(q)_m$ is a polynomial with positive coefficients
in the variable $t=q-1$. Thus, when acting with $T_{1,f}$ on
a basis vector $\epsilon_{1,m}$ one finds
$T_{1,f}\epsilon_{1,m} =  \sum_{k\geq 0} \frac{1}{k!} ev_1(f^{(k)}) \epsilon_{1,m+k}$,
with $\frac{1}{k!} ev_1(f^{(k)}) \geq n_k$.
Thus, the element $T_{1,f}\epsilon_{1,m}$ is in $\cH$ only if all but finitely  many of the
$n_k$ are zero. More generally, for any element $f=\sum_m P_m(q)\in \widehat{\Z[q]}$,
the coefficients $\frac{1}{k!} ev_1(f^{(k)})$ are integers, hence $T_{1,f}$ is bounded
only if all but finitely many of these terms are zero, which implies that $f$ is in fact
a polynomial, since the Taylor expansion at $\zeta=1$ is an injective homomorphism.
\endproof

We just write $\Z[q]\subset \widehat{\Z[q]}$ for the functions 
$f=\sum_m P_m(q)\in \widehat{\Z[q]}$ that are in fact polynomial, as
above. We let $C^*_{\ft,\cZ}(\Z[q])$ denote the $C^*$-subalgebra of the
algebra of bounded operators on $\cH=\cB(\ell^2(\N)\otimes \ell^2(\Z_{\geq 0}))$ 
generated by the operators $T_{\zeta,f}$, for $f\in \Z[q]\subset \widehat{\Z[q]}$ and $\zeta\in \cZ$. 
Moreover, we denote as above with $\mu_n$ the
isometries $\mu_n \epsilon_{k,m}= \epsilon_{nk,m}$.
The $\mu_n$ satisfy the relations \eqref{relsmuBC}
and the commutation relations
\begin{equation}\label{zerocomm}
[\mu_n,\delta_k]=[\mu_n,\delta_k^*]=[\mu_n^*,\delta_k]=[\mu_n^*,\delta_k^*]=0, \ \ 
\forall, n \in \N, \forall k\in \Z_{\geq 0},
\end{equation}
with the shift operators $\delta_k$ of \eqref{deltak}. 
We denote by $\cA_{\ft,\cZ}$ the $C^*$-subalgebra of $\cB(\cH)$ generated by the 
$T_{\zeta,f}$ as above, together with the operators $\mu_n$, $\mu_n^*$.

\begin{prop}\label{repAZqtaui}
The $C^*$-algebra $\cA_{\ft,\cZ}$ is a semigroup crossed product
\begin{equation}\label{tausemicross}
\cA_{\ft,\cZ} = C^*_{\ft,\cZ}(\Z[q])\rtimes_{\rho} \N,
\end{equation}
where $\rho_n(T_{\zeta,f})=\mu_n T_{\zeta,f} \mu_n^*$.
The endomorphisms $\rho_n$ satisfy
\begin{equation}\label{rhohbaract}
\rho_n (T_{\zeta,f})  =e_n  T_{\zeta , \eta_n(f)} e_n = 
\frac{1}{n} \sum_{\xi^n =\zeta} ev_\xi(u) T_{\zeta,\eta_n(f)} ,
\end{equation}
with $e_n =\mu_n\mu_n^*$ and with $u\in \widehat{\Z[q]}$ the class of $u(q)=q$. 
\end{prop}

\proof The semigroup crossed product of $C^*_{\ft,\cZ}(\Z[q])\rtimes_{\rho} \N$ is 
implemented by the action
$\rho_n (T_{\zeta,f}) \epsilon_{m,\ell} = \mu_n T_{\zeta,f} \mu_n^* \epsilon_{m,\ell} 
 = e_n \sum_{k\geq 0} \frac{1}{k!} (f\circ \sigma_{m/n})^{(k)}(\zeta) \, e_n \,
\epsilon_{m,\ell+k}$,
where $e_n=\mu_n \mu_n^*$ acts as
$e_n \epsilon_{m,\ell} = 0$ if $n\not| m$ and  as
$e_n \epsilon_{m,\ell} = \epsilon_{m,\ell}$ if $n | m$.
Thus, we obtain \eqref{rhohbaract}, where the right hand side is well defined
because it is nonzero only on the elements $\epsilon_{m,\ell}$ where $n|m$,
where the operator $T_{\zeta,f}$ of $\ft_\zeta(f)$ acts as $\ft_\zeta(h \circ \sigma_{m/n})$ for 
$h=f\circ \sigma_n \in \hat R_n$, so that it makes sense to apply the morphism
$\eta_n$. Notice then that the action of $e_n$ on the given Hilbert space coincides
with the action of the element of $C^*_{\ft,\cZ}(\Z[q])$ given by
$\sum_{\xi^n =\zeta} ev_\xi(u)$.  In fact, we have
$$ \frac{1}{n} \sum_{\xi^n =\zeta} ev_\xi(u) \epsilon_{m,\ell} = 
\frac{1}{n} \sum_{\xi^n =\zeta} \xi^m \epsilon_{m,\ell} = 
\left\{ \begin{array}{ll} 1 & n|m \\ 0 & n\not| m.  \end{array}\right. $$

One then argues as in Theorem \ref{AinftyZq} to show that both
$\cA_{\ft,\cZ}$  and $C^*_{\ft,\cZ}(\Z[q])\rtimes_\rho \N$
are linearly generated by elements of the form $\mu_N T_{\zeta,f} \mu_M^*$, for some
$N,M\in \N$, and for some $f\in \Z[q]\subset \widehat{\Z[q]}$ and some
$\zeta\in \cZ$, and that the product in $\cA_{\ft,\cZ}$ agrees with the product given
by the semigroup action by the relation $\rho_n(T_{\zeta,f})=\mu_n T_{\zeta,f} \mu_n^*$.
\endproof

In the following, we denote by $\cA_{\ft,\cZ,\delta}$ the $C^*$-subalgebra of $\cB(\cH)$ 
generated by $\cA_{\ft,\cZ}$ together with the shift operators $\delta_\ell$ and $\delta_\ell^*$ 
defined by \eqref{deltak}. We also use of the shorthand notation
\begin{equation}\label{shorttzetaf}
\ft_{\zeta,k}(f):=\frac{1}{k!}  ev_\zeta((f\circ\sigma_n)^k).
\end{equation}

\begin{lem}\label{Adelta}
The algebra $\cA_{\ft,\cZ,\delta}$ is generated by the $T_{\zeta,f}$, the $\mu_n$, $\mu_n^*$ 
and the $\delta_\ell$, $\delta_\ell^*$, subject to all the relations previously described
for the algebra $\cA_{\ft,\cZ}$, with the additional relations \eqref{zerocomm} and
$\delta_\ell T_{\zeta,f}=T_{\zeta,f}\delta_\ell$ and
\begin{equation}\label{deltastarTf}
[\delta_\ell^*, T_{\zeta,f}] = Y_{\zeta,f,\ell},
\end{equation}
where, for  $f\in \Z[q]\subset \widehat{\Z[q]}$, the operator
$Y_{\zeta,f,\ell}$ is the bounded operator on $\cH$ acting as
\begin{equation}\label{Yzetafell}
Y_{\zeta,f,\ell} \epsilon_{n,m} = \sum_{k\geq 0} \frac{1}{(k+\ell)!} 
ev_\zeta((f\circ\sigma_n)^{k+\ell}) \epsilon_{n,m+k}.
\end{equation}
\end{lem}

\proof The only statement that does not follow immediately from the definitions is
the last. We check
$$ \delta_\ell^* T_{\zeta,f}\, \epsilon_{n,m} = \left\{\begin{array}{ll} \sum_{k\geq 0} 
\ft_{\zeta,k}(f\circ\sigma_n) \epsilon_{n, k+m-\ell} & k+m-\ell \geq 0 \\
0 & k+m-\ell <0, \end{array}\right. $$
while
$$ T_{\zeta,f} \delta_\ell^*\, \epsilon_{n,m} =\left\{\begin{array}{ll} \sum_{k\geq 0} 
\ft_{\zeta,k}(f\circ\sigma_n) \epsilon_{n, k+m-\ell} & m-\ell \geq 0 \\
0 & m-\ell <0, \end{array}\right.  . $$
This gives
$(\delta_\ell^* T_{\zeta,f}-T_{\zeta,f} \delta_\ell^*)\, \epsilon_{n,m} =\sum_{k\geq \ell-m}
\ft_{\zeta,k}(f\circ\sigma_n) \, \epsilon_{n,m+k-\ell}$. 
\endproof

We now construct a time evolution on the algebra $\cA_{\ft,\cZ,\delta}$ introduced 
above, defined by analogy to the case of the BC algebra.

\begin{prop}\label{BChbar}
Fix a choice of a parameter $\hbar$ with $0< \hbar <1$, such that
the only solution $(r,x)\in \Q^*_+\times \Z$ of the equation
$\log r - x \log \hbar =0$ is given by the pair $r=1$ and $x=0$.
Then setting
\begin{equation}\label{sigmathbar}
\sigma_t(\delta_k) = \hbar^{-ikt} \delta_k , \ \ \ \ \sigma_t(\mu_n)=
n^{it}\, \mu_n, 
\end{equation}
\begin{equation}\label{sigmathbartau}
\sigma_t (T_{\zeta,f}) \epsilon_{n,m}= \sigma_t \left( \sum_{k\geq 0}
\ft_{\zeta,k}(f\circ \sigma_n) 
\delta_k \right) \epsilon_{n,m} = \sum_{k\geq 0} \frac{1}{k!} (P\circ
\sigma_n)^{(k)}(\zeta)  
\, \hbar^{-ikt} \, \delta_k \epsilon_{n,m}.
\end{equation}
defines a time evolution on the algebra $\cA_{\ft,\cZ,\delta}$. 
The operator $H$ on $\cH$ given by
\begin{equation}\label{Hhbar}
H \, \epsilon_{n,m} = \left( \log(n) - m \log(\hbar) \right)\,
\epsilon_{n,m}  
\end{equation}
is a Hamiltonian for the time evolution $\sigma_t$. The partition
function is of the form 
\begin{equation}\label{Zbetahbar}
Z_\hbar(\beta) = \frac{\zeta(\beta)}{1-\hbar^\beta},
\end{equation}
where $\zeta(\beta)$ is the Riemann zeta function.
\end{prop}

\proof To see that \eqref{sigmathbar} and \eqref{sigmathbartau} 
define a time evolution 
$\sigma: \R \to \Aut(\cA_{\ft,\cZ,\delta})$ we need 
to check the compatibility with the relations in the algebra. The
action \eqref{sigmathbartau} 
is compatible with \eqref{sigmathbar} and both are compatible with
\eqref{rhohbaract} 
and \eqref{deltastarTf}. To see that $H$ is a Hamiltonian we check that
$e^{itH} \mu_k e^{-itH} \epsilon_{n,m} = (nk)^{it} \hbar^{-itm}
n^{-it} \hbar^{itm} \epsilon_{nk,m} = 
k^{it} \mu_k \epsilon_{n,m} =\sigma_t (\mu_k) \, \epsilon_{n,m}$ and that
$e^{itH} \delta_k e^{-itH} \epsilon_{n,m} = n^{it} \hbar^{-it(m+k)}
n^{-it} \hbar^{itm} \epsilon_{n,m+k}= 
\hbar^{-itk} \delta_k \epsilon_{n,m} =\sigma_t (\delta_k) \epsilon_{n,m}$.
Similarly, we have
$e^{itH}  T_{\zeta, f} e^{-itH} \epsilon_{n,m} =
\sum_{k\geq 0} \hbar^{-it(m+k)}   \ft_{\zeta,k}(f\circ\sigma_n)
\hbar^{itm} \epsilon_{n,m+k}= 
 \sum_{k\geq 0} \hbar^{-ikt}  \ft_{\zeta,k}(f\circ\sigma_n)
\delta_k \epsilon_{n,m}= 
\sigma_t(T_{\zeta,f}) \epsilon_{n,m}$.
Notice that, since $k\geq 1$ and $0< \hbar < 1$, we have
$\log(k) - m \log(\hbar) \geq 0$ and it is equal to zero
only for the pair $k=1$ and $m=0$. Moreover, for all $a\in \R_+$, 
the equation $\log(k) - m \log(\hbar)=a$ has at most one 
solution with $k\in \N$ and $m\in \Z_{\geq 0}$, else we 
would obtain a nontrivial solution $(r,x)\in \Q^*_+\times \Z$ 
of $\log(r) - x \log(\hbar) =0$.
This implies that the Hamiltonian has positive energy and
all the eigenvalues of $H$ have multiplicity one. 
The partition function is
then given by
$$ Z_\hbar(\beta)= \Tr(e^{-\beta H}) = \sum_{n,m} e^{-\beta
(\log(k)-m\log(\hbar))} = 
(\sum_{n\geq 1} n^{-\beta}) (\sum_{m\geq 0} \hbar^{\beta m}) =
\zeta(\beta) \frac{1}{1-\hbar^\beta}. $$  
\endproof

\smallskip

Since the partition function $Z_\hbar(\beta)$ converges absolutely for $\beta >1$, in
this low temperature range of the thermodynamic parameter we have KMS states 
that are of the Gibbs form
\begin{equation}\label{GibbsKMShbar}
\varphi_\beta(a) = \frac{\Tr(a e^{-\beta H})}{\Tr(e^{-\beta H})} .
\end{equation}
These satisfy the following property.

\begin{prop}\label{KMSinfty}
In the zero temperature limit the Gibbs states \eqref{GibbsKMShbar}
converge weakly to KMS$_\infty$ states satisfying
\begin{equation}\label{KMSinftyhbar}
\varphi_\infty(T_{\zeta,f})=\lim_{\beta\to
\infty}\varphi_\beta(T_{\zeta,f}) 
= ev_\zeta(f).
\end{equation}
Moreover, one obtains limits
\begin{equation}\label{limKMSbeta}
\lim_{\beta\to \infty} \frac{\varphi_\beta(\delta_\ell^*
T_{\zeta,f})}{\hbar^{\beta \ell}} = \ft_{\zeta, \ell}(f),
\end{equation}
that is, the $\ell$-th coefficient of the Taylor expansion 
at $\zeta$ of the Habiro function $f$.
\end{prop}

\proof First notice that, for an element in the algebra of the 
form $\delta_\ell^* T_{\zeta,f}$, we have 
\begin{equation}\label{Gibbstauzeta} 
\varphi_\beta( \delta_\ell^* T_{\zeta,f}) = \hbar^{\beta \ell}
(1-\hbar^\beta) \zeta(\beta)^{-1} 
\sum_n \ft_{\zeta, \ell}(f\circ \sigma_n) n^{-\beta},
\end{equation}
since
$\varphi_\beta( \delta_\ell^* T_{\zeta,f}) = 
Z_\hbar(\beta)^{-1} \, \sum_{n,m}  \langle \epsilon_{n,m},
 \delta_\ell^*  T_{\zeta,f} \epsilon_{n,m} \rangle \,  n^{-\beta}\,
\hbar^{\beta m} $ which gives
 $Z_\hbar(\beta)^{-1} \, \sum_n \ft_{\zeta, \ell}(f\circ \sigma_n)
\,  n^{-\beta}\, \hbar^{\beta \ell} 
 = \hbar^{\beta \ell} (1-\hbar^\beta) \zeta(\beta)^{-1}  
\sum_n \ft_{\zeta, \ell}(f\circ \sigma_n)
 n^{-\beta}$.
In particular, from \eqref{Gibbstauzeta} we have as a special case
\begin{equation}\label{KMSforBC}
\varphi_\beta(T_{\zeta,f}) =(1-\hbar^\beta) 
\zeta(\beta)^{-1} \sum_n ev_{\zeta}(f\circ \sigma_n)\,\, n^{-\beta} .
\end{equation}
One then sees from this that in the weak limit as $\beta \to \infty$
one obtains
$\lim_{\beta\to \infty} \varphi_\beta(\delta_\ell^*
T_{\zeta,f}) =  \lim_{\beta\to \infty}  \hbar^{\beta \ell}
(1-\hbar^\beta) \zeta(\beta)^{-1} 
\sum_n \ft_{\zeta, \ell}(f\circ \sigma_n) n^{-\beta}$.
This limit is zero if $\ell\neq 0$, while for $\ell=0$ it gives
the projection onto the kernel of the Hamiltonian $H$, which 
is the span of the vacuum 
vector $\epsilon_{1,0}$, hence KMS$_\infty$ states are of the form 
$\langle \epsilon_{1,0}, T_{\zeta,f} \epsilon_{1,0} \rangle$, 
which gives \eqref{KMSinftyhbar}.
Similarly, the limits of \eqref{limKMSbeta} give
$\lim_{\beta\to \infty} (1-\hbar^\beta) \zeta(\beta)^{-1}  
\sum_n \ft_{\zeta, \ell}(f\circ \sigma_n)
 n^{-\beta} = \langle \epsilon_{1,0}, 
\delta_\ell^* T_{\zeta,f} \epsilon_{1,0} \rangle = \ft_{\zeta,
\ell}(f)$.
\endproof

We then consider the action of $\hat\Z^*\cong\Gal(\Q^{cycl}/\Q)$ on
the algebra $\cA_{\ft,\cZ}$ given by
$$ \alpha: T_{\zeta, f} \mapsto T_{\alpha(\zeta),f},  $$
with $\alpha(\mu_n)=\mu_n$ and $\alpha(\mu_n^*)=\mu_n^*$.  We recover
the intertwining property of symmetries and Galois action as in the BC
system in the case of KMS$_\infty$ states evaluated at elements of the
rational subalgebra $\cA_{\ft,\cZ,\Q}$, by
\begin{equation}\label{KMSGalTzeta}
\varphi_\infty( \alpha(T_{\zeta,f}) ) = ev_{\alpha(\zeta)} (f) =
\alpha( ev_\zeta(f)) = 
\alpha (\varphi_\infty(T_{\zeta,f}) ).
\end{equation}
As $\langle \epsilon_{1,0}, \delta_\ell^* T_{\zeta,f} \epsilon_{1,0} \rangle 
= \ft_{\zeta, \ell}(f)$ also holds in the case of unbounded
$T_{\zeta,f}$, these last observations apply to non-polynomial 
$f\in \widehat{\Z[q]}$.

\section{Multivariable Bost--Connes systems}\label{MultivarSec}

In \cite{Man}, Manin introduced multivariable versions of the Habiro
ring.  We recall here the construction and we show that these also
admit semigroups of endomorphisms coming from self maps of
higher dimensional algebraic tori. This then serves as a guideline
for a construction we present in this section of multivariable BC 
endomotives, obtained from iterations of self-maps of higher 
dimensional tori.

\subsection{Multivariable Habiro rings and endomorphisms}

The multivariable versions of the Habiro ring are defined in \cite{Man} 
as projective limits, as in the one variable case, by setting
\begin{equation}\label{hatZqi}
\widehat{\Z[q_1,\ldots,q_n]}=\varprojlim_N \Z[q_1,\ldots,q_n]/I_{n,N},
\end{equation}
where $I_{n,N}$ is the ideal 
\begin{equation}\label{InN}
I_{n,N}= \left( (q_1-1)(q_1^2-1)\cdots(q_1^N-1), \ldots,
(q_n-1)(q_n^2-1)\cdots (q_n^N-1) \right). 
\end{equation}
As shown in \cite{Man}, like in the single variable case, the rings \eqref{hatZqi} have 
evaluations at roots of unity. Namely, given a vector $Z=(\zeta_1,\ldots,\zeta_n)$ in $\cZ^n$, 
there is a ring homomorphism
\begin{equation}\label{evzetai}
ev_{Z} : \widehat{\Z[q_1,\ldots,q_n]} \to \Z[\zeta_1,\ldots,\zeta_n].
\end{equation}
Moreover, the functions in the multivariable Habiro rings have a Taylor expansion at
points $Z=(\zeta_1,\ldots,\zeta_n)$ in $\cZ^n$. The Taylor expansion is an {\em injective}
ring homomorphism
\begin{equation}\label{Taylorqizetai}
\ft_{n,Z}: \widehat{\Z[q_1,\ldots,q_n]} \to \Z[\zeta_1,\ldots,\zeta_n][[q_1-\zeta_1,\ldots,q_n-\zeta_n]].
\end{equation}

We introduce the following notation. Given an 
integer $n\times n$-matrices with positive determinant
$\alpha\in M_n(\Z)^+$, and given
$q=(q_1,\ldots,q_n)$, we write $q^\alpha$ for 
\begin{equation}\label{qalpha}
q^\alpha =( q^\alpha_i )_{i=1,\ldots,n}, \ \ \ \text{ with } \ \ \  q^\alpha_i =\prod_j q_j^{\alpha_{ij}}.
\end{equation}
We show now that the semigroup $S_n=M_n(\Z)^+$ acts by endomorphisms 
on the multivariable Habiro ring.

\begin{prop}\label{actSnZqi}
For $\alpha \in M_n(\Z)^+$, the ring homomorphisms $$\sigma_\alpha: \Z[q_1,\ldots,q_n, 
q_1^{-1},\ldots, q_n^{-1}]\to  \Z[q_1,\ldots,q_n,q_1^{-1},\ldots, q_n^{-1}]$$ given by
\begin{equation}\label{sigmaqalpha}
\sigma_\alpha: q\mapsto q^\alpha,
\end{equation}
with $q^\alpha$ defined as in \eqref{qalpha}, induces a ring homomorphism
$$ \sigma_\alpha : \widehat{\Z[q_1,\ldots,q_n]}\to  \widehat{\Z[q_1,\ldots,q_n]} $$
of the multivariable Habiro ring.
\end{prop}

\proof First notice that an equivalent presentation of the Habiro ring can be given as 
\begin{equation}\label{Habiroqinv}
\begin{array}{rl}
\widehat{\Z[q_1,\ldots,q_n]} =& \varprojlim_N \Z[q_1,\ldots,q_n, q_1^{-1}, \ldots, q_n^{-1}]/\cJ_{n,N} 
\\[2mm] = & \varprojlim_N \Z[q_1,\ldots,q_n, q_1^{-1}, \ldots, q_n^{-1}]/\cI_{n,N} ,
\end{array}
\end{equation}
where $\cJ_{n,N}$ is the ideal generated by the polynomials $(q_i-1)\cdots (q_i^N-1)$, for
$i=1,\ldots, n$ and $(q_i^{-1}-1)\cdots (q_i^{-N}-1)$ and $\cI_{n,N}$ is the ideal generated
by the $(q_i-1)\cdots (q_i^N-1)$, for $i=1,\ldots, n$.

For $\alpha\in M_n(\Z)^+$, the map 
\begin{equation}\label{sigmaalpha}
q\mapsto \sigma_\alpha(q)=\sigma_\alpha (q_1,\ldots,q_n) =(q_1^{\alpha_{11}} q_2^{\alpha_{12}}\cdots q_n^{\alpha_{1n}},\ldots,
q_1^{\alpha_{n1}}  q_2^{\alpha_{n2}} \cdots q_n^{\alpha_{nn}})=q^\alpha
\end{equation}
defines a ring homomorphism of $\Z[q_1,\ldots,q_n, q_1^{-1}, \ldots, q_n^{-1}]$. We want to show that it
induces a morphism of the projective limits of the quotients by the ideals $\cI_{n,N}$. To this purpose
it is useful to make another preliminary observation, namely that, as shown in \cite{Hab1}, the ring
$\widehat{\Z[q]}$ can equivalently be described as
\begin{equation}\label{PhiZqhat}
\widehat{\Z[q]} = \varprojlim_{f\in \Phi^*} \Z[q]/(f(q)),
\end{equation}
where $\Phi^*$ is the multiplicative subset of $\Z[q]$ generated by the cyclotomic polynomials
$\Phi_N(q)$. Similarly, we can describe the multivariable versions of the Habiro ring also in
terms of the cyclotomic polynomials instead of the polynomials $(q)_N$ by
\begin{equation}\label{multiPhiZqhat}
\widehat{\Z[q_1,\ldots,q_n]} = \varprojlim_{f_i\in \Phi^*} \Z[q_1,\ldots,q_n]/(f_1(q_1),\ldots,f_n(q_n)).
\end{equation}
The image under $\sigma_\alpha$ of the polynomials $(q_i-1)\cdots (q_i^N-1)$ is given by the polynomials 
$$ P_i^\alpha(q_1,\ldots,q_n,q_1^{-1},\ldots,q_n^{-1})=
(\prod_{j=1}^n q_j^{\alpha_{ij}} -1)\cdots (\prod_{j=1}^n q_j^{N\alpha_{ij}} -1). $$
We show that there are $N_i\in \N$ 
such that there is an inclusion of ideals
\begin{equation}\label{inclideals}
(\Phi_{N_1}(q_1),\ldots,\Phi_{N_n}(q_n)) \supset (P_1^\alpha(q_1,\ldots,q_n,q_1^{-1},\ldots,q_n^{-1}), \ldots, P_n^\alpha(q_1,\ldots,q_n,q_1^{-1},\ldots,q_n^{-1})),
\end{equation}
so that we can view the map $\sigma_\alpha$ as inducing a homomorphism
$$ \sigma_\alpha : \Z[q_1,\ldots,q_n]/\cI_{n,N} \to \Z[q_1,\ldots,q_n]/(\Phi_{N_1}(q_1),\ldots,\Phi_{N_n}(q_n)). $$
Let $\zeta_i$ be a solution of $\Phi_{N_i}(q_i)=0$. Then $\zeta_i^{N_i}=1$ and all these
roots appear with multiplicity one, which is the advantage of using the cyclotomic polynomials
as opposed to the $(q)_N$ here. Substituting in the polynomial $P^\alpha_i$ we obtain
the expressions
$$ \prod_{k=1}^N (\prod_{j=1}^n \zeta_j^{\alpha_{ij}k}-1), $$
for $i=1,\ldots, n$.  Thus, if we choose each $N_j$ in such a way that $N_j | \alpha_{ij} k_i$ for some 
$k_i\in \{1,\ldots, N\}$ and for all $\alpha_{ij}$, we see that the $\zeta_i$
that are solutions of $\Phi_{N_i}(q_i)=0$ are also solutions of 
$P^\alpha_i(q_1,\ldots,q_n,q_1^{-1},\ldots, q_n^{-1})=0$.
\endproof

It is then possible to consider, as in the one-variable case, an associated 
non commutative ring obtained as the (semi)group crossed product
\begin{equation}\label{Habirocrossed}
\cA_{\Z,q,n}=\widehat{\Z[q_1,\ldots,q_n]}_\infty \rtimes \GL_n(\Z)^+,
\end{equation}
where $\cA_{\Z,q,n}$
is the ring generated by the multivariable Habiro functions 
$f\in \widehat{\Z[q_1,\ldots,q_n]}$, together with additional generators 
$\mu_\alpha$ and $\mu_\alpha^*$ subject to relations
\begin{equation}\label{mualpharels}
\begin{array}{ll}
(\mu_\alpha^*)^* =\mu_\alpha & \forall \alpha \in M_n(\Z)^+ \\[2mm] 
\mu_\alpha \mu_\beta = \mu_{\alpha\beta}, & \forall \alpha,\beta \in M_n(\Z)^+ \\[2mm]
\mu_\alpha^* \mu_\beta = \mu_\gamma & \text{for } \, \beta=\alpha\gamma \in M_n(\Z)^+ .
\end{array}
\end{equation}
In particular, the $\mu_\alpha$ are isometries satisfying 
$\mu_\alpha^*\mu_\alpha=1$ for all $\alpha \in M_n(\Z)^+$.
We also have in $\cA_{\Z,q,n}$ the additional relation
\begin{equation}\label{mualphasigmarelHab}
\mu_\alpha^* \, f = \sigma_\alpha(f)\, \mu_\alpha^*, \ \ \ \text{ and } \ \ \  f\, \mu_\alpha =
\mu_\alpha \, \sigma_\alpha(f),
\end{equation}
for all $f\in \widehat{\Z[q_1,\ldots,q_n]}$ and $\alpha \in M_n(\Z)^+$. 

The semigroup crossed product structure \eqref{Habirocrossed} is given as in
Theorem \ref{AinftyZq} for the one-variable case, with the ring $A_{\infty,n}=
\widehat{\Z[q_1,\ldots,q_n]}_\infty$ the inductive limit of the system of
the $\sigma_\alpha : \widehat{\Z[q_1,\ldots,q_n]} \to \widehat{\Z[q_1,\ldots,q_n]}$.
As in the 1-dimensional case, one obtains induced homomorphisms
$\sigma_\alpha: \widehat{\Z[q_1,\ldots,q_n]}_\infty \to \widehat{\Z[q_1,\ldots,q_n]}_\infty$
that are now invertible with inverses implementing the 
group action in  \eqref{Habirocrossed}, given by
\begin{equation}\label{rhoalphaHab}
\rho_\alpha (x) = \mu_\alpha\, x\, \mu_\alpha^*, \ \ \ \text{ for } x \in A_{\infty,n}.
\end{equation}
As in Proposition \ref{polyfrac}, this can also be seen by considering the ring
$\cP_{\Z,n}$ of polynomials in $n$ variables $(q_1^{r_1}, \cdots, q_n^{r_n})$
that are fractional powers of the variables $q_i$, with $r_i\in \Q_+^*$, and
the projective limit 
$$ \hat\cP_{\Z,n} = \varprojlim_N \cP_{\Z,n}/ \fI_{n,N}, $$
with respect to the ideals 
$$ \fI_{n,N}=( (q_1^{r_1}-1)\cdots (q_1^{r_1 N}-1), \ldots, (q_n^{r_n}-1) \cdots (q_n^{r_n N}-1)). $$
for $r_i \in \Q_+^*$, endowed with the action $\rho_\alpha(f)(q^r)=f(q^{\alpha^{-1}r})$ and
$\sigma_\alpha(f)(q^r) =f(q^{\alpha r})$, with $q^r=(q_1^{r_1},\ldots, q_n^{r_n})$, where
for $\alpha\in \GL_n(\Q)^+$ we write 
$$ q^{\alpha r}= (q_1^{\alpha_{11}r_1} q_2^{\alpha_{12}r_2}\cdots q_n^{\alpha_{1n}r_n},\ldots,
q_1^{\alpha_{n1} r_1}  q_2^{\alpha_{n2}r_2} \cdots q_n^{\alpha_{nn}r_n}) .   $$   

\medskip

\subsection{Multivariable BC endomotives}

The construction described above using the multivariable Habiro rings
immediately suggests the existence of an interesting class of ``multivariable"
generalizations of the BC endomotive. Just like the original BC endomotive
is obtained in \cite{CCM} using self-maps of the multiplicative group $\bG_m$,
as recalled in \S \ref{BCSec} above, the multivariable versions are similarly
associated to self-maps of higher rank algebraic tori.

\smallskip

Unlike the case of the action on the Habiro rings, in this case as in the case
of the original BC system, we find semigroup actions of $M_n(\Z)^+$ that
are implemented by isometries hence they do not extend to group actions
of $\GL_n(\Q)^+$.

\smallskip

We let $\bT^n=(\bG_m)^n$ be an $n$-dimensional (split) 
algebraic torus. 
We then consider the multiplicative semigroup of 
self-maps of $\bT^n$ that is the
analog of the semigroup $\N$ acting by \eqref{ttk} 
on the 1-dimensional torus $\bG_m$.

\begin{lem}\label{actMnZplus}
The semigroup $M_n(\Z)^+$ 
acts by endomorphisms of the torus $\bT^n=(\bG_m)^n$. 
The elements of $\SL_n(\Z)\subset
M_n(\Z)^+$ act by automorphisms. 
The point $t_0=(1,1,\ldots,1)\in \bT^n$ 
is a fixed point of all 
the endomorphisms in the semigroup $S_n= M_n(\Z)^+$.
\end{lem}

\proof This can be seen directly using the exponential map
$$ 0\to \Z \to \C \stackrel{\exp(2\pi i \cdot)}{\to} 
\C^* \to 1 $$ 
and the action of $M_n(\Z)^+$ on $\C^n$, preserving 
the lattice $\Z^n$ and the orientation. 
The action of $\alpha \in M_n(\Z)^+$ on the algebra 
$\C[t_i,t_i^{-1}]$ of
the n-torus is given by
\begin{equation}\label{salpha}
 s_\alpha: P(t_i,t_i^{-1}) \mapsto 
P(\exp(2\pi i \sum_j \alpha_{ij} u_j),
\exp(-2\pi i \sum_j \alpha_{ij} u_j)), 
\end{equation}
where $t_i = \exp(2\pi i u_i)$ for $i=1,\ldots,n$. 
Since all the endomorphisms 
are induced by linear maps on $\C^n$, the point $t_0$ 
is a fixed point for all $\alpha\in S_n$.
\endproof

We define as above $t^\alpha$ to be the result of applying
$\alpha\in M_n(\Z)^+$ to the coordinates $t=(t_i)$ of 
$\bT^n$,
\begin{equation}\label{talpha}
t^\alpha := (\exp(2\pi i \sum_j 
\alpha_{ij} u_j))_{j=1,\ldots, n}, \ \ \ 
\text{ for } t_i = \exp(2\pi i u_i),
\end{equation}
that is, equivalently
\begin{equation}\label{talpha2}
t^\alpha =(t^\alpha_i)_{i=1,\ldots,n} \ \ \ 
\text{ with } \ \ \ t^\alpha_i = \prod_j t_j^{\alpha_{ij}}.
\end{equation}
Under multiplication in the semigroup 
$M_n(\Z)^+$ this satisfies
\begin{equation}\label{talphabeta}
t^{\alpha\beta}=(t^\beta)^\alpha.
\end{equation}
For example, in the case of $\alpha\in M_2(\Z)^+$ 
this gives
$$ \alpha=\left(\begin{array}{cc} a & b \\ c & d 
\end{array}\right) : \,\, t=(t_1,t_2) 
\mapsto t^\alpha=(t_1^a t_2^b, t_1^c t_2^d). $$ 
The action \eqref{salpha} is then given
equivalently as the endomorphisms of the
rings $\Q[t_i,t_i^{-1}]$ determined by setting 
\begin{equation}\label{salpha2}
s_\alpha: t\mapsto t^\alpha,
\end{equation}
for $t=(t_1,\ldots, t_n)$ and $t^\alpha$ as in
\eqref{talpha2}.

Notice that, unlike in the setting of endomotives defined in \cite{CCM}, here we allow semigroups
of endomotives that are non-abelian, like our $S_n=M_n(\Z)^+$. This contains interesting abelian
subsemigroups, in particular the diagonal subsemigroup
$S_{n,diag}=\N^n$ acting by
\begin{equation}\label{titiki}
s_k: P(t_i,t_i^{-1}) \mapsto P(t_i^{k_i},t_i^{-k_i}), \ \ \ i=1,\ldots,n, \ \ \  k=(k_1,\ldots,k_n)\in \N^n .
\end{equation}
on the ring of functions of $\bT^n$. 
We can similarly restrict to other interesting abelian subsemigroups of $M_n(\Z)^+$.
For example, one can consider the set of $n$-tuples
\begin{equation}\label{Sorder}
S_{n,ord}:=\{ (k_1,\ldots,k_n)\in \N^n\,|\,k_1\geq \cdots \geq k_n \}.
\end{equation}

\medskip

For $\alpha\in M_n(\Z)^+$, we denote by $X_\alpha$ the preimage
\begin{equation}\label{Xalpha}
X_\alpha =\{ t=(t_1,\ldots,t_n)\in \bT^n\,|\, s_\alpha(t)=t_0 \}.
\end{equation}
The $X_\alpha$ form an inverse system as in \cite{CCM} with the maps
\begin{equation}\label{xialphabeta}
 \xi_{\alpha,\beta}: X_\beta\to X_\alpha, \ \ \  t \mapsto t^\gamma, 
\end{equation} 
as in \eqref{talpha},  whenever $\alpha=\beta\gamma \in M_n(\Z)^+$. 
We denote by $X=\varprojlim_\alpha X_\alpha$ the inverse limit of this 
family and by $A=\varinjlim_\alpha A_\alpha$ the corresponding commutative 
algebra, with $X=\Spec(A)$.
We have 
\begin{equation}\label{Xalpha2}
X_\alpha = \Spec(\Q[t_i,t_i^{-1}]/(t^\alpha -t_0)) ,
\end{equation}
with $t^\alpha$ as in \eqref{talpha2}.

The multivariable analog of the BC 
endomotive is then given by the action on $X$
of the multiplicative semigroup $S_n =M_n(\Z)^+$. 
We have the following more concrete description
of the space $X(\bar\Q)$.

\begin{lem}\label{SpecAZ}
Any choice of an element $\rho\in GL_n(\hat\Z)$ determines 
an identification between the
set of algebraic points $X(\bar\Q)$ of the abelian part $X=\Spec(A)$ of the
multivariable Bost--Connes endomotive and the set
\begin{equation}\label{XbarQZ}
\rho(X(\bar\Q))= \cZ^n= \underbrace{\cZ\times\cdots\times\cZ}_{n-{\rm times}},
\end{equation}
where, as above, $\cZ$ is the set of all roots of unity in $\C$.
\end{lem}

\proof We show that the projective limit $X$ computed over the system of maps
$\xi_{\alpha,\beta}$ of \eqref{xialphabeta} is in fact the same as computing it over the
smaller system given by considering only the diagonal maps of \eqref{titiki}. Since 
$S_{n,diag}$ is a subsemigroup of $M_n(\Z)^+$, we have a surjective map 
of the projective limits
$$ X=\varprojlim_{\alpha\in M_n(\Z)^+} X_\alpha \rightarrow X'=\varprojlim_{k\in S_{n,diag}} X_k. $$
On the other hand, to see that we also have a map in the opposite direction, notice that, for any
given $\alpha \in M_n(\Z)^+$ we can find a $k=k(\alpha) \in S_{n,diag}$ such that there is a map
in the projective system \eqref{xialphabeta} 
\begin{equation}\label{xialphak}
 \xi_{\alpha,k(\alpha)}: X_{k(\alpha)} \to X_\alpha. 
\end{equation}
Let $k(\alpha)=(k_1,\ldots,k_n)$ be such that $k_i =\gcd\{ |\alpha_{ij}| \}_{j=1,\ldots, n}$. Then
there exists $\beta \in M_n(\Z)^+$ such that $\alpha = k(\alpha) \beta$. This implies that there is
a map in the projective system of the $X_\alpha$ of the form \eqref{xialphak}. The $X_{k(\alpha)}$
form an inverse system, since for $\alpha = \beta \gamma \in M_n(\Z)^+$ we have $\alpha k(\beta)^{-1}
k(\gamma)^{-1} \in M_n(\Z)^+$ hence $(k(\beta) k(\gamma))_i | \alpha_{ij}$ for all $j=1,\ldots,n$,
which means that $(k(\beta) k(\gamma))_i | k(\alpha)_i$ hence there is some $k\in S_{n,diag}$
such that $k(\alpha)=k(\beta) k$, which gives a map
\begin{equation}\label{Xkalpha}
 \xi_{k(\alpha),k(\beta)}: X_{k(\beta)} \to X_{k(\alpha)}. 
\end{equation}
This shows that we also have a map of inverse limits $X'\to X$ and that therefore the inverse
limits agree.

The inverse system \eqref{Xkalpha} of the $X_k$ for $k\in S_{n,diag}$ corresponds dually to
a direct system of algebras $A_k=\Q[t_i,t_i^{-1}]/(t^{k_i}-1)$, with $k\in \N^n$,
whose direct limit is then just $\Q[\Q/\Z]^{\otimes n}$. We then use an element of $\GL_n(\hat\Z)$
to give an identification of $(\Q/\Z)^n$ with $\cZ^n\subset \C^n$, as one does with $\hat\Z^*$ in the
case of the one-dimensional BC endomotive.
\endproof

The result of Lemma \ref{SpecAZ} is simply showing that the set $X(\bar\Q)$ that defines the
analytic endomotive of the multivariable Bost--Connes system consists of all the torsion points
of the $n$-torus $\bT^n$.
In particular, it follows that
the algebra $C(X(\bar\Q))$ has an explicit 
presentation with generators 
$e(r_1)\otimes \cdots \otimes e(r_n)$, 
for $r_i\in \Q/\Z$, where the $e(r_i)$
are the generators of the abelian part 
of the 1-dimensional BC algebra and
satisfying the same relations.

To construct the noncommutative crossed product algebra 
for the multivariable BC endomotives, we
consider, as in the one-variable case, isometries 
implementing the action of the semigroup of
endomorphisms. As in the previous discussion for
the multivariable Habiro ring, we consider 
isometries $\mu_\alpha$ as in \eqref{mualpharels}.
Notice that these can be represented as isometries
on the Hilbert space $\cH_n = \ell^2(M_n(\Z)^+)$
by setting
\begin{equation}\label{mualpha}
\mu_\alpha \epsilon_\beta = \epsilon_{\alpha\beta}.
\end{equation}
with the adjoints acting as
\begin{equation}\label{mualphastar}
\mu_\alpha^* \epsilon_\beta = \left\{ \begin{array}{ll} \epsilon_\gamma & 
\beta=\alpha\gamma \in M_n(\Z)^+ \\ 0 & \text{otherwise}.\end{array}\right.
\end{equation}

Similarly, we have representations of the abelian
part $C(X(\bar\Q))$ of the algebra on the same
Hilbert space by the following.

\begin{lem}\label{reperi}
The choice of an element $\rho\in \GL_n(\hat\Z^*)$ 
determines a representation of the
algebra $C(X(\bar\Q))$ on $\ell^2(M_n(\Z)^+)$. 
\end{lem}

\proof This is given on the generators
$e(r_1)\otimes \cdots \otimes e(r_n)$ by
\begin{equation}\label{XactMn}
(e(r_1)\otimes \cdots \otimes e(r_n))\, \epsilon_\beta = 
\left(\prod_{i=1}^n \zeta^{\tilde\beta}_i\right) \, 
\epsilon_\beta,
\end{equation}
where $\zeta=(\zeta_{r_1},\ldots,\zeta_{r_n})=
\rho(r_1,\ldots,r_n)\in \cZ^n$, with the transpose
$\tilde\beta = \beta^\dag \in M_n(\Z)^+$
and $\zeta^{\tilde\beta}=(\zeta^{\tilde\beta}_i)$ 
is defined as in \eqref{talpha2}. 
\endproof

Let $\sigma_\alpha$ for $\alpha \in M_n(\Z)^+$
denote the endomorphisms of the algebra 
$C(X(\bar\Q))$ induced by the $s_\alpha$ of \eqref{salpha}.
The relation between the isometries $\mu_\alpha$ and
the endomorphisms $\sigma_\alpha$ is as follows.

\begin{lem}\label{sigmamualpha}
The endomorphisms $\sigma_\alpha$ of the multivariable 
BC endomotive defined by the torus maps \eqref{salpha} 
are given by
\begin{equation}\label{sigmaalphamu}
\sigma_\alpha (e(\underline{r})) = 
\mu_\alpha^* \, e(\underline{r}) \, \mu_\alpha,
\end{equation}
for $\underline{r}=(r_1,\ldots,r_n) \in (\Q/\Z)^n$.
\end{lem}

\proof For $e(\underline{r})=e(r_1)\otimes \cdots 
\otimes e(r_n)$, 
the operator 
$\mu_\alpha^* \,e(\underline{r}) \, 
\mu_\alpha$ 
acts on a basis element $\epsilon_\beta$ by
$$ \mu_\alpha^* \, e(\underline{r}) \,
\mu_\alpha\, \epsilon_\beta =
\left(\prod_{i=1}^n \zeta^{\widetilde{\alpha\beta}}_i\right) \, \epsilon_\beta = 
 \left(\prod_{i=1}^n \zeta^{\tilde\beta\tilde\alpha}_i\right) \, \epsilon_\beta =
\left(\prod_{i=1}^n (\zeta^{\tilde\alpha})^{\tilde\beta}_i\right) \, \epsilon_\beta = 
 \left(\prod_{i=1}^n \xi^{\tilde\beta}_i\right) \, \epsilon_\beta , $$
where $\xi=\zeta^{\tilde\alpha}$.
\endproof

Similarly, the following result gives the direct analog of the relation
$$ \mu_k e(r) \mu_k^* =\frac{1}{k} \sum_{ks=r} e(s) $$
of the 1-dimensional case.

\begin{lem}\label{mualphaer}
The operators $\mu_\alpha$ and $e(\underline{r})=e(r_1)\otimes\cdots\otimes e(r_n)$ satisfy 
the relation
\begin{equation}\label{crossermualpha}
\mu_\alpha \, e(\underline{r})\,
 \mu_\alpha^* = \frac{1}{\det\alpha} 
\sum_{\alpha(\underline{s})=\underline{r}} e(\underline{s}).
\end{equation}
for $\underline{r}=(r_1,\ldots,r_n)$ and 
$\alpha(\underline{s})=(\alpha(\underline{s})_i)
\in (\Q/\Z)^n$.
\end{lem}

\proof We first show that the number of solutions to $\alpha(\underline{s})=\underline{r}$ is
equal to $\det(\alpha)\in \N$. In fact, the number of such solutions is equal to the number of 
solutions of the equation $\alpha(\underline{s})=0$. Any $\underline{s}=(s_i)$ such that all the
$s_i$ satisfy $\det(\alpha) s_i =0$ is a solution, since we can write $\alpha = \det(\alpha) \hat\alpha$
so that if $\det(\alpha) \underline{s}=0$ then also $\alpha(\underline{s})=0$. Conversely, a solution
of $\alpha(\underline{s})=0$ will also satisfy $\det(\alpha) s_i =0$. In fact, the entries $s_i$ are  
torsion of orders that divide $\det(\alpha)$. This one can see easily in the case of 
$\alpha\in M_2(\Z)^+$, where
$$ \left\{ \begin{array}{ll} a s_1 + b s_2 =0 \\ c s_1 + d s_2 =0 \end{array}\right. \Rightarrow 
(ad-bc) s_1 =0 \ \ \text{ and } \ \  (ad-bc)s_2 =0. $$
The general case is completely analogous. We then have 
$$ \mu_\alpha \, e(\underline{r})\, 
\mu_\alpha^*\, \epsilon_\beta = 0 $$
if $\alpha^{-1}\beta \notin M_n(\Z)^+$, while otherwise
$$ \mu_\alpha \, e(\underline{r})\, \mu_\alpha^*\, 
\epsilon_\beta = 
\left(\prod_{i=1}^n \zeta^{\widetilde{\alpha^{-1}\beta}}_i\right) \, \epsilon_\beta = 
 \left(\prod_{i=1}^n \zeta^{\tilde\beta\widetilde{\alpha^{-1}}}_i\right) \, \epsilon_\beta =
\left(\prod_{i=1}^n (\zeta^{\widetilde{\alpha^{-1}}})^{\tilde\beta}_i\right) \, \epsilon_\beta = 
\left(\prod_{i=1}^n \xi^{\tilde\beta}_i\right) \, \epsilon_\beta , $$
where $\xi$ is such that $\xi^{\tilde\alpha}=\zeta$. Notice that, if $\xi_1$ and $\xi_2$ both
satisfy $\xi_k^{\tilde\alpha}=\zeta$, for $k=1,2$, we have
$$ \zeta^{\widetilde{\alpha^{-1}\beta}} = (\xi_k^{\tilde\alpha})^{\widetilde{\alpha^{-1}\beta}}
= (\xi_k^{\tilde\alpha})^{\tilde\beta \widetilde{\alpha^{-1}}} = \xi_k^\beta . $$
Thus, the expression above for $\mu_\alpha e(\underline{r})\mu_\alpha^*$ does
not depend on the choice of a $\xi$ with $\xi^\alpha =\zeta$, \ie of an $\underline{s}$ with
$\alpha(\underline{s})=\underline{r}$, so that it can be written equivalently in the form
\eqref{crossermualpha}.
\endproof

The noncommutative algebra of the multivariable
BC endomotive is  the rational semigroup crossed product
algebra
\begin{equation}\label{BCncross}
\cA_{\Q,BC,n}= A\rtimes_\rho S_n \cong \Q[\Q/\Z]^{\otimes n} 
\rtimes_\rho M_n(\Z)^+, 
\end{equation}
and the $C^*$-algebra that gives the 
analytic endomotive is obtained as 
\begin{equation}\label{BCncrossan}
\cA_{BC,n}= C(X(\bar\Q))\rtimes_\rho S_n 
\cong C^*(\Q/\Z)^{\otimes n} 
\rtimes_\rho M_n(\Z)^+,
\end{equation}
where the semigroup action is given by
$\rho_\alpha(e(\underline{r}))=\mu_\alpha\, e(\underline{r})
\, \mu_\alpha^*$, as in \eqref{crossermualpha}.

It is natural to consider on the algebra $\cA_{BC,n}$ a time
evolution and quantum statistical mechanical properties that
generalize the corresponding ones of the BC system. However,
as we discuss briefly here and more in detail in \S 
\ref{ConvolAlgSec} below, the natural generalization of
the BC dynamics runs into a problem coming from the
presence of infinite multiplicities in the spectrum of 
the Hamiltonian. We show in \S \ref{ConvolAlgSec} that
this problem can be resolved following the same method
used in the case of 2-dimensional $\Q$-lattices in 
\cite{CM1}.

\begin{lem}\label{multitevol}
The multivariable Bost--Connes algebra $\cA_{BC,n}=
C^*(\Q/\Z)^{\otimes n} \rtimes_\rho M_n(\Z)^+$
has a time evolution $\sigma:\R \to \Aut(\cA_n)$ 
defined by setting
\begin{equation}\label{sigmatmuer}
\sigma_t(\mu_\alpha) = \det(\alpha)^{it}\, \mu_\alpha , \ \ \forall \alpha \in M_n(\Z)^+, \ \ \ \text{ and } \ \ \  \sigma_t(e(\underline{r}))=e(\underline{r}) \ \  \forall \underline{r}\in (\Q/\Z)^n. 
\end{equation}
Let $H$ be the Hamiltonian $H$ that generates this time 
evolution in a representation $\pi_\rho$ of $\cA_{BC,n}$ 
on $\cH_n$, with $\rho\in \GL_n(\hat\Z)$. The operator $H$ 
has eigenvalues $\log \det(\alpha)$ with infinite 
multiplicities.
\end{lem}

\proof It is easy to check that \eqref{sigmatmuer} 
defines a 1-parameter family
$\sigma:\R \to \Aut(\cA_{BC,n})$ of automorphisms. 
Recall that a self-adjoint positive (unbounded)
linear operator $H$ on $\cH_n$ is a Hamiltonian for 
the time evolution $\sigma_t$, in a representation
$\pi_\rho$ of the algebra $\cA_{BC,n}$, if it satisfies
$$ \pi_\rho(\sigma_t(a)) =e^{it H} \, \pi_\rho(a) \, 
e^{-itH}, \ \ \  \forall a\in \cA_{BC,n}, \ \ \ 
\forall t\in \R. $$
Let $H$ be the operator on $\cH_n$ defined by
\begin{equation}\label{Halpha}
 H \, \epsilon_\beta = \log \det(\beta)\, \epsilon_\beta.
\end{equation} 
We have
$$ e^{itH}\, \mu_\alpha e^{-itH} \epsilon_\beta = \det(\beta)^{-it} \det(\alpha\beta)^{it} \epsilon_{\alpha\beta}= \det(\alpha)^{it} \epsilon_{\alpha\beta} = \sigma_t(\mu_\alpha) \epsilon_\beta. $$
Similarly
$$ e^{itH}\, \pi_\rho(e(\underline{r}))\, e^{-itH} 
\epsilon_\beta = \pi_\rho(e(\underline{r})) \epsilon_\beta =
\pi_\rho(\sigma_t(e(\underline{r}))) \epsilon_\beta. $$
Thus $H$ is a Hamiltonian for the time evolution in all the
representations $\pi_\rho$ with $\rho\in \GL_n(\hat\Z)$. 
Any other choice of Hamiltonian differs by a constant,
$H' = H + \log \lambda$, for some $\lambda \in \R^*_+$. 
The Hamiltonian has infinite multiplicities
in the spectrum, since all basis vectors 
$\epsilon_{\gamma\beta}$ and $\epsilon_{\beta\gamma}$
with $\gamma\in \SL_n(\Z)$ are eigenvalues with the 
same eigenvector $\log\det(\beta)$.
\endproof

\subsection{Multivariable BC endomotive and automorphisms}\label{ConvolAlgSec}

One new feature of the multivariable case is the 
presence of a large subgroup of automorphisms
in the semigroup of endomorphisms defining the 
endomotive, here given by 
$\SL_n(\Z)\subset M_n(\Z)^+$. We have seen above how 
this creates a problem of infinite multiplicities 
in the spectrum of the Hamiltonian generating the 
time evolution. The same phenomenon occurred already in
the construction of the quantum statistical mechanical 
system for 2-dimensional $\Q$-lattices in \cite{CM1} and 
we will treat it here in a similar manner. Namely, instead
of working with the $C^*$-algebra 
$C(X(\bar\Q))\rtimes M_n(\Z)^+$, for the purpose of
quantum statistical mechanics we can replace it with a
convolution algebra where we mod out the automorphisms.

\smallskip

The algebra, representations, and time evolution we describe here are
a minor variant over the one introduced in \cite{CM1},  as ``determinant
part of the $\GL_2$-system". Here we deal with the higher rank case of
$\GL_n$ and with vectors in $\hat\Z^n$ instead of matrices in $M_n(\hat\Z)$,
but we make use of the same technique of passing to a type II$_1$ factor
and corresponding trace to eliminate the infinite multiplicities, as in 
\cite{CM1}, \S 1.7. 

\subsubsection{Groupoid, convolution algebra, and type II$_1$ factor}

We let
\begin{equation}\label{cU}
\cU =\{ (\alpha,\rho)\in \GL_n(\Q)^+ \times \hat\Z^n \,|\, \alpha(\rho)\in \hat\Z^n \}.
\end{equation}
The set $\cU$ is a groupoid with source and 
target maps $s(\alpha,\rho)=\rho$ and $t(\alpha,\rho)=\alpha(\rho)$ 
and with composition law $(\alpha_2,\rho_2)\circ (\alpha_1,\rho_1)=
(\alpha_2\alpha_1, \rho_1)$ if $\rho_2= \alpha_1(\rho_1)$. The inverse
of the arrow $(\alpha,\rho)$ in the groupoid $\cU$ is the arrow $(\alpha^{-1}, \alpha(\rho))$.

We also define a 2-groupoid $\cU_\Gamma$ where the objects and 1-morphisms are those
of the groupoid $\cU$ and the 2-morphisms are given by the action
$$ (\gamma_1,\gamma_2): (\alpha, \rho) \mapsto (\gamma_1 \alpha \gamma_2^{-1}, 
\gamma_2(\rho)), $$
of $\Gamma \times \Gamma$, with $\Gamma =\SL_n(\Z)$. This 2-groupoid is a way to represent
the quotient
\begin{equation}\label{UGamma}
\{ (\alpha,\rho) \in \Gamma\backslash \GL_n(\Q)^+ \times_\Gamma \hat\Z^n \,|\,
\alpha(\rho) \in \hat\Z^n \},
\end{equation}
in a way that avoids actually taking the quotient by the $\Gamma \times \Gamma$-action,
which would not lead to a good quotient in the ordinary sense.
The 2-groupoid $\cU_\Gamma$ provides
a version of the multivariable Bost--Connes algebras that eliminates the infinite multiplicities
due to the action of $\Gamma$.  In the case of a 2-groupoid, as in the case of more general
2-categories (see also \S \ref{semiZhsSec} below and \cite{MaZa}) one can construct an
associated convolution algebra by considering compactly supported continuous 
functions $f: \cU_\Gamma \to \C$, \ie functions $f(\alpha,\rho,\gamma_1,\gamma_2)$
with the two convolution products of vertical and horizontal composition of 2-morphisms.

Following the analogous case of \cite{CM1}, on can instead work with the groupoid $\cU$
itself, and eliminate the infinite multiplicities in the spectrum of the Hamiltonian 
caused by the $\Gamma$-symmetry by working with an associated type II$_1$ factor and
the corresponding trace instead of the ordinary trace.

The convolution algebra $C_c(\cU)$ with the product dictated by the groupoid
composition law
\begin{equation}\label{convprodU}
(f_1 \star f_2)(\alpha,\rho) =\sum_{(\alpha,\rho)=(\alpha_1,\rho_1)
\circ (\alpha_2,\rho_2) \in \cU_\Gamma} f_1(\alpha_1,\rho_1) f_2 (\alpha_2,\rho_2),
\end{equation}
is just of the form
\begin{equation}\label{convprodU2}
(f_1 \star f_2)(\alpha,\rho) =\sum_{\beta \in \GL_n(\Q)^+\,:\, \beta\rho \in \hat\Z^n} 
f_1(\alpha\beta^{-1}, \beta(\rho)) f_2(\beta,\rho), 
\end{equation}
and with the adjoint $f^*(\alpha,\rho)= \overline{f(\alpha^{-1}, \alpha(\rho))}$. 
This algebras has a family of representations determined by the choice of an element
$\rho \in \hat\Z^n$, the set of units of the groupoid $\cU$, by setting
\begin{equation}\label{repUGamma}
(\pi_\rho (f) \xi) (\alpha) =\sum_{\beta \in \Gamma \backslash \GL_n(\Q)^+\,:\, 
\beta\rho\in \hat\Z^*} f( \alpha\beta^{-1} , \beta (\rho)    ) \xi(\beta).
\end{equation}
on the Hilbert space $\ell^2(G_\rho)$, where 
$G_\rho=\{ \alpha\in \GL_n(\Q)^+\,|\, \alpha(\rho) \in \hat\Z^n \}$.  
The algebra $C_c(\cU)$ can be completed with respect to the norm
$$\| f\| =\sup_{\rho\in \hat\Z^n} \| \pi_\rho(f) \|_{\cB(\ell^2(G_\rho))} $$
and one denotes by $C^*(\cU)$ the resulting $C^*$-algebra. This gives
an equivalent description of the multivariable BC endomotive. The time
evolution described above corresponds to
\begin{equation}\label{timeUGamma}
\sigma_t(f) (\alpha,\rho) = \det(\alpha)^{it} f(\alpha,\rho).
\end{equation}
In fact, it is clear from \eqref{convprodU} that 
$\sigma_t (f_1\star f_2)= \sigma_t(f_1)\star \sigma_t(f_2)$
and that $\sigma_t(f^*)=\sigma_t(f)^*$. 

We consider as before the representations $\pi_\rho$ of the algebra of the multivariable
BC system, described here as the groupoid algebra $C^*(\cU)$ on the 
Hilbert space $\ell^2(M_n(\Z)^+)$, where $\rho \in (\hat\Z^*)^n$. 
The group $\Gamma =\SL_n(\Z)$ acts
on the right on the same Hilbert space by $\epsilon_\alpha \mapsto \epsilon_{\alpha\gamma^{-1}}$,
so that this action commutes with the action of $C^*(\cU)$. This action of $\Gamma$
determines a type II$_1$ factor given by the image of the group algebra $C^*(\Gamma)$ in
the algebra of bounded operators $\cB(\ell^2(M_n(\Z)^+))$. Thus, instead of the usual trace,
we can use as in \cite{CM1} the trace $\Tr_\Gamma$ associated to this type II$_1$ factor,
which eliminates the infinite multiplicities in the spectrum of the Hamiltonian. The same
argument given in \S 1.7 of \cite{CM1} gives then the following result.

\begin{lem}\label{timeZstar}
In the representation $\pi_\rho$ on $\ell^2(M_n(\Z)^+)$
the time evolution is generated by the Hamiltonian
\begin{equation}\label{HamGamma}
H \, \epsilon_m = \log \det(m) \,\, \epsilon_m, \ \ \ \text{ for } \ m\in M_n(\Z)^+/\Gamma.
\end{equation}
The partition function, computed with respect to the type II$_1$-trace, is then of the form
\begin{equation}\label{ZbetaGamma}
Z(\beta) =\Tr_\Gamma(e^{-\beta H})= \sum_{m\in M_n(\Z)^+/\Gamma} \det(m)^{-\beta} =\prod_{k=0}^{n-1} \zeta(\beta-k).
\end{equation}
\end{lem}

\proof The argument is completely analogous to that of Lemma 1.18 of \cite{CM1} 
for the $\GL_2$-case.
\endproof

Notice that, in these representations the Hamiltonian has positive energy, since on $M_n(\Z)^+$ we have $\log\det(m)\geq 0$.

In particular, at sufficiently low temperature, KMS states of Gibbs type can be
constructed as in Proposition 1.19 of \cite{CM1}. Roughly, they will be of the form
\begin{equation}\label{KMSlow}
\varphi_{\beta,\rho}(f) = Z(\beta)^{-1} \Tr_\Gamma(\pi_\rho(f) e^{-\beta H})= Z(\beta)^{-1} \sum_{m\in M_n(\Z)^+/\Gamma} f(1, m(\rho)) \det(m)^{-\beta} ,
\end{equation}
parameterized by the elements $\rho\in (\hat\Z^*)^n$, though more precisely one first
defines them on functions that depend on $\hat\Z^n$ through a projection $p_N: \hat\Z \to \Z/N\Z$
and then extend them to the algebra $C^*(\cU)$.
At zero temperature these give rise, by weak limits, to KMS states of the form
\begin{equation}\label{ZeroTKMSn}
\varphi_{\infty,\rho}(f) = f(1,\rho).
\end{equation}

\subsection{The dual system and zeta functions}\label{KuroSec}

As in the case of the original Bost--Connes system, one
can consider for the multivariable cases the dual system
obtained by taking the crossed product with the time evolution.
As shown in \cite{CCM} and \cite{CCM3}, as well as in Connes'
earlier work on the Riemann zeta function \cite{CoZeta}, the
dual system of the Bost--Connes system is the natural 
noncommutative space that supports the spectral realization
of the zeros of the Riemann zeta function, on a ``noncommutative
motive'' describing the complement of the classical points of
this space (see \cite{CCM3} for a detailed analysis). It is
then natural to ask what the corresponding procedure (called
``cooling and distillation'' in \cite{CCM}) applied to the
multivariable Bost--Connes systems introduce here will give
in terms of L-functions. It is natural to expect that something
like the Kurokawa $L$-function \cite{Kuro} should appear in this setting.
These motivically should correspond to the conjectural tensor 
product over $\Spec \Z$ as explained in Manin's \cite{ManZeta}.
This is beyond the scope of the present paper and we hope to
return to this question in future work.

In particular, as explained in \cite{CCM3} (see also Chapter 4 of 
\cite{CMbook}), there is a well developed dictionary of analogies
between the Weil proof of the Riemann hypothesis for function
fields and the approach for number fields via noncommutative
geometry developed in \cite{CoZeta}, \cite{CCM}, \cite{CCM3}.
The action of Frobenius of \'etale cohomology is replaced by
the scaling action on the cyclic homology of a cokernel (in the
abelian category of cyclic modules defined in \cite{CoExt}) of
the restriction map from the dual system $\cA_{BC} \rtimes_\sigma \R$
to its ``space of classical points''. 

Geometrically, if one hopes to transpose the main ideas in the Weil proof 
from the function field to the number field case, one needs 
a good analog for number fields of the geometry of the curve $C$ for 
function fields $\K=\F_q(C)$, and in particular an analog of 
correspondences on $C\times_{\Sp(\F_q)} C$. Possible analogs are
identified in \cite{CCM3} in terms of the noncommutative geometry
of the ad\`ele class space. 

The idea of developing geometry over $\F_1$ is also aimed at a
similar purpose, and that's why identifying precisely its relation to 
noncommutative geometry via the BC endomotive is an important
ingredient of the general picture.  From the point of view of geometry 
over $\F_1$, one would like to view $\Sp(\Z)$ as a direct analog of 
the curve $C$. However, as already observed in \cite{Man} and
\cite{Soule}, $\Sp(\Z)$ is  not of finite type over $\F_1$. The geometric
space where a possible analog of the Weil proof would be taking 
place should then be $\Sp(\Z) \otimes_{\Sp(\F_1)} \Sp(\Z)$. 
Although this object does  not have a good definition within 
$\F_1$-geometry, some of its expected properties were 
identified by Manin in \cite{ManZeta} in terms of the zeta 
functions of Kurokawa \cite{Kuro}.

Kurokawa's approach to $L$-functions that would be associated
to spaces $\Sp(\Z) \times_{\F_1}\times \cdots \times_{\F_1} \Sp(\Z)$
is by sums over zeros of the Riemann zeta function. More precisely,
the Kurokawa tensor product of zeta functions corresponds to
regularized infinite products of the form (\cf \cite{Man} (2.30)-(2.31))
\begin{equation}\label{Kurokawa}
\prod_{\lambda \in \Xi} (s-\lambda)^{m_\lambda} \otimes \prod_{\mu \in \Theta}
(s-\mu)^{n_\mu} := \prod_{(\lambda,\mu)} (s-\lambda-\mu)^{m_\lambda + n_\mu},
\end{equation}
where the infinite products are computed as
\begin{equation}\label{regprod}
\prod_{\lambda \in \Xi} (s-\lambda)^{m_\lambda} = \exp ( -\frac{d}{dz} (\Gamma(z)^{-1}
\int_0^\infty \sum_\lambda m_\lambda e^{(s-\lambda)t} t^{z-1} \, dt )|_{z=0}).
\end{equation}
In particular, for $\rho$ ranging over the critical zeros of the Riemann zeta
function one has (see \cite{Man}, (1.5)) 
\begin{equation}\label{zetazeros}
\prod_\rho \frac{s-\rho}{2\pi} = \frac{s(s-1)}{4 \pi^2} 2^{-1/2} \pi^{-s/2} \Gamma(s/2) \zeta(s)
\end{equation}
Kurokawa uses a splitting into the product (\cite{Man}, \S 5.2)
$$ \xi_{\pm}(s):=\prod_{\sign \Im(\rho) =\pm 1} (1-\frac{s}{\rho}) e^{s/\rho} $$
to deal with the infinite products and defined the tensor product of the completed
Riemann zeta function $\xi(s) = s(s-1) \Gamma_\R(s) \zeta(s)$ to be
\begin{equation}\label{Kuroprodzeta}
\xi(s)^{\otimes n} = \xi_+(s)^{\otimes n}\left( \xi_-(s)^{\otimes n} \right)^{(-1)^{n-1}}.
\end{equation}
In the noncommutative geometry approach to the Riemann hypothesis developed
in \cite{CoZeta} the zeros of the Riemann zeta function arise naturally from the 
scaling action on the cokernel of the restriction map from the noncommutative
adeles class space to the classical spaces of ideles classes. In the formulation
of \cite{CCM}, \cite{CCM3}, this action is on the cyclic homology of the
cyclic module describing this cokernel, with the adeles class space arising from
the dual system of the BC quantum statistical mechanical system. It is therefore
natural to view the multivariable versions of the BC system considered here as
a possible candidate for endomotives where the sums over zeros of zeta used
in the construction of the Kurokawa zeta functions should arise by following the
general procedure of ``cooling and distillation" of endomotives described in
\cite{CCM}. We intend to return to a more detailed analysis of this problem in
future work. 

Another source of evidence for the fact that the multivariable BC 
endomotives introduced here should be related to the missing 
geometry of the products $\Sp(\Z) \times_{\F_1}\times \cdots \times_{\F_1} \Sp(\Z)$
comes from a different source, namely the relation between BC endomotives
and $\Lambda$-rings, which we discuss in detail in \S \ref{LambdaSec} here
below. In fact, recently, an approach to geometry over $\F_1$ in terms of
$\Lambda$-rings proposed in \cite{Bor}, \cite{Bor2}, \cite{BoSm}
suggests that the missing $\Sp(\Z) \times_{\F_1}\times \cdots \times_{\F_1} \Sp(\Z)$
should be related to the space of big Witt vectors $W(\Sp(\Z))$.

\section{Endomotives and $\Lambda$-rings}\label{LambdaSec}

The notion of $\Lambda$-ring was introduced by Grothendieck \cite{Groth}
to formalize the properties of the exterior power operations on Grothendieck
groups and $K$-theory. They provide a useful notion in the context of 
characteristic classes and the Riemann--Roch theorem, see \cite{FuLa}.

\smallskip

It was shown recently in Borger--de Smit \cite{BoSm} that the study of 
abstract properties of $\Lambda$-rings has interesting number theoretic 
implications. The formulation of the $\Lambda$-ring structure adopted 
in \cite{BoSm} is based on an equivalent reformulation given in \cite{Wil} 
of the original Grothendieck definition. 

\smallskip

In particular, as argued in greater detail in \cite{Bor}, the structure of $\Lambda$-ring
can be thought of as a descent condition for a ring from $\Z$ to $\F_1$, in the sense
that it imposes a strong compatibility condition on the Frobenius at all the different
reductions modulo primes. This gives a different approach to defining varieties and
schemes over $\F_1$, where the compatibility is expressed in terms of a semigroup
action, rather than through the cyclotomic points as in the Soul\'e approach. One of
the main points we want to stress in this paper is the fact that the BC endomotive
encodes the compatibility between these two approaches, in the sense that not
only it fits the Soul\'e definition of affine (pro)varieties over $\F_1$ as proved in
\cite{CCM2}, but it also fits the $\Lambda$-ring approach to the definition of varieties 
over $\F_1$ in the sense of \cite{Bor}, \cite{BoSm}. In fact, it is a universal object both 
for the Soul\'e approach, where it encodes the tower of extensions $\F_{1^n}$ of
\cite{KaSmi} \cite{Soule}, with their
Frobenius actions, and also in the sense of \cite{BoSm}, because (together with
its multivariable versions introduced here) it gives universal $\Lambda$ rings into
which any $\Lambda$-ring that is torsion free and finite over $\Z$ embeds.

\smallskip

We are going to show here that the integral model of the BC
endomotive constructed in \cite{CCM2} is a an {\em integral $\Lambda$ model}
in the sense of  \cite{BoSm}. Moreover, it is universal in the sense that any
integral $\Lambda$ model admits a map of $\Lambda$-rings to a (multivariable)
integral BC endomotive.

\smallskip

Following  \cite{BoSm}, we give the following definition of $\Lambda$-structure
on a commutative ring and on finite dimensional reduced $\Q$-algebra.

\begin{defn}\label{Lambdaring}
Let $R$ be a commutative ring such that the underlying
abelian group is torsion free. An (integral) $\Lambda$-ring structure on $R$ is an action of the
abelian multiplicative semigroup $\N$ generated by the primes as endomorphisms
of the ring $R$, with the property that the endomorphism $\sigma_p$ lifts the
Frobenius map modulo $p$, \ie
\begin{equation}\label{liftFrob}
s_p(x) -x^p \in pR, \ \ \  \forall x\in R.
\end{equation}
A map of (torsion free) $\Lambda$-rings is a ring homomorphism compatible with
the semigroup actions, \ie satisfying $f\circ s_k = s_k \circ f$ for all $k\in \N$.
\end{defn}

A typical example of a $\Lambda$ structure on a ring is the following case.

\begin{ex}\label{ex1}  Consider the ring $R=\Z[t,t^{-1}]/(t^n-1)$
with the endomorphisms $s_k(P)(t,t^{-1})=P(t^k,t^{-k})$ for $k\in \N$, as in \eqref{ttk}.
It is a $\Lambda$-ring in the sense of Definition \ref{Lambdaring}, satisfying
\eqref{liftFrob}. 
\end{ex}

This will be the basic example we are interested in to establish the connection between
$\Lambda$-rings and the multivariable BC endomotives.

As in \cite{BoSm}, after tensoring with $\Q$, one obtains a $\Q$-algebra $A=R\otimes \Q$, with an
action of the semigroup $\N$ by endomorphisms, now without the condition \eqref{liftFrob}.
We assume here, as in \cite{BoSm}, to be working with $\Lambda$ rings that are torsion free
and finite over $\Z$. As in the case of endomotives, these correspond to zero-dimensional
geometries.

\begin{defn}\label{LambdaQalg}
Let $A$ be a finite dimensional reduced algebra over $\Q$. Then $A$ has a $\Lambda$-ring
structure if it has an action of the multiplicative group $\N$ by endomorphisms. This $\Lambda$-ring
structure has an {\em integral model} if $A=R\otimes \Q$, where $R$ is a $\Lambda$ ring as in
Definition \ref{Lambdaring} and the inclusion of $R$ as a subring of $A$ is a morphism of
$\Lambda$-rings. 
\end{defn}

We also give the following definition that will be useful in comparing $\Lambda$-rings
and endomotives.

\begin{defn}\label{LambdaLim}
A direct limit of $\Lambda$-rings is a direct system of $\Q$-algebras $A_\alpha$ such
that all the maps in the direct system are morphisms of $\Lambda$-rings.
\end{defn}

This induces a $\Lambda$-ring structure (\ie an action of $\N$ by endomorphisms)
on the direct limit $A=\varinjlim A_\alpha$. 
Consider then the BC endomotive. We have the following result. 

\begin{lem}\label{AnLambda}
The Bost--Connes endomotives is a direct limit of $\Lambda$-rings with compatible 
integral models.
\end{lem}

\proof The abelian algebra of the BC endomotive is a direct limit
$A=\Q[\Q/\Z]= \varinjlim_n A_n$, where $A_n=\Q[\Z/n\Z]$. At each finite 
level the algebra $A_n=\Q[\Z/n\Z]$ has a $\Lambda$-ring structure
given by the semigroup action $\sigma_k: A_n \to A_n$ is the one given by 
\eqref{sigmanBC}, which correspond to the $s_k$ of \eqref{ttk}. The actions
$\sigma_k$ are compatible with the maps of the direct system since we have
$$ \xi_{m,n} \circ \sigma_k =\sigma_k \circ \xi_{m,n}, \ \ \  \text{for } n|m \ \text{ and } \ \forall k\in \N, $$
for $s_k$ as in \eqref{ttk} and with  $\xi_{m,n}: X_n \to X_m$ the maps of the inverse system with 
$X_n=\Spec(\Q[t]/(t^n-1))$. 

\smallskip

The integral model of the BC endomotive constructed in \cite{CCM2} is obtained 
by taking as abelian part the ring $\Z[\Q/\Z]=\varinjlim \Z[\Z/n\Z]$ and the action of 
$\N$ by endomotives given by the same $\sigma_n$ as above. 
(Unlike the $\rho_n$ these do not involve denominators and are 
therefore defined over $\Z$.) Thus, at each finite level the $A_n=R_n \otimes \Q$, with
$R_n =\Z[\Z/n\Z]$ are $\Lambda$-rings with integral models. The compatibility with the
maps in the direct system give the result for $A=R\otimes \Q$ with the resulting $\N$
action by endomorphisms. To see that the integral model of the BC endomotive
defines, at each finite level $R_n$ an integral $\Lambda$-ring in the sense of Definition
\ref{Lambdaring}, we only need to check that the compatibility with the Frobenius
is satisfied for all primes $p$. The compatibility of the action of $\sigma_{p^\ell}$ on the
integral model of the BC system with the Frobenius action follows from 
the case of Example \ref{ex1}.
\endproof

The fact that the endomorphisms action of $\N$ on the integral model of the
abelian part of the BC algebra is compatible with the action of the
Frobenius in all the reductions modulo primes is discussed in detail in \cite{CCM2},
where this compatibility is an important part of the structure in viewing the BC
endomotive as defined over $\bF_1$, with the abelian part providing the system of
extensions $\bF_{1^n}$ and the endomorphisms $\sigma_p$ providing a compatible 
action that induces the Frobenius correspondences in characteristic $p>0$. In fact,
the compatibility with Frobenius action as formulated for the BC system
in \cite{CCM2} is stronger than the condition \eqref{liftFrob} used in the context 
of $\Lambda$-rings. It would be interesting to investigate whether stronger versions
of \eqref{liftFrob} of the kind considered in \cite{CCM2}, involving reductions mod p
of the endomotive, may also be implemented for $\Lambda$-rings and what 
kind of structure they give rise to.

\smallskip

We also have the following general result that produces endomotives from
$\Lambda$-rings.

\begin{lem}\label{LambdaEndomotives}
A direct limit of $\Lambda$-rings defines an endomotive over $\Q$. This has an
integral model if at each finite level the $\Lambda$-rings have an integral model
in the sense of Definition \ref{LambdaQalg}.
\end{lem}

\proof We have a direct system of finite dimensional reduced algebras $A_\alpha$ 
over $\Q$. Consider the crossed product algebra $A\rtimes \N$, where the 
commutative $A$ is the direct limit $A=\varinjlim_\alpha A_\alpha$ and the crossed
product with $\N$ is implemented by 
$$ \rho_n(a) = \mu_n \, a \, \mu_n^*, $$
where the $\mu_n$ are isometries satisfying
$$ \mu_n^* \, a \, \mu_n = \sigma_n(a) . $$ 
Unlike the $\sigma_n$, the $\rho_n$ do not preserve levels $A_\alpha$.
This is the noncommutative algebra of an algebraic endomotive over $\Q$,
according to the definition given in \cite{CCM}. Notice moreover that
$\Lambda$-rings given by finite-dimensional reduced commutative algebras 
over $\Q$ have a continuous action of the Galois group $G=\Gal(\bar\Q/\Q)$ on
$$ \cX_\alpha = \Hom(A_\alpha, \bar\Q), $$
as observed in \cite{BoSm}.  The induced action on $\cX=\Hom(A,\bar\Q)$ is
the action of $G$ defined for endomotives in \cite{CCM} given by composition
$$ A\stackrel{\chi}{\to} \bar\Q \stackrel{\gamma}{\to} \bar\Q, $$
for $\chi \in \cX$ and $\gamma\in G$. The analytic endomotive is then
obtained as the $C^*$-algebra $C(\cX)\rtimes \N$ with the induced action of $G$
by automorphisms.

If the $\Lambda$ rings $A_\alpha$ have an integral model $R_\alpha$, then the
compatibility of the action $\sigma_n$ with the maps of the injective system
implies that we have an induced action on $R=\varinjlim R_\alpha$. Now the
noncommutative algebra for the integral model can be defined by adding
generators $\mu_n^*$ and $\tilde \mu_n$ with the relations \eqref{relstildemu}
and 
$$ \mu_n^* \, a = \sigma_n(a) \, \mu_n^* \ \ \ \text{ and } \ \ \  a\, \tilde\mu_n =\tilde\mu_n \sigma_n(a), $$
for all $a\in R$ and all $n\in \N$. This is an integral model for the algebraic endomotive algebra
$A\rtimes \N$.
\endproof

It would be interesting to investigate more closely what this relation between $\Lambda$-rings
and endomotives can say in the original case of $\Lambda$-rings in the context of Chern classes
as introduced by Grothendieck in \cite{Groth}. 

\smallskip

An interesting characterization is given in \cite{BoSm} of when a $\Lambda$-ring
defined over $\Q$ has an integral model in the sense of Definition \ref{LambdaQalg}.
Namely, this is the case if and only if the action of $G\times \N$ on $\cX=\Hom(A,\bar\Q)$ 
factors through a (multiplicative) action of $\hat\Z$. In the case of the BC endomotive, this
property can be seen explicitly as follows.

\begin{lem}\label{hatZact}
The multiplicative action of $\hat\Z$ on $\Hom(\Z[\Q/\Z],\bar\Q)$ is given by the
symmetries of the BC system, together with the action of $\N$ by 
the endomorphisms $\sigma_n$.
\end{lem}

\proof For $\alpha \in \hat\Z$ we set
\begin{equation}\label{alphaer}
\alpha: e(r) \mapsto e(\alpha(r)), \ \ \ \forall r\in \Q/\Z,
\end{equation}
where we identify $\hat\Z=\Hom(\Q/\Z,\Q/\Z)$. This gives a multiplicative action 
of $\hat\Z$ on the abelian part of the BC endomotive, $\Q[\Q/\Z]$ and its integral
model $\Z[\Q/\Z]$. For $\alpha = n\in \N$ this action reduces to the action of the
endomorphisms $\sigma_n$ described above, while for $\alpha\in \hat\Z^*$
it agrees with the action of the automorphisms group of the Bost--Connes system.
This action of $\hat\Z^*$ corresponds to the Galois action of $G=\Gal(\bar\Q/\Q)$
as the latter factorizes through the action of $\Gal(\Q^{ab}/\Q)$ given by the
cyclotomic character, \ie through the action of $\hat\Z^*$ of the form \eqref{alphaer},
under the class field theory isomorphism $\Gal(\Q^{ab}/\Q)\simeq \hat\Z^*$.
\endproof

Thus, the action of $\hat\Z$ is given by the
action of the symmetries $\hat\Z^*$ of the BC algebra (which preserve the
integral model) combined with the action of the endomorphisms $\N$.  
The following observation then links the notion of integral structure of a
$\Lambda$-ring to Haran's proposal \cite{Haran} for a notion of Frobenius over 
the field with one element.

\begin{cor}\label{HaranFrob}
For the BC endomotive, the action of $\hat\Z$ on $\Hom(\Z[\Q/\Z],\bar\Q)$ 
that determines the integral structure on the inductive limit of $\Lambda$-rings 
agrees with the ``Frobenius over $\bF_{1^\infty}$" of \cite{Haran}.
\end{cor}

\proof As observed in \S 6.1 of \cite{CCM2}, the action \eqref{alphaer} corresponds to the
notion of Frobenius over $\bF_{1^\infty}$ as defined in \cite{Haran} as follows.
Given a set with a free action of roots of unity (regarded as a ``vector space
over $\bF_{1^\infty}$), and an element $\alpha\in \hat\Z$, one obtains a new 
action on the same set by 
\begin{equation}\label{zetalpha}
 \zeta: x \mapsto \zeta^\alpha x,
\end{equation} 
which corresponds in fact to the action \eqref{alphaer}.
\endproof

In order to relate the multivariable BC endomotives to $\Lambda$-rings,
we first observe that, as in the one-variable case, these admit integral
models constructed as in \cite{CCM2}. Namely, 

the integer model of the multivariable BC endomotive algebra $\cA_{BC,n}$ is the
ring $\cA_{\Z,BC,n}$ generated by the elements
$e(\underline{r})$ of $\Z[\Q/\Z]^{\otimes n}$ and by generators $\tilde\mu_\alpha$ and
$\mu_\alpha^*$, with the relations
\begin{equation}\label{tildemualpharels}
\begin{array}{ll}
\tilde\mu_\alpha \tilde\mu_\beta = \tilde\mu_{\alpha\beta}, & \forall \alpha,\beta \in M_n(\Z)^+ \\[2mm]
\mu_\alpha^* \mu_\beta^* = \mu_{\beta\alpha}^*, & \forall \alpha,\beta \in M_n(\Z)^+ \\[2mm]
\mu_\alpha^* \tilde\mu_\alpha = \det(\alpha) & \forall \alpha \in M_n(\Z)^+ \\[2mm]
\mu_\alpha^* \tilde\mu_\beta = \det(\alpha) \tilde\mu_\gamma & 
\text{for } \, \beta=\alpha\gamma \in M_n(\Z)^+ .
\end{array}
\end{equation}
replacing the relations \eqref{mualpharels}, and with the additional relations
\begin{equation}\label{musigmarelalpha}
\mu_\alpha^* \, x = \sigma_\alpha(x)\, \mu_\alpha^*, \ \ \ \text{ and } \ \ \  x\, \tilde\mu_\alpha =
\tilde\mu_\alpha \, \sigma_\alpha(x).
\end{equation}
for all $x\in \Z[\Q/\Z]^{\otimes n}$ and all $\alpha\in M_n(\Z)^+$. As in the 
one-variable case, the $\tilde\mu_\alpha$
are no longer adjoints to the $\mu_\alpha^*$: they satisfy 
$\tilde\mu_\alpha=\det(\alpha) \mu_\alpha$.

\smallskip

We can then prove that the results of \cite{BoSm}, when rephrased in our setting, 
show that the BC endomotive and its multivariable analogs considered in
the previous section are universal with respect to the (torsion free)
$\Lambda$-rings of Definition \ref{Lambdaring} that are finite over $\Z$.  

\begin{thm}\label{BCunivLambda}
Let $A=\varinjlim_\alpha A_\alpha$ be a direct limit of $\Lambda$-rings, with
an integral model $R=\varinjlim_\alpha R_\alpha$, where the $R_\alpha$
have finite rank as abelian groups an have no non-zero nilpotent elements. 
Then the associated endomotive $A$ with the semigroup action of $\N$
that gives the $\Lambda$-ring structure embeds in a product of copies of
the BC endomotive $\Q[\Q/\Z]$ compatibly with the action of 
$S_{n,diag}\subset M_n(\Z)^+$.
\end{thm}

\proof Rephrased in our setting, Corollary 0.3 of \cite{BoSm} shows 
that every torsion free finite rank $\Lambda$ ring embeds in a $\Lambda$-ring that
is a product of copies of $\Z[\Q/\Z]$, which is the integral model of the
abelian part of a multivariable BC endomotive.  The fact that this is an 
embedding of $\Lambda$-rings implies that the semigroup action of $\N$ on $R$
is compatible with the action of the semigroup $S_{n,diag}\subset M_n(\Z)^+$
of the multivariable BC endomotive, so that one obtains an induced
embedding of the crossed product algebras compatible with the integer model
$\cA_{\Z,BC,n}$ described above.
\endproof

This result shows how the (multivariable) BC endomotives bridge between the
Soul\'e approach to varieties over $\F_1$, with respect to which the BC endomotive
embodies the tower of $\F_{1^n}$ as in \cite{CCM2}, and the approach of Borger
via $\Lambda$-rings, where they provide universal $\Lambda$-rings in the sense of
Theorem \ref{BCunivLambda} above.

\medskip

The examples considered in this paper of crossed product constructions both in
the case of Habiro rings and of the original BC system and its multivariable versions
suggests that, more generally, one should be able to associate an endomotive to
a $\Lambda$-ring, in the form of a crossed product, which we loosely write here
as $A\rtimes \N$.  While we do not give here a general construction for such crossed
products, one can see that the statement of Theorem \ref{BCunivLambda} should
then give a natural map, in the case of the crossed product of a $\Lambda$-ring that
is finite rank as an abelian groups an has no non-zero nilpotent elements, from the
crossed product algebra $A\rtimes \N$ to a multivariable BC algebra, compatibly with
the integral models.

\medskip

As recalled in \S 4.1 of \cite{Man} (see also \cite{Bor2} and \cite{BoWie}, \S 2) 
the affine ring scheme of big Witt vectors is defined by considering
the polynomial ring $\Z[\lambda_1,\lambda_2,\ldots,\lambda_n,\ldots]$ 
in the elementary symmetric functions $\lambda_1 =x_1+x_2+\cdots$, $\lambda_2=x_1x_2
+x_1x_3+x_2x_3+\cdots$, etc. This can be equivalently written in the set of generators
$\Z[ u_1, u_2,\ldots, u_n, \ldots]$
\begin{equation}\label{ghost}
\psi_n =\sum_{d|n} d\, u_d^{n/d},
\end{equation}
where the ghost coordinates are $\psi_n = x_1^n + x_2^n +\cdots$ 
and the $u_n$ are determined by \eqref{ghost}. For a commutative ring $R$,
the ring of big Witt vectors is $W(R)=\prod_{n\geq 1} R$, with the addition and
multiplication given by componentwise addition and multipication in the
ghost coordinates.  The truncated Witt scheme $W^{(N)}$ is obtained by
considering the subring $\Z[u_1,\ldots,u_N]$. As shown in \S 4.1 of \cite{Man}
the $W^{(N)}$ define gadgets $\cW^{(N)}_\Z$ over $\F_1$ in the sense of Soul\'e \cite{Soule},
by setting $\cW^{(N)}_\Z(R)$ to be the points of $W^{(N)}(R)$ with ghost coordinates
that are either zero or roots of unity. The $\Lambda$-ring structure is given by the
Adams operations $\psi_n$, which by \eqref{ghost} satisfy the requirement \eqref{liftFrob}.
Thus, in particular, one can use the result of Theorem \ref{BCunivLambda} to relate the
gadgets $\cW^{(N)}_\Z$ to the gadgets over $\bF_1$ associated to the multivariable BC
endomotives in the same way as done in \cite{CCM2} for the original BC algebra, with
compatible $\Lambda$-ring structure. This may provide the compatibility between the
idea of realizing the missing $\Sp(\Z)\times_{\bF_1}\cdots\times_{\bF_1} \Sp(\Z)$ in
terms of the quantum statistical mechanics of multivariable BC endomotives, as sketched
in \S \ref{KuroSec} above and the idea of \cite{Bor} to realize them via the spaces of
big Witt vectors. 

\smallskip

In this light, the multivariable BC endomotives constructed here should also be regarded
as basic building blocks for a more general class of noncommutative spaces with quantum 
statistical mechanical properties that can be constructed from the data of toric varieties and
possibly of other varieties over $\bF_1$ in the sense of \cite{CoCons}, 
\cite{CoCons2}, \cite{LL}. We intend to return to this topic in future work.

\section{Quantum channels and endomotives}\label{QcompSec}

This section contains a quick side remark on the general quantum
statistical mechanical formalism of endomotives, aimed at showing that, while one usually
concentrates on studying equilibrium (KMS) states for quantum statistical
mechanical systems such as the BC system and its generalizations, one can
also get some interesting information by considering states that are not
equilibrium states under the dynamics, but which are defined by particular
density operators. Given endomotives consisting of an abelian
part of the algebra and a semigroup of isometries acting on it, it is interesting
to see how density matrices of states transform under the semigroup action.

\smallskip

In quantum computing a formalism that describes the phenomenon of decoherence, the evolution
of pure states into mixed states, is obtained by replacing the group action by unitaries that
gives the evolution of density matrices by an evolution given not by a semigroup action by
isometries, which is given in the operator-sum representation or Kraus representation in the
form $\rho \mapsto \sum_{i=1}^N \mu_i \rho \mu_i^*$, where the $\mu_i$ are isometries satisfying
$\sum_i \mu_i^* \mu_i =1$ (see \eg \cite{Preskill}). 
The case of a single $\mu$ with $\mu^*\mu=\mu\mu^*=1$ gives 
back the unitary evolution. In the terminology of quantum computing
such data are said to describe a quantum channel. We show here how a similar
formalism can be applied in the context of endomotives using the semigroup action
and the induced action on states to obtain a transformation of the convex space of
density operators. It is important to distinguish the finite dimensional case, where
a quantum channel realized by a single isometry would have to be given by a unitary, 
and the infinite dimensional case we are considering here, where it is possible to have
a single isometry with $\mu^*\mu=1$ but with $\mu\mu^* =e= e^2\neq 1$ a projector. We 
consider here the isometries in the family $\mu_s$ with $s\in S$ a semigroup acting
by endomorphisms on an abelian algebra of observables and we consider the
quantum channels $\rho \mapsto \mu_s \rho \mu_s^*$ determined by these 
isometries. 

\smallskip

Given a semigroup $S$ acting by endomorphisms on the abelian part $A$ of an endomotive
in the sense of \cite{CCM}, consider the action of endomorphisms on states given as in
\cite{CM1} by
\begin{equation}\label{sphiact}
(s^*\varphi) (a) = \frac{\varphi(s(a))}{\varphi(s(1))},
\end{equation}
where $s(a)=\mu_s \, a \, \mu_s^*$ and $\mu_s \mu_s^* =e_s$ is a non-zero idempotent in the
algebra $A$, $e_s^2=e_s=e_s^*$, such that $\varphi(e_s)\neq 0$. 
Suppose that we consider states on $A$ of the form
\begin{equation}\label{phirho}
\varphi(a) = \frac{\Tr(a\varrho)}{\Tr(\varrho)},
\end{equation}
where $\varrho$ is a positive density operator, \ie a positive trace class operator.
Then we have
$$ \frac{\varphi(\mu_s \, a \, \mu_s^* \varrho)}{\varphi(e_s\, \varrho)} =\frac{\Tr(a\,\mu_s^* \, \varrho\, \mu_s)}
{\Tr(\varrho)} \frac{\Tr(\varrho)}{\Tr(\mu_s^* \varrho \mu_s)}. $$
This means that we have
$$ (s^*\varphi)(a)= \frac{\Tr( a\varrho_s)}{\Tr(\varrho_s)}, \ \ \ \text{ with } 
\varrho_s = \mu_s^* \varrho  \mu_s. $$ 
Notice that, while the action $a\mapsto \mu_s^* a \mu_s$ is an algebra homomorphism
only on the compressed algebra $a=e_s x e_s$, with $x\in A$ and $e_s=\mu_s \mu_s^*$ 
the idempotent that gives the projection onto the range of $\mu_s$, so that 
$(\mu_s^* a \mu_s) (\mu_s^* b \mu_s)=
\mu_s^* a e_s b \mu_s = \mu_s^* ab \mu_s$, for $a=e_s x e_s$ and $b=e_s y e_s$ with $x,y\in A$,
using $e_s^2=e_s$. However, when considering an action on density operators, one does
not need the compatibility with the multiplicative structure, as the only operation one performs 
on density operators is that of taking convex linear combinations.  Thus, one obtains a 
semigroup action of $S$ on the convex space of density operators by
\begin{equation}\label{musrho}
\varrho \mapsto \mu_s^* \varrho \mu_s .
\end{equation}
We now sketch some examples of how this point of view may be useful. 

Recall that a rational convex cone $\cC$ is a convex subset of $\R^n$ is the span
\begin{equation}\label{cone}
\cC=\R_+ v_1 +\cdots + \R_+ v_r, \ \ \ \text{ with } \ \ v_i \in \Q^n, \, i=1,\ldots, r.
\end{equation}
The cone is simplicial if the $v_i$ are linearly independent. Terasoma \cite{Tera} 
constructed multiple zeta values for rational convex cones by setting
\begin{equation}\label{TeraMZVcone}
\zeta_\cC (\ell_1,\ldots, \ell_k, \chi) = \sum_{v\in \cC^0 \cap \Z^n} \frac{\chi(v)}{\ell_1(v)\cdots \ell_k(v)},
\end{equation}
where the $\ell_i$ are $\Q$-linear forms on $\Q^n$ that are strictly positive on the interior $\cC^0$
of the cone $\cC$, $\ell_i(v)>0$ for all $v\in \cC^0$, and $\chi\in \Hom(\Z^n,\C^*)$ is a character 
of $\Z^n$. 

Let then $A=\Q[\Q/\Z]^{\otimes n}$
be the abelian part of one of the multivariable Bost--Connes systems considered earlier in this paper.
Let $\cC$ be a rational convex cone in $\R^n$ with the property that $\cC^0\cap \Z^n$ is preserved
by the action of a subsemigroup $S_\cC \subset M_n(\Z)^+$.
We consider special cases of states of the form \eqref{phirho} where we set
\begin{equation}\label{ConeMZVstate}
\Tr(a\varrho) =\sum_{v\in \cC^0 \cap \Z^n} 
\frac{\chi_a(v)}{\ell_1(v)\cdots \ell_k(v)},
\end{equation}
where the density matrix $\varrho$ is defined as the operator on the Hilbert space 
$\ell^2(\cC^0\cap \Z^n)$ given by
$$ \varrho \, \epsilon_v = \frac{1}{\ell_1(v)\cdots \ell_k(v)} \epsilon_v, $$
and the character $\chi_a$ is defined by the choice of an element $\alpha \in \GL_n(\hat\Z)$
that maps $(\Q/\Z)^n$ to $\cZ^n$ and then defining, for $a=e(r_1)\otimes\cdots\otimes e(r_n)$,
$$ \chi_a(v) = \zeta_1^{k_1}\cdots \zeta_n^{k_n}, $$
for $v=(k_1,\ldots, k_n)\in \cC^0\cap \Z^n$ and $\zeta_i=\alpha(e(r_i))\in \cZ$.

\smallskip

Under the action of \eqref{musrho}, for $s\in S_{\cC}$, we have 
$$ (\ell_1, \ldots, \ell_k) \mapsto (\ell_1 \circ s, \ldots, \ell_k\circ s) $$
so that
$$ \Tr(a\varrho_s) =\sum_{v\in \cC^0 \cap \Z^n} 
\frac{\chi_a(v)}{\ell_1(s(v))\cdots \ell_k(s(v))}. $$

\smallskip

A quantum channels $\varrho \mapsto \varrho_s$ thus
gives a ``change of variables'' in the multiple zeta 
values of cones, which preserves relations, in the sense that
if a polynomial relation of the form
\begin{equation}\label{polyrel}
\cR(\zeta_\cC(\ell_{1,i},\ldots,\ell_{k_i,i}, \chi_{a_i}))=\sum_I \alpha_I \prod_{i\in I} \Tr(a_i \varrho_i)=0
\end{equation}
holds for some $a_i\in (\Q/\Z)^n$, and if $a_i = s_i(b_i)$
for some $s_i \in \cS_\cC$ then 
\begin{equation}\label{polyrel2}
\cR(\zeta_\cC(\ell_{1,i}\circ s_i,\ldots,\ell_{k_i,i}\circ s_i, \chi_{b_i}))=0.
\end{equation}
This follows immediately from the fact that, for
$a_i=s_i(b_i)=\mu_{s_i}^* b_i \mu_{s_i}$, we have
$\Tr(a_i \varrho_i)=\Tr(b_i \varrho_{s_i,i})$ with 
$\varrho_{s_i,i} = \mu_{s_i}^* \varrho_i \mu_{s_i}$.

\section{Homology 3-spheres}\label{WRTSec}

We now discuss briefly how Habiro's universal Witten--Reshetikhin--Turaev 
(WRT) invariant may be used to lift some of the results described in the
previous sections to the world of 3-manifolds. This section will mostly contain
some questions about 3-manifolds motivated by the circle of ideas considered in
the previous sections. 

The idea of relating certain categories of 3-manifolds to 
convolution algebras and quantum statistical mechanical systems
was developed, in a different context, in \cite{MaZa}, where instead 
of relying on the description of 3-manifolds in terms of surgeries on 
links in the 3-sphere, which one uses to define WRT invariants,
one works with the description as branched coverings of the 3-sphere. 
While the considerations we outline here are different from the 
setting of \cite{MaZa}, the underlying principle and philosophy 
are similar.

\subsection{Habiro's universal WRT invariant}

The WRT invariants of 3-manifolds were originally introduced by Witten \cite{Wi}
in terms of Chern--Simons path integrals, and then defined rigorously by
Reshetikhin and Turaev \cite{RT} using representations of quantum groups 
at roots of unity. Restricting to the case of 3-manifolds that are integral homology 
3-spheres, one has a family of WRT invariants parameterized by roots of unity
\begin{equation}\label{WRTzeta}
\tau(M): \cZ \to \C, \ \ \zeta \mapsto \tau_\zeta(M),
\end{equation}
with consistency conditions, such as Galois equivariance 
$$ \tau_{\alpha(\zeta)}(M)=\alpha(\tau_\zeta(M)), \ \ \  
\forall \alpha \in \Gal(\Q^{ab}/\Q), \ \forall \zeta\in \cZ . $$
Moreover, the Ohtsuki series provides an expansion near $\zeta =1$
of the WRT invariant, with coefficients that in turn are invariants of
3-manifolds,
\begin{equation}\label{Oseries}
\tau^O(M)= 1 + \sum_{n=1}^\infty \lambda_n(M) (q-1)^n.
\end{equation}

The recent work of Habiro \cite{Hab2} provides a universal formulation for
the WRT invariants, in terms of a single functions $J_M(q)$ in the Habiro ring
\begin{equation}\label{HabJq}
J_M(q) \in \widehat{\Z[q]},
\end{equation}
with the properties that its evaluations at roots of unity recover the usual WRT
invariants for 3-dimensional integral homology spheres,
\begin{equation}\label{Jandtau}
ev_\zeta(J_M) = \tau_\zeta(M),
\end{equation}
and such that one has an analog of the Ohtsuki series at every root of unity,
given by the Taylor expansion of Habiro functions,  which recovers the
Ohtsuki series at $\zeta =1$,
\begin{equation}\label{JandO}
\ft_1(J_M)=\tau^O(M),
\end{equation}
where $\ft_\zeta(f)$ denotes, as before, the Taylor expansion of $f\in \widehat{\Z[q]}$
at $\zeta\in \cZ$.

\smallskip

The universal WRT invariant is obtained in \cite{Hab2} by defining the invariant
for an algebraically split link $L$ in $S^3$ with framing $\pm 1$. Then one uses 
the fact that any integral homology $3$-sphere has a presentation as surgery on
an algebraically split link in the 3-sphere with framing $\pm 1$, where two such
presentations give rise to the same homology 3-spheres if and only is they differ
by Fenn--Rourke moves. One then shows the invariance under Fenn--Rourke
moves to obtain a well defined $J_M$.

\smallskip

Let $\Z HS$ denote the free abelian group generated by the orientation-preserving 
homeomorphism classes of integral homology 3-spheres. This is a ring with the 
product given by the connected sum $M_1 \# M_2$.

\smallskip

The universal WRT invariant defines a ring homomorphism
\begin{equation}\label{WRThom}
J: \Z HS \to \widehat{\Z[q]},
\end{equation}
since it satisfies the properties
\begin{equation}\label{WRThom2}
J_{M_1\# M_2}(q) = J_{M_1}(q) J_{M_2}(q), \ \ \  J_{S^3}(q)=1, \ \ \  J_{-M}(q) =J_M(q^{-1}).
\end{equation}

\smallskip

Moreover, to express the fact that the WRT invariants have some ``finite type" properties,
one can consider, as in \cite{Hab2}, the Ohtsuki filtration on $\Z HS$,
\begin{equation}\label{Ofiltr}
\Z HS = F_0 \supset F_1 \supset \cdots \supset F_k \supset \cdots,
\end{equation} 
where $F_k$ is the $\Z$-submodule spanned by the alternating sums
\begin{equation}\label{spanFk}
[M, L_1, \ldots, L_k] = \sum_{L'\subset\{ L_1,\ldots, L_k\}} (-1)^{|L'|} M_{L'},
\end{equation}
with $|L|$ the number of components and the $L_i$ algebraically split links
of framing $\pm 1$. This is generalized in \cite{Hab2} to a filtration by 
submodules $F_d$, for $d: \N \to \Z_{\geq 0}$ finitely supported, where now 
$L=\{ L_1,\ldots L_k\}$ is an algebraically split link with framings in the set
$\{ \pm 1/m\,|\, m\in \Z \smallsetminus \{ 0 \}\}$ and in the expressions 
\eqref{spanFk} there are $d(n)$ components of the link that
have framing with $m=\pm n$.

\smallskip

Habiro conjectures in \cite{Hab2} that the universal WRT invariant will induce a
ring homomorphism
\begin{equation}\label{HabConj}
J : \widehat{\Z HS}\to \widehat{\Z[q]},
\end{equation}
where 
\begin{equation}\label{hatZHS}
\widehat{\Z HS} = \varprojlim_d \Z HS /F_d.
\end{equation}

\smallskip

Notice that the universal WRT invariant takes values in the single-variable Habiro
ring and one may wonder whether there are similar 3-manifolds or link invariants
that give rise naturally to functions in the multivariable Habiro rings of \cite{Man}.
The multivariable link invariants constructed in \cite{Geer} may provide such objects.

\subsection{Homology 3-spheres and the field with one element}

We recall the notion of {\em truc} or {\em gadget} over the
field with one element $\F_1$, as defined by Soul\'e in \cite{Soule}.

\begin{defn}\label{truc} 
A gadget over $\F_1$ is the datum $(X,\cA_X,e_{x,\sigma})$ 
of a covariant functor $X: \cR \to \cS$
from the category of commutative finite flat rings $\cR$ to the category 
$\cS$ of sets, a complex algebra $\cA_X$ and evaluation maps given by
$\C$-algebra homomorphisms $e_{x,\sigma}:\cA_X \to \C$,
for all $x\in X(R)$ and all homomorphisms $\sigma: R \to \C$, such
that, for any given ring homomorphism $f: R'\to R$,
\begin{equation}\label{exsigma}
 e_{f(y),\sigma} = e_{y,\sigma\circ f}. 
\end{equation}
\end{defn}

Affine varieties $V_\Z$ over $\Z$ give rise to associated gadgets $X=G(V_\Z)$ over
$\F_1$, by taking $X(R)=\Hom(O(V),R)$ and $\cA_X=O(V)\otimes \C$. 

According to the definition of Soul\'e \cite{Soule}, affine variety over $\F_1$ are
gadget over $\F_1$ where all the $X(R)$ are finite and for which there exists a
variety $X_\Z$ over $\Z$ and a morphism of gadgets
$$ X \to G(X_\Z) $$
such that all morphisms of gadgets $X \to G(V_\Z)$ come from morphisms of
varieties $X_\Z \to V_\Z$ over $\Z$. 

In particular, in the Soul\'e approach to $\F_1$ geometry, it is a consistency 
condition on the cyclotomic points of a scheme over $\Z$ that determines
whether it comes from an underlying $\F_1$ structure. This replaces the
consistency condition on lifts of Frobenius that is used in the $\Lambda$-ring
approach we discussed in \S \ref{LambdaSec} above. In particular, as
observed in \S 1.7 of \cite{Man}, one can probe the presence of an $\F_1$ structure
on a $\Z$-scheme, in the Soul\'e setting, by restricting to the case where the rings
$R$ are group rings of finite abelian groups, such as the $\Z[\Z/n\Z]$.

For example, as discussed in \cite{CCM2}, the affine varieties $\mu^{(n)}$ over 
$\F_1$ given by roots of unity assemble as a direct limit, with $\cA_X =C(S^1)$,
to give the multiplicative group $\bG_m$ over $\F_1$ (as in \cite{Soule}) or as
an inverse limit, with $\cA_X = C^*(\Q/\Z)$, the abelian part of the BC algebra,
to give the structure of (pro)variety over $\F_1$ to the BC endomotive.

A generalization of this notion of variety over $\F_1$ that covers certain important 
classes of non-affine cases was recently given in \cite{CoCons}.

\smallskip

We now show how to use integral homology 3-spheres and the universal
Witten--Reshetikhin--Turaev invariant to construct a gadget over $\F_1$
in the sense described above.

\begin{lem}\label{XZHS}
For $R\in \cR$ let $\Hom(\Z HS,R)$ denote the set of ring homomorphisms
and define
\begin{equation}\label{ZHSR}
X_{\Z HS}(R):= \{ \phi\in \Hom(\Z HS,R)\,|\, \exists \tilde\phi: 
\widehat{\Z[q]} \to R, \,\, \phi = \tilde\phi \circ J \}.
\end{equation}
Then $X_{\Z HS}$ defines a covariant functor $X_{\Z HS}: \cR \to \cS$.
\end{lem}

\proof Let $f: R' \to R$ be a ring homomorphism, for $R,R'\in \cR$.
Let $X_{\Z HS}(f) : X_{\Z HS}(R') \to X_{\Z HS}(R)$ be given by
$X_{\Z HS}(f)(\psi)=f \circ \psi$. If there exists $\tilde\psi:
\widehat{\Z[q]} \to R'$ such that $\psi =\tilde \psi  \circ J$,
then setting $\tilde \phi =f \circ \tilde \psi$ satisfies
$\phi = \tilde\phi \circ J$ for $\phi= f\circ \psi$. This suffices
to see that $X_{\Z HS}$ defines a covariant functor as described.
\endproof

The set $X_{\Z HS}(R)$ defined as in \eqref{ZHSR} describes the
set of all $R$-valued invariants of integral homology 3-spheres that
are coarser than the Witten--Reshetikhin--Turaev invariant, in the
sense that they factor through this invariant. This means that they
do not distinguish homology 3-spheres that are not already
distinguished by the WRT invariant.

\smallskip

We define the complex algebra $\cA_{\Z HS}$ as follows. Consider
the complex algebra $\Z HS\otimes \C$ and let $\cA_{\Z HS}$ be the 
$C^*$-completion with respect to the norm $\| M \| = \| J_M \|_{\cZ}$,
where the latter is the norm defined as in Lemma \ref{ZzetaRep},
using the evaluations $ev_\zeta(J_M) \in \Z[\zeta]$ viewed as complex numbers
after embedding $\cZ\subset \C$. The supremum is taken over $\cZ$ 
identified with $\hat\Z$, so that it is a compact set. 
With the involution defined by $M \mapsto -M$, the algebra $\cA_{\Z HS}$
is an abelian $C^*$-algebra. 

\begin{lem}\label{evxsigmaZHS}
Given $\phi \in X_{\Z HS}(R)$ and $\sigma:R \to \C$, then there
exists a unique $C^*$-algebra homomorphism
\begin{equation}\label{evZHS}
e_{\phi,\sigma}: \cA_{\Z HS} \to \C, 
\end{equation}
that extends $\sigma\circ \phi$ to $\Z HS\otimes \C$ and to $\cA_{\Z HS}$.
This satisfies the property \eqref{exsigma}.
\end{lem}

\proof Suppose given a ring hoomorphism $\sigma: R \to \C$. Then
$\sigma \circ \phi$ determines a homomorphism $\Xi =\sigma \circ \tilde\phi :
\widehat{\Z[q]}\to \C$. Notice that any homomorphism $\Xi:  \widehat{\Z[q]}\to \C$
factors through the evaluations at roots of unity. Namely, assigning such a 
homomorphism is equivalent to assigning a family of homomorphisms 
$\Xi_m : \Z[q]/((q)_m) \to \C$, compatible with the maps of the projective system.
A homomorphism $\Xi_m :\Z[q] \to \C$ such that $\Xi_m((q)_m)=0$ has the
property that $\Xi_m(q) =\zeta \in \C$ is a root of unity, since
$$ \Xi_m ((1-q)\cdots (1-q^m))= (1-\zeta) \cdots (1-\zeta^m) =0. $$
Thus,  $\Xi_m$ factors through $ev_\zeta$ and so, by the compatibility of the
$\Xi_m$ with the maps of the projective system, does $\Xi$. One can then define
$e_{\phi,\sigma}: \cA_{\Z HS} \to \C$ to be given by $e_{\phi,\sigma}(M)=ev_\zeta(J_M)$. 
By the choice of the norm $\| M \| = \| J_M \|_{\cZ}$ on $\cA_{\Z HS}$ we see that
these are continuous maps. Moreover, they satisfy by construction the property 
\eqref{exsigma}.
\endproof

This means that one uses the WRT invariants as cyclotomic coordinates on the
gadget $X_{\Z HS}$. The results of Lemma \ref{XZHS}  and \ref{evxsigmaZHS} 
directly imply the following.

\begin{prop}\label{WRTgadget}
The data $(X_{\Z HS}, \cA_{\Z HS}, e_{\phi,\sigma})$ define a gadget over $\F_1$.
\end{prop}

It seems then reasonable to ask whether, as in the Habiro conjecture recalled
above, the analogous data $(X_{\Z HS/F_d}, \cA_{\Z HS/F_d},e_{\phi,\sigma})$ 
constructed using the quotients of the Ohtsuki filtration, can be used to define
a directed system of {\em affine varieties} over $\F_1$, which gives rise to
a direct limit $(X_{\widehat{\Z HS}}, \cA_{\widehat{\Z HS}}, \hat e_{\phi,\sigma})$
as varieties over $\F_1$. 

\subsection{Semigroup actions on homology 3-spheres}\label{semiZhsSec}

It is also natural in this context to look for semigroup actions of the
multiplicative semigroup $\N$ of positive integers on the ring $\Z HS$
of integral homology 3-spheres. The question is whether one can
produce an endomotive, constructed out of the ring $\Z HS$ (or of
the $\Z HS/ F_d$ and $\widehat{\Z HS}$) in such a way that a suitable
invariant $\alpha: \Z HS \to \widehat{\Z[q]}$, obtained from the universal
WRT invariant, gives an algebra homomorphism of the semigroup
crossed product to the Bost--Connes algebra.

\smallskip

The construction of the universal WRT invariants given in  \cite{Hab2}
is based on the description of homology 3-spheres via a special kind
of surgery presentations, as in \cite{Hab3}.
Using a surgery presentation of integral homology 3-spheres as
surgeries $M=S^3_{(L,\fr)}$, with $L=L_1 \cup \cdots \cup L_\ell$ 
an algebraically split link in $S^3$ and framings $\fr=(1/m_i)_{i=1,\ldots, \ell}$, $m_i\in \Z$.

\begin{center}
\begin{figure}
\includegraphics[scale=0.6]{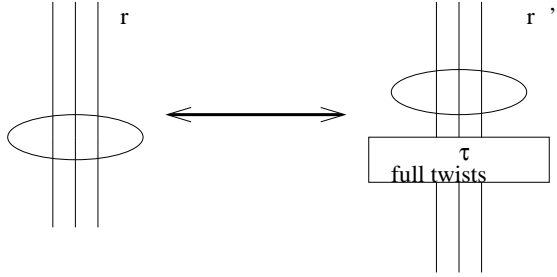}
\caption{Rolfsen move for links with rational framings. \label{Rolf2fig}}
\end{figure}
\end{center}

In the case of integer framings, it is well known that different surgery 
presentations that give rise to the same  3-manifold up to an orientation
preserving homeomorphism differ by Fenn--Rourke moves or 
equivalently by the two Kirby moves (see \S VI.19 of \cite{PraSo}). 
In the case of rational instead of
integer framings, the analog of these Kirby moves are the Rolfsen moves
described in \cite{Rolf}. The main result of \cite{Rolf} states that two surgery
presentations with links $L$ and $L'$ and rational framings $r$ and $r'$
give rise to orientation-preserving homeomorphic 3-manifolds if and only
if the surgery data are related by a sequence of moves of the following types
(Rolfsen moves):
\begin{itemize}
\item Introduce or delete a component of the link with framing $\infty$.
\item Replace a linked component $L_1$ with framing $r_1$ with a 
number $\tau$ of full twists of the remaining components with new
framings
\begin{equation}\label{newframesR}
\left\{\begin{array}{rl}
r_1' = & \frac{1}{ r_1^{-1} +\tau } \\[3mm]
r_i' = & r_i + \tau \ell(L_1,L_i)^2, \ \ \  i\neq 1,
\end{array}\right.
\end{equation}
where $\ell(L_1,L_i)$ is the linking number. This move is illustrated in Figure \ref{Rolf2fig}.
\end{itemize}

\begin{center}
\begin{figure}
\includegraphics[scale=0.5]{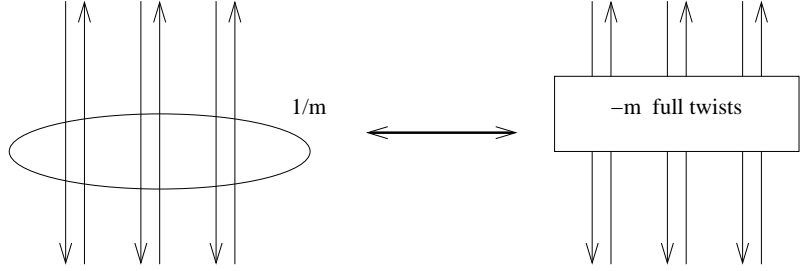}
\caption{Rational Hoste move between admissible links. \label{HosteFig}}
\end{figure}
\end{center}

It was conjectured by Hoste in \cite{Hoste}
that this should in fact be possible, that is, that two surgery presentations by
algebraically split links with framings in the set $\{ 1/m \,|\, m\in \Z, m\neq 0 \}$ give
rise to orientation-preserving homeomorphic homology 3-spheres if and only if
they are related by a sequence of Rolfsen moves {\em through links of this same
type}.  The conjecture was then proved by Habiro in \cite{Hab3}, Corollary 5.2. More
precisely, Habiro defines an {\em admissible} link to be an algebraically split
link with framings in the set $\{ 1/m \,|\, m\in \Z \}$ and defines a {\em rational 
Hoste move} as a Rolfsen move between a pair of admissible links (see Figure 
\ref{HosteFig}).  It is shown in Corollary 5.2 of \cite{Hab3} that two admissible links 
yield a pair of orientation-preserving diffeomorphic homology 3-spheres if and 
only if they are related by a sequence of rational Hoste moves. In these moves a
component $L_1$ with framing $1/m$ surrounding a number of other components 
$L_i$, arranged in such a way that $\ell(L_i,L_1)=0$ as in the figure, is replaced 
by $-m$ full twists of the remaining pairs of strands. The framings of the remaining
components are unchanged by \eqref{newframesR} with $\tau = -m$ and 
$\ell(L_i,L_1)=0$.

We can then formulate the following question, suggested by the 
constructions described in the previous sections in terms of endomotives 
and Habiro rings.

\begin{ques}\label{question}
Is there any natural (semi)group action
\begin{equation}\label{sigmanZHS}
M \mapsto \sigma_n(M), \ \ \  \sigma_n(M_1\# M_2) = \sigma_n(M_1)\# \sigma_n(M_2)
\end{equation}
of the multiplicative semigroup $\N$ or of the group $\Q^*_+$,
by symmetries of the ring $\Z HS$?
If so, is it then possible to obtain a suitable 
modification $\Upsilon_M$ of the WRT invariant $J_M$ such that
it is still a ring homomorphism $\Upsilon: \Z HS \to \widehat{\Z[q]}$,
or possibly a ring homomorphism 
$\Upsilon: \widehat{\Z HS} \to \widehat{\Z[q]}$, with the property 
\begin{equation}\label{Upsilon}
\Upsilon_{\sigma_n(M)} = \sigma_n(\Upsilon_M) \in  \widehat{\Z[q]}
\end{equation}
that it intertwines the semigroup actions? More generally, are there
semigroup actions (not necessarily by $\N$) on 3-manifolds that
induce corresponding actions on the WRT invariants?
\end{ques}

In the more general form stated above, one knows that there are
such examples, at least for suitable classes of 3-manifolds. For instance,
in the case of mapping cylinders and Seifert 3-manifolds (see \cite{Andersen},
\cite{Hansen}, \cite{LaRo}) an action of the mapping class group on the
surface translates into a linear action on the corresponding space of
invariants. In the more restrictive form of the question formulated above,
that is, whether one can construct an action of the multiplicative semigroup 
$\N$, and whether this can be done in a way related to the action on the
Habiro ring, the question is more complicated to address.

\medskip

More precisely, the problem of defining an action $M \mapsto \sigma_n(M)$ 
can be formulated as the question of defining an action on link diagrams,
which is well defined under the equivalence via Reidemeister moves
and in addition preserves the Hoste moves, so that it induces a well
defined action on the homology 3-spheres. 

In view of the result of \cite{Hab2}, \S 15.3, on the universal WRT 
invariant of $1/m$ surgery on a knot in a homology 3-sphere,
and of the effect of changes of framings on the WRT invariants as
analyzed in \cite{KiMe}, it would seem that an ingredient
in the construction of a semigroup action of the form \eqref{sigmanZHS}
may be given by a multiplicative change of the framings, that is,
a semigroup action on the set of admissible framed links of the form
$(L,\fr) \mapsto (L, n^{-1} \fr)$, 
where $L=L_1\cup\cdots, L_\ell$ with framings $\fr=(1/m_1,\ldots, 1/m_\ell)$,
$m_i\in \Z$, $m_i\neq 0$, and $n^{-1} \fr =(1/(m_1n) ,\ldots, 1/(m_\ell n))$. 
However, while changing the framings multiplicatively as in 
gives a well defined action on the set of admissible links, the Hoste move 
shows that this does not translate directly into a well defined semigroup 
action on the set of homology 3-spheres
unless one operates also on the full twists. In fact, different surgery presentations
for the same homology 3-sphere can replace a component with $1/m$ framings
by $-m$ full twists of the pairs of strands linked by the components, and the
operations can only become independent of the surgery presentation if the
number of full twists is also changed multiplicatively to $-mn$ full twists. However,
the number of twists is a property of a planar projection of a link and can be
changed by Reidemeister moves: by a Reidemeister move one can change a
crossing into a full twist followed by a crossing in the opposite direction. (This
simple observation leads to the description of crossings changes via Hoste moves.) 

\medskip

One can, in fact, approach the problem of Question \ref{question} 
in a slightly different way. In the formulation we gave of Question 
\ref{question} we considered ring homomorphisms $\Z HS \to \fR$,
in particular with ring of values $\fR=\widehat{\Z[q]}$ the Habiro
ring. The generators of $\Z HS$ are homology 3-spheres up to oriented
homeomorphisms, or equivalently admissible links up to Reidemeister
and Hoste moves, and the multiplication operation is connected sums,
or equivalently disjoint unions of links. Instead of working with 
equivalence classes of admissible links, both under Reidemeister 
and Hoste moves, one can follow the usual method for handling 
equivalence relations in noncommutative geometry, which is to keep the equivalence
relations explicitly in the form of a groupoid or a category, instead of passing to the set
of equivalence classes. In this sense, the situation resembles more closely the
approach followed in \cite{MaZa} using 3-manifolds as coverings and then cobordisms
between them, assembled first in a 2-category and then into an associative algebra.

\medskip

In terms of admissible links with equivalences given by Reidemeister and Hoste
moves, one can proceed in the following way.

\begin{defn}\label{2catlinks}
Let $\cC\cL^{(2)}$ be the 2-category whose objects are admissible framed link diagrams 
$(D,\fr)$, the 1-morphisms are sequences of  Hoste moves $(D,\fr)\stackrel{H}{\to} (D',\fr')$
between diagrams, including the identity move, and 2-morphisms $H_1 \stackrel{\cR}{\to} H_2$
are sequences of pairs of Reidemeister moves $\cR=(R,R')$ between the resulting diagrams so that 
the squares commute
\begin{equation}\label{xypic1}
 \xymatrix{
(D_1,\fr_1) \rto^{H_1} \dto^{R} & (D'_1,\fr_1') \dto^{R'} \\
(D_2,\fr_2) \rto^{H_2}  & (D'_2,\fr_2')
} \end{equation}
There is an associated convolution algebra of functions $\cA_{\cC\cL^{(2)}}$ given
by functions with compact support $f: Mor^{(2)}(\cC\cL^{(2)}) \to \C$ with
two different associative products corresponding, respectively, to the horizontal
and vertical composition of 2-morphisms,
\begin{equation}\label{prod1cat2}
(f_1\circ f_2)(\cR)=\sum_{\cR=\cR_1 \circ \cR_2} f_1(\cR_1) f_2(\cR_2),
\end{equation}
\begin{equation}\label{prod2cat2}
(f_1\bullet f_2)(\cR)=\sum_{\cR=\cR_1\bullet \cR_2} f_1(\cR_1) f_2(\cR_2),
\end{equation}
where the horizontal composition is given by $\cR_1=(R,R')$ and $\cR_2=(R',R'')$ with
$$ \xymatrix{ 
(D_1,\fr_1) \rto^{H_1} \dto^{R_1} & (D'_1,\fr_1') \dto^{R_1'} \rto^{H_1'} & (D''_1,\fr_1'') \dto^{R_1''} \\
(D_2,\fr_2) \rto^{H_2} & (D'_2,\fr_2') \rto^{H_2'} & (D''_2,\fr_2'') } $$ 
and the vertical composition is given by $\cR_1=(R_1,R_1')$ and $\cR_2=(R_2,R_2')$ with
$$ \xymatrix{
(D_1,\fr_1) \rto^{H_1} \dto^{R_1} & (D'_1,\fr_1') \dto^{R_1'} \\
(D_2,\fr_2) \rto^{H_2} \dto^{R_2} & (D'_2,\fr_2') \dto^{R_2'} \\
(D_3,\fr_3) \rto^{H_3} & (D'_3,\fr_3') }
$$
\end{defn}

The algebra $\cA_{\cC\cL^{(2)}}$ has an involution $f^* (\cR)= \overline{f(\cR^{-1})}$,
where $\cR^{-1}=(R'^{-1},R^{-1})$ is defined by the diagram 
$$ \xymatrix{ (D_2',\fr_2') \rto^{H_2^{-1}} \dto^{R'^{-1}} & (D_2,\fr_2) \dto^{R^{-1}} \\
(D_1',\fr_1') \rto^{H_1^{-1}} & (D_1,\fr_1) }
$$
This involution satisfies
$$ (f_1\circ f_2)^\vee = f_2^\vee \circ f_1^\vee, \ \ \ \text{ and } \ \ \  
(f_1\bullet f_2)^\vee = f_2^\vee \bullet f_1^\vee. $$

\medskip

One can pass from the 2-category $\cC\cL^{(2)}$ to a category $\cC\cL$ 
where one ``collapses the 2-morphisms", that is, one uses the 2-morphisms
to generate an equivalence relation on objects and 1-morphisms. The
category one obtains in this way is the following.

\begin{defn}\label{catlinks}
Let $\cC\cL$ be the category whose objects are admissible framed link types $(L,\fr)$ 
in $S^3$ and whose morphisms are sequences of Hoste moves between link types (that is,
Hoste moves up to Reidemeister moves), $(L,\fr)\stackrel{H}{\to} (L',\fr')$.  
The corresponding convolution algebra $\cA_{\cC\cL}$ is given by functions 
with compact support $f: Mor(\cC\cL) \to \C$ with the convolution product
\begin{equation}\label{algCL}
(f_1\circ f_2)(H)=\sum_{H = H_1 \circ H_2} f_1(H_1)
f_2(H_2),
\end{equation}
where the sum is over the compositions of arrows given by Hoste moves:
$$ (L,\fr)\stackrel{H}{\to}(L',\fr') = ((L,\fr)\stackrel{H_1}{\to} (\tilde L,\tilde\fr)) \circ ((\tilde L,\tilde\fr) \stackrel{H_2}{\to} (L',\fr')). $$
\end{defn}

\medskip

One can, for instance, construct a time evolution on the algebra $\cA_{\cC\cL}$  by setting
\begin{equation}\label{evolcomponents}
 \sigma_t(f)(H) = e^{i m t} f(H), \ \ \  \text{ with } \ \ \  m = \ell - \ell',
 \end{equation}
 if $H$ is a sequence of Hoste moves between links $L=L_1\cup\cdots\cup L_\ell$ and 
 $L'=L'_1\cup\cdots\cup L'_{\ell'}$.
 This manifestly satisfies $\sigma_f(f_1\circ f_2)=\sigma_t(f_1)\circ \sigma_t(f_2)$.

Let $\cS_{(L,\fr)}$ be the set of sequences of Hoste moves $H$ transforming some
link types $(\tilde L, \tilde \fr)$ into $(L,\fr)$. Let 
$\cH_{(L,\fr)}=\ell^2(\cS_{(L,\fr)})$ and define the action of $\cA_{\cC\cL}$ on $\cH_{(L,\fr)}$ 
by setting
\begin{equation}\label{piLACL}
(\pi_{(L,\fr)}(f) \xi)((\tilde L,\tilde\fr)\stackrel{H}{\to} (L,\fr))= \sum_{H= H_1 \circ H_2}
f((\tilde L,\tilde\fr)\stackrel{H_1}{\to} (L'',\fr'')) \,\, \xi((L'',\fr'') \stackrel{H_2}{\to} (L,\fr)).
\end{equation}

In this representation the time evolution \eqref{evolcomponents} is generated by a
Hamiltonian $\fH$ of the form
\begin{equation}\label{HamHoste}
\fH \delta_H = m\, \delta_H, \ \ \   H\in \cS_{(L,\fr)},
\end{equation}
where $(\tilde L, \tilde \fr) \stackrel{H}{\to} (L,\fr)$ with $m=\tilde\ell- \ell$ the difference
between the number of components of $\tilde L$ and $L$. Thus, $\Sp(\fH)=\Z$
and the eigenspace $E_m$ of $\fH$ is spanned by all the sequences of Hoste moves
that have final result $(L,\fr)$ and change the total number of components by $m$.
To see that $\fH$ is the Hamiltonian in the representation $\pi_{(L,\fr)}$, notice that we have
$$ (\pi_{(L,\fr)}(\sigma_t(f)) \xi)(H)= \sum_{H=H_1\circ  H_2} e^{i m_1 t} f(H_1) \xi(H_2) $$
which agrees with
$$ (e^{i\fH t} \pi_{(L,\fr)}(f) e^{-i\fH t} \xi)(H)= e^{i m t} \sum_{H=H_1 \circ H_2} 
f(H_1) e^{-i m_2 t} \xi(H_2). $$

In \cite{MaZa} we introduced the notion of horizontal and vertical time evolutions
on convolution algebras associated to 2-categories, \ie time evolutions with
respect to either the horizontal or the vertical convolution product. 
Since Reidemeister moves do not change the number of components, the time
evolution $\sigma_t$ of \eqref{evolcomponents} can be lifted to a horizontal
time evolution on the algebra $\cA_{\cC\cL^{(2)}}$. Vertical time evolutions can
be obtained similarly, by considering, instead of the number of components,
quantities such as the self-linking number with respect to the blackboard framing
of the link diagram, which is altered along Reidemeister moves.

\medskip

In view of passing from links and Hoste moves to 3-manifolds,
we can introduce a refinement of the algebras $\cA_{\cC\cL^{(2)}}$ 
and $\cA_{\cC\cL}$ considered above. In fact, the convolution 
products of these algebras are designed to account for the compositions
of Hoste and Reidemeister moves but not for the product
of links given by disjoint union, which in turn will determine the
product in the ring $\Z HS$ of integral homology 3-spheres given
by the connected sum.  One can introduced variants $\cA_{\cC\cL^{(2)},\cup}$ 
and $\cA_{\cC\cL,\cup}$ where the convolution product is modified to
include also the product operation of disjoint union. We explain in detail
how this can be done in the simpler case of $\cA_{\cC\cL,\cup}$, while the
case of $\cA_{\cC\cL^{(2)},\cup}$ is similar.

We define $\cA_{\cC\cL,\cup}$ to be the algebra of functions with finite support
$f: Mor(\cC\cL)\to \C$ with the modified product
\begin{equation}\label{prod2ACL}
(f_1\star f_2)(L\stackrel{H}{\to} L')=\sum_{H=(H_{1,1}\cup H_{1,2})(H_{2,1}\cup H_{2,2})} 
f_1(L_1 \stackrel{H_{1,1}}{\to} \tilde L_1) f_2(\tilde L_2 \stackrel{H_{2,2}}{\to} L'_2),
\end{equation}
where one sums over all possible splittings into connected components as well
as over compositions of Hoste moves, with
$$ L=L_1\cup L_2  \stackrel{H= H_1\cup H_2}{\longrightarrow} L'=L_1'\cup L_2'. $$

Notice in particular, that in the case where $H=id$ one obtains in this way the product
$$ (f_1\star f_2)(L)=\sum_{id=(H_{1,1}\cup H_{1,2})(H_{2,1}\cup H_{2,2})} f_1(H_1) f_2(H_2) , $$
and the sum splits into a term of the form
$\sum_{L=L_1\cup L_2} f_1(L_1) f_2(L_2)$, corresponding to the case where
all the terms have $H_{i,j}=id$, and a term containing nontrivial Hoste moves
and their inverses $H_{1,1}=H_{2,1}^{-1}$ and $H_{1,2}=H_{2,2}^{-1}$.
Thus, this can be seen as a modification of the group ring product
$$ (f_1 * f_2)(L)=\sum_{L=L_1\cup L_2} f_1(L_1) f_2(L_2) $$
of the group of links with the product $\cup$ given by disjoint union.

\begin{lem}\label{assprod2ACL}
The product \eqref{prod2ACL} is associative.
\end{lem}

\proof Consider a triple composition of morphisms $L_j \stackrel{H^{(j,j+1)}}{\to} L_{j+1}$,
$j=1,2,3$, with corresponding possible splittings into connected components of the form
$$ \cup_i L_{1,i} \stackrel{ \cup H^{(12)}_{i} }{\to} \cup_i L_{2,i} \stackrel{ \cup H^{(23)}_i }{\to}
\cup_i L_{3,i} \stackrel{ \cup_i H^{(34)}_i }{\to} \cup_i L_{4,i}, $$
for $i=1,2,3$. We then have
$$ (f_1\star(f_2\star f_3))(L_1\stackrel{H}{\to} L_4) = \sum_{H= H^{(12)}\circ H^{(24)}} 
f_1(H^{(12)}_1) \,\, (f_2\star f_3)(H^{(24)}_2), $$
where $H^{(24)}=H^{(23)}\circ H^{(34)}$ with
$L_{2,2}\cup L_{2,3} \stackrel{H^{(24)}_{23}}{\to} L_{4,2}\cup L_{4,3}$.
This gives then
$$  (f_1\star(f_2\star f_3))(L_1\stackrel{H}{\to} L_4) = \sum_{H= H^{(12)}\circ H^{(23)}\circ H^{(34)}}
f_1(H^{(12)}_1) \,\, f_2(H^{(23)}_2) \,\, f_3(H^{(34)}_3), $$
where
$$ L_{2,2}\cup L_{2,3} \stackrel{H^{(24)}_{23}}{\to} L_{4,2}\cup L_{4,3} =
L_{2,2}\cup L_{2,3} \stackrel{H^{(23)}_2 \cup H^{(23)}_3}{\to} L_{3,2}\cup L_{3,3} 
\stackrel{H^{(34)}_2 \cup H^{(34)}_3}{\to} L_{4,2}\cup L_{4,3}. $$
On the other hand, we have 
$$ ((f_1\star f_2)\star f_3)(L_1\stackrel{H}{\to} L_4) = \sum_{H=H^{(13)}\circ H^{(34)}}
(f_1\star f_2)(H^{(23)}_{12}) \,\, f_3(H^{(34)}_3) $$
with $L_{1,1}\cup L_{1,2} \stackrel{H^{(13)}_{12}}{\to} L_{3,1}\cup L_{3,2}$,
so that we have
$$ ((f_1\star f_2)\star f_3)(L_1\stackrel{H}{\to} L_4) = \sum_{H= H^{(12)}\circ H^{(23)}\circ H^{(34)}}
f_1(H^{(12)}_1)\,\, f_2(H^{(23)}_2)\,\, f_3(H^{(34)}_3), $$
with 
$$ L_{1,1}\cup L_{1,2} \stackrel{H^{(13)}_{12}}{\to} L_{3,1}\cup L_{3,2} =
L_{1,1}\cup L_{1,2} \stackrel{H^{(12)}_1 \cup H^{(12)}_2}{\to}
L_{2,1}\cup L_{2,2} \stackrel{H^{(23)}_1 \cup H^{(23)}_2}{\to} L_{3,1}\cup L_{3,2}. $$
Thus, we still have an associative product.
\endproof

One can similarly refine the algebra $\cA_{\cC\cL^{(2)}}$ to a version 
$\cA_{\cC\cL^{(2)},\cup}$ that incorporates the product of link
diagrams by disjoint union. Then, in this setting, an analog of 
Question \ref{question} can be formulated in terms of
homomorphisms from one of the algebras
$\cA_{\cC\cL}$, $\cA_{\cC\cL,\cup}$, $\cA_{\cC\cL^{(2)}}$, $\cA_{\cC\cL^{(2)},\cup}$  
to a suitable ring $\fR$ of values (\eg the Habiro ring), viewed 
as ``categorified" versions of invariants of homology 3-spheres usually defined 
as ring homomorphisms $\Z HS \to \fR$.

\smallskip

One can introduce other variants of the convolution algebras of
link diagrams, Hoste and Reidemeister moves. We discuss here
a variant that is related to the possibility of defining an action on 
links by a multiplicative change of framings. 

\smallskip

We recall that given a link $L$ with $\ell$
components and a vector $\n=(n_i)\in \N^\ell$, the 
cabling $c_\n(L)$ has each component of $L$ replaced
by $n_i$ parallel components, where parallel is defined
with respect to the Seifert framing of the link, \ie the one obtained using
a Seifert surface $\Sigma$ for $L$ and moving the link components to parallel
copies lying on a collar neighborhood of the boundary of the Seifert 
surface. Notice that if $L$ is an admissible
link, so that all the linking numbers are trivial, $lk(L_i,L_j)=0$,
then $c_\n(L)$ is also an admissible link. In fact,
the Seifert framing has the property that the parallel components
$L_{i,k}$ for $k=1,\ldots, n_i$ have zero linking numbers
$lk(L_{i,k},L_{i,k'})=0$, while  for $i\neq j$ the linking numbers 
remain $lk(L_{i,k},L_{j,k'})=lk(L_i,L_j)=0$.

\smallskip

In the variant we describe here below, we also impose an orientation
on the Hoste moves, so that an oriented Hoste move is directed
towards reducing the number of components of the link. We still
consider the identity move as part of the oriented Hoste moves.
Working with morphisms that are oriented Hoste moves will give
rise to a convolution algebra that is only a semigroupoid and
not a groupoid algebra and is no longer involutive.

\begin{defn}\label{cable2cat}
The category of cabled links $\cC\cC\cL$ has objects
of the form $(L,\fr, \n)$ where $L$ is an admissible link 
with components $L=L_1 \cup \cdots \cup L_\ell$, and framings
$\fr=(1/m_i)_{i=1,\ldots,\ell}$, and with $\n \in \N^\ell$ the datum
of a cabling $c_\n(L)$ of $L$ obtained as above, where all the components
$L_{i,k}$ of $c_\n(L)$ are given the same framing $1/m_i$. The morphisms
are given by sequences of ``cabled oriented Hoste moves". These are oriented
Hoste moves $c_\n(L)\stackrel{H}{\to} c_{\n'}(L')$ that
come from an oriented Hoste move $L\stackrel{\tilde H}{\to} L'$ in the sense
that one has a commutative diagram
$$ \xymatrix{ c_\n(L)\rto^H & c_{\n'}(L') \\
L \uto^{c_\n} \rto^{\tilde H} & L' \uto^{c_{\n'}} }$$
The corresponding convolution algebra $\cA_{\cC\cC\cL}$ is given by 
functions with finite support on the objects of this category with the
associative product
$$ (f_1 \star f_2)(H)=\sum_{H=H_1\circ H_2} f_1(H_1) f_2(H_2) $$
where the sum is over composition of cabled oriented Hoste moves.
\end{defn}

Notice that we can also modify the product to obtain an associative
algebra $\cA_{\cC\cC\cL,\cup}$ where the product also involves
the product by disjoint union of links as in \eqref{prod2ACL}. In particular, notice
that in this case, since we work with oriented Hoste moves, if $H=id$
and $H=H_1 H_2$ also $H_1=H_2=id$. Thus, in $\cA_{\cC\cC\cL,\cup}$,
when we restrict to $H=id$ we recover the group ring product
$$ (f_1\star f_2)(L,\n,\fr)=\sum_{L=L_1\cup L_2} f_1(L_1,\n_1,\fr_1) 
f_2(L_2,\n_2,\fr_2) $$
of the group of (cabled) framed links with the disjoint union operation.

\medskip

We then make the following observation.

\begin{lem}\label{divHoste}
Suppose given an admissible framed link $L=L_1\cup\cdots\cup L_\ell$,
and a positive integer $n\geq 1$ such that $n|m_i$ for all the integers $m_i$
such that $L_i$ has framing $1/m_i$. Then if a link $L'$ is obtained from $L$
by applying a sequence of oriented Hoste moves, one still has $n|m_i'$
for all the framings $1/m_i'$ of the components of $L'$. 
\end{lem}

\proof This is an immediate consequence of the fact that an oriented Hoste
move eliminates a component, inserting a number of full twists of the remaining
components, but does not alter the framings of all the remaining components.
\endproof

Notice that the same property clearly does not work for Hoste moves 
with the opposite orientation, namely those that replace a number of
full twists with an additional component. In fact, even if
$n|m_i$ for all components of $L'$ this does not imply that the same condition
will be satisfied for $L$ because, depending on the number of twists on the components
of $L'$, the extra component of $L$ may have a framing $1/m$ with $n\not| m$. 

\begin{lem}\label{Hostecable}
Let $(L,\fr)$ be an admissible framed link. Suppose given $n\geq 1$ such that
$n|m_i$ for all the framings $1/m_i$ of the components $L_i$. Let $L_n$ denote
the cabling $c_\n(L)$ of the framed link $L$ with $\n=(n,\ldots,n)$. Then, if there 
is an oriented Hoste move $H$ transforming $(L,\fr)$ into another framed link $(L',\fr')$, 
there exists an oriented cabled Hoste move $H_n$ relating the link $(L_n, n \fr)$
to the link $(L_n', n\fr')$.
\end{lem}

\proof It suffices to take the case of two links that differ by a single Hoste move
as in Figure \ref{HosteFig}, with $L$ having one more component
than $L'$, since the Hoste move is oriented. If we denote by $L_0$ the 
component of $L$ that is removed in the Hoste move, with framing $m_0$, 
then the remaining components  of the link acquire $-m_0$ full twists 
by effect of the Hoste move. 

By Lemma \ref{divHoste}  $n|m_i$ for all components of both $L$ and $L'$.
In the link $(L_n, n \fr)$, the component $L_0$ of $L$ is replaced by $n$ 
parallel components, all with framing $n/m_0$. On the other hand, the link 
$(L_n', n\fr')$ consists of the $n$-cabled $L'_n$ with components of framings 
$n\fr'=(n/m_i)$. Notice then that, if we apply a Hoste move to each of the $n$ 
components parallel to $L_0$ in $(L_n, n \fr)$ with framing $n/m_0$ we obtain a link
where we replace each of these components by $- m_0/n$ full twists. After
performing this operation on each of the components parallel to $L_0$
we therefore obtain $-m_0$ full twists on the $n$-cabled version of the
remaining components of $L\smallsetminus L_0$, with framings $n/m_i$,
which is $(L'_n, n\fr')$. 
\endproof

\bigskip
\bigskip

We can now define a transformation on the algebra $\cA_{\cC\cC\cL,\cup}$ given by
\begin{equation}\label{rhonLink}
\rho_n(f)(H) = \left\{ \begin{array}{ll} f(H_n) & n| m_i, \,\, \forall i, 
\\[2mm]
0 & \text{otherwise,} \end{array} \right. 
\end{equation}
where $H$ is a cabled oriented Hoste move
$$ (L,\n,\fr) \stackrel{H}{\to} (L',\n',\fr'), $$
the framing of $L$ is $\fr=(1/m_i)$, and 
$H_n$ is the cabled oriented Hoste move obtained as
in Lemma \ref{Hostecable} with
$$ (L, n\n, n\fr) \stackrel{H_n}{\to} (L', n\n', n \fr'). $$
Notice that $n\fr'$ is still an admissible framing by Lemma \ref{divHoste}.

\begin{prop}\label{rhonACL}
The $\rho_n$ of \eqref{rhonLink} give an action of
$\N$ by endomorphisms of $\cA_{\cC\cC\cL,\cup}$.
\end{prop}

\proof The product $(f_1\star f_2)(H)$ in $\cA_{\cC\cC\cL,\cup}$ 
is given by a sum over compositions $H=H_1\circ H_2$ and splittings into 
disjoint unions of components $H_1=H_{1,1}\cup H_{2,1}$ and
$H_2=H_{2,1}\cup H_{2,2}$. Consider a composition
$$ (L,\n,\fr) \stackrel{H_1}{\to} (\tilde L,\tilde\n,\tilde\fr) \stackrel{H_2}{\to} (L',\n',\fr'). $$
We show that $\rho_n (f_1 \star f_2)=\rho_n(f_1)\star \rho_n(f_2)$. One first can observe
that $\rho_n (f_1 \star f_2)$ is nonzero if and only if $L$ has the property that $n|m_i$ for all
the component. This is the same property that determines whether $\rho_n(f_1)$ is non-zero.
If $\rho_n(f_1)\neq 0$ then Lemma \ref{divHoste} implies that $\tilde L$ also has the
property that $n| \tilde m_i$ hence $\rho_n(f_2)\neq 0$. 
One then sees that, when non-zero, the product $\rho_n(f_1)\star \rho_n(f_2)$ gives
$$(\rho_n(f_1)\star \rho_n(f_2))(H)=\sum f_1(H_{n,1,1}) f_2(H_{n,2,2}),$$ 
where the sum is over the splittings
$$ H_n = H_{n,1}\circ H_{n,2}= (H_{n,1,1}\cup H_{n,1,2})\circ (H_{n,2,1}\cup H_{n,2,2}), $$
while $\rho_n (f_1 \star f_2)$ gives
$$ \rho_n (f_1 \star f_2)(H)= \sum f_1(\tilde H_{1,1}) f_2(\tilde H_{2,2}), $$
where here the sum is over all splittings
$$ H_n = \tilde H_1 \circ \tilde H_2 =(\tilde H_{1,1}\cup \tilde H_{1,2}) \circ
(\tilde H_{2,1} \cup \tilde H_{2,2}). $$
Since we use {\em cabled} oriented Hoste moves as morphisms, we know that
$\tilde H_{i,j}= H_{n,i,j}$ must be obtained as the n-cabling of a Hoste move
$H_{i,j}$, so that the two sums agree.
\endproof

This gives a partial answer, in this more general setting of convolution
algebras of links and Hoste moves, to the question of defining an action
of $\N$ by endomorphisms of a ring generalizing $\Z HS$. It would be
interesting to construct actions of $\N$ on the other algebras
$\cA_{\cC\cL}$, $\cA_{\cC\cL,\cup}$, $\cA_{\cC\cL^{(2)}}$, $\cA_{\cC\cL^{(2)},\cup}$
described above.

\bigskip

We discuss briefly another possible approach towards constructing
endomorphisms actions on 3-manifolds, using closed braid representations
of links and endomorphisms of braid groups. We only sketch the
construction, but we do not fully develop these ideas in the present paper.

The problem can once again be summarized as trying to construct an action on (framed)
links, which is well defined under the usual equivalence of planar diagrams
via Reidemeister moves, and which can be constructed locally at crossings,
and which moreover preserves the Fenn--Rourke moves (or Rolfsen moves 
in the case of rational framings), so that it induces a 
well defined action at the level of the 3-manifolds, mapping equivalent
surgery presentations to equivalent surgery presentations.

First recall quickly the relation between links and braids (see the references
in \cite{BirMe1}). It is well known by a
classical result of Alexander that all link types can be represented by
closed $N$-braids, for some $N$. A closed braid representation of a link is
obtained by a choice of a braid axis $A$ and an open braid is obtained from the
link by cutting it open along a half-plane $H_\theta$ originating from the axis $A$. 
Link types correspond to equivalences of braid representations via conjugation and
Markov stabilization/destabilization operations. These give a reformulation
in braid terms of the usual equivalence of link diagrams via Reidemeister moves.
It is known by a result of Morton that an equivalence class of closed braid
diagrams corresponding to a link type consists of an infinite family of conjugacy
classes in braid groups $B_N$, for varying $N$, or in $B_\infty$, the union 
of the $B_N$. The main result of \cite{BirMe1} extracts from this set of conjugacy
classes a finite number of elements of {\em minimal complexity}, with respect to
a complexity function. The complexity function constructed in \cite{BirMe1}
is obtained from geometric data, which in turn admit a combinatorial formulation.
Namely, given a link $L$ and a choice of an axis $A$, one considers an oriented
(in general not connected) embedded Seifert surface of maximal Euler characteristic with 
boundary $\partial \Sigma = L$ and a fibration $H$ of  $\R^3\smallsetminus A$ with
fibers $H_\theta$, with good assumptions on how tubular neighborhoods of $A$ and 
the fibers $H_\theta$ intersect $\Sigma$. The complexity function is then defined
as the triple $(\# (L\cap H_\theta), \# (A\cap \Sigma), \# T(H,\Sigma))$, where
$T(H,\Sigma)$ is the set of singular points of the foliation induced on $\Sigma$
by $H$. The latter is known as the {\em braid foliations}. The complexity function 
defined in this way takes values in $\Z_+ \times \Z_+ \times \Z_+$
ordered lexicographically, and it only depends on the conjugacy class in $B_N$
determined by the data $(L,A)$, with $N=\# (L\cap H_\theta)$. 
The advantage of using representatives of minimal complexity is that one can
avoid Markov stabilization moves, which is the basis for the more recent ``Markov's theorem
without stabilization" of \cite{BirMe2}. We discuss below the reason why this approach
may be useful to our problem.

We first need to recall also some properties of braid groups and their endomorphisms.
It is known by \cite{BeMa} that braid groups (for $N\geq 4$) 
are Hopfian but not co-Hopfian, where the Hopfian property follows 
from the fact that braid groups are linear (see \cite{Bigelow}) and it means that a
surjective endomorphism of a braid group $B_N$ is in fact an automorphism.
Automorphisms are rather easy to describe, as they are combinations of
the inner ones and the inversion that sends each generator to its inverse.
The fact that the braid groups are not co-Hopfian means that there are
injective endomorphisms that are not automorphisms. In fact, a complete
description of injective endomorphisms is given in \cite{BeMa}. They come
from homeomorphisms $h: D_N \to D_N$ of the $N$-punctured disk, by setting
\begin{equation}\label{rhohm}
 \rho_{h,m}(s_a)= s^{\epsilon(h)}_{h(a)} T_N^m, 
\end{equation}
where $\epsilon(h)=\pm 1$ according to the orientation preserving/reversing 
property of $h$, while $s_a$ are the half-twists switching two punctures 
along an arc $a$,  which generate the group $B_N$, and $T_N$ is the 
generator of the center of $B_N$, namely the full twist of the $n$ strands,
$T_N=(s_1\cdots s_{N-1})^N$. 

Suppose then that one defines an operation on closed braid diagrams of links
by considering a simple class of such endomorphisms of braid groups, namely
those of the form $\rho(s_a)=s_a T_N^m$, for some $m\in \Z$.  This means that
one obtains a new braid diagram from a given one by replacing the open braid
$\gamma\in B_N$ by $\rho(\gamma) \in B_N$ and closing the diagram again.
A different opening and closing of the braid diagram gives rise of an element
in the same conjugacy class $C(\gamma)$ in $B_N$, which will in turn be mapped 
to an element in the conjugacy class $C(\rho(\gamma))$ so that the action is
well defined on conjugacy classes. In fact, even more precisely, since the
element $T_N$ is in the center, one has $\rho(\alpha \gamma \alpha^{-1}) =
\alpha \rho(\gamma) \alpha^{-1}$, since $\rho(\alpha) = \alpha T_N^{k(\alpha)}$
with $k(\alpha)= m \ell(\alpha)$, where $\ell(\alpha)$ is the length of $\alpha$ defined
as the image under the homomorphism $\ell: B_N \to \Z$ that sends all the generators
$s_i$ to $1$.

Thus, we can define in this way an operation of $\rho \in \End(B_N)$ on 
conjugacy classes in $B_N$. For it to induce an action on the link types
we would need it to be compatible with the Markov equivalence between
braid representatives, realized by conjugations and Markov de/stabilizations.
It is easy to see that the operation is well defined under conjugation, or
equivalently under the second Reidemeister move. In fact, if two braids
$\gamma$ and $\gamma'$ differ by a twist $s_i$ and an opposite twist $s_i^{-1}$,
their images $\rho(\gamma)$ and $\rho(\gamma')$ will again differ in the
same way, again because of the fact that $T_N$ is in the center.
However, the case of Markov de/stabilizations
is more delicate, since these change the number of braids and the extra
crossing that corresponds to the first Reidemeister move has the effect
that, if $\gamma \in B_N$ and $\gamma' \in B_{N+1}$ differ by a Markov
stabilization, that is $\gamma' = \gamma s_N$, 
then the images satisfy
$$ \rho_{N+1}(\gamma') = \rho_{N+1}(\gamma) \rho_{N+1}(s_N) 
= \rho_{N+1}(\gamma) s_N T_{N+1}^{m_{N+1}} $$
which is off from a Markov stabilization by $m_{N+1}$ full twists.

This is where the existence of minimal complexity representatives may
become useful in trying to make this kind of action well defined on link
types. Selecting the finitely many minimal complexity representative
conjugacy classes, one eliminates the need for Markov de/stabilizations.
However, iterating the transformation defined by $\rho$ on these classes
requires the image of a minimal complexity representative to be still
of minimal complexity for the link type obtained after applying the 
transformations. The new link is in each case obtained from the old
one simply by inserting (it does not matter where, since it is an element in
the center of $B_N$) a number $m_N \ell(\gamma)$ of full twists.

A first question is the following: if $\gamma$ and $\gamma'$ in $B_N$
represent the same link type, do the corresponding elements 
$\gamma T_N^{m_N \ell(\gamma)}$ and $\gamma' T_N^{m_N \ell(\gamma')}$ 
also represent the same link type? Notice that a positive answer to this
question would follow, if we knew that two equivalent
$\gamma$ and $\gamma'$ will have the same length. This is certainly 
not the case if one is allowed to use Markov de/stabilizations.  
However, it is a well known conjecture of Jones (see \cite{Jones} p.357) 
that the algebraic crossing number (or writhe) of a minimal braid (with braid 
index  $N$ equal to the minimum value for the given knot type) is a link 
invariant, hence if $\gamma$ and $\gamma'$ in $B_N$ are braid representatives 
with $N$ minimal, then $\ell(\gamma)=\ell(\gamma')$. Evidence for this 
conjecture is given, for instance, in \cite{Kawa}. 

If we assume the Jones conjecture that minimal braids have a unique writhe,
then it is clear that, if $\gamma$ and $\gamma'$ are two such braid representatives
of a given link type, then the transformed $\gamma T_N^{m_N \varpi}$ and 
$\gamma' T_N^{m_N \varpi}$ also represent the same link type, 
where $\varpi$ is the writhe of both $\gamma$ and $\gamma'$. In fact, both
are obtained from equivalent braids by composing them with the same number $m_N \varpi$
of full twists of the $N$ strands, which still yields equivalent braids.  Minimal
complexity braid representatives are in particular minimal braids, so if we apply
our transformation to these elements only and we assume the Jones conjecture
holds, we do get a positive answer to the question above, which in turn implies
that the link type of the image is well defined regardless of which minimal
complexity representative we choose among the finitely many available for
a given link type.
 
We can then ask the next question, which is needed in order to be able to
iterate the transformation we defined by applying the braid group endomorphism 
$\rho$ to braids that are of minimal complexity in the sense of \cite{BirMe1}.

\begin{ques}\label{quescomplexity}
For what classes of links does the operation 
$\gamma \mapsto \gamma T_N^{m_N \ell(\gamma)}$ 
preserve the minimal complexity property?
\end{ques}

To understand the answer to this question one needs to understand
how inserting full twists in a given minimal complexity braid diagram
alters the properties of the Seifert surface $\Sigma$ and of the braid foliation 
determined on it by the fixed fibration $H$. Notice that neither the choice of the
axis $A$ nor the number $N$ are changed by the operation above.
However, the fact that $N$ is unchanged does not a priori imply that
it will still be minimal for the new link type. We do not at present know of an
answer to Question \ref{quescomplexity}. Thus, we can, for the sake of this 
discussion, define $\cL_{\min}$ to be the class of link types that admit a minimal
complexity braid representative for which the answer to the question is
positive. Notice that the set $\cL_{\min}$ is non-empty. In fact, it contains at 
least the torus knots $T(a,b)$, which have a unique minimal braid conjugacy class,
which is determined by the braid $(s_1 \cdots s_{a-1})^b$. On these minimal braids
the action described above is of the form
$$ \rho: (s_1\cdots s_{a-1})^b \mapsto  (s_1\cdots s_{a-1})^b T_a^{m b (a-1)} =
(s_1 \cdots s_a)^{b (1+m a(a-1))}, $$
which has the effect of transforming $T(a,b)$ into $T(a,b (1+m a(a-1)))$.

Thus, modulo obtaining a better answer to Question \ref{quescomplexity},
we will be restricting our attention to the subring $\cL_{\min}$ of $\cL$.
One can use the complexity arguments of \cite{BirMe4} to check that the
ring multiplication, which is given by the disjoint union of links, is well defined 
on this subset. 

Notice then that if we have two endomorphisms of the same type
$$ \rho_{N,n_i}(\gamma) = \gamma T_N^{n_i \ell(\gamma)}, \ \ \  i=1,2, $$
they satisfy
$$ \rho_{N,n_2}\rho_{N,n_1}(\gamma) =
\gamma T_N^{(n_1 +n_2 + n_1 n_2  N(N-1))\ell(\gamma)} =\rho_{N,n_1}\rho_{N,n_2}(\gamma). $$
Thus, upon choosing an endomorphism $\rho_{N,n_p}$ for each generator $p$ of the
multiplicative semigroup $\N$  one obtains an action of the abelian semigroup $\N$ on
the links of $\cL_{\min}$.

We are then left with the question of implementing the surgery equivalence of links.
Since we are working with closed braid representatives, it is best to think of the
equivalence relation generated by the Fenn--Rourke moves (let us restrict to the original
case of integer framings for simplicity here) in terms of braid groups. This was done 
in \cite{Levine}, see also \cite{GaLe} and \cite{Habeg}. We give here the formulation
of \cite{Habeg}, which is closer to the form we need. In this form, the set of
surgery equivalence classes of $N$-strand braid representatives is the quotient of
the braid group $B_N$ by its normal subgroup $P(N)_3$, which is the third term in
the lower central series of the pure braid group $P(N)$ defined by the exact sequence
$$ 1 \to P(N) \to B_N \to S_N \to 1, $$
which maps the braid group to the symmetric group acting as permutations of the
endpoints of the braid. This point of view has the advantage of formulating the
equivalence relation in a way that is manifestly compatible with the braid group
structure and compositions of braids. Thus, in particular, if our transformation is
given by the braid group endomorphism $\rho_N: \gamma \mapsto 
\gamma T_N^{m_N \ell(\gamma)}$, then if $\gamma$ and $\gamma'$ are surgery
equivalent so that $\gamma^{-1}\gamma' \in P(N)_3$, we have $\rho_N(\gamma)^{-1}
\rho_N(\gamma')= \gamma^{-1} \gamma' T_N^{m_N (\ell(\gamma')-\ell(\gamma))}$.
In particular, if $\gamma$ and $\gamma'$ are surgery equivalent and have 
$\ell(\gamma)=\ell(\gamma')$ then $\rho_N(\gamma^{-1}\gamma')= 
\gamma^{-1}\gamma' \in P(N)_3$, so their images under $\rho_N$ are also 
surgery equivalent. This only takes care of implementing surgery equivalence 
among representatives with the same number of strands and with the same writhe.
Since $\gamma$ and $\gamma'$ here belong to different link types, the Jones
conjecture for their minimal braid representatives need no longer imply that
$\ell(\gamma)=\ell(\gamma')$ even if they have the same $N$. Thus, this argument
alone does not suffice to implement surgery equivalence. A more refined argument
along these lines may be obtained, in the case of algebraically split links, 
in terms of the Milnor invariants obtained from triple intersections of Seifert surfaces
as in \cite{Levine}. 

A different way to approach the question of implementing the surgery moves is
to proceed as in \cite{Wenzl}, where one works with multiplicative invariants of
links and obtains from these invariants of 3-manifolds only after a suitable averaging
over repeated cabling operations and a thermodynamical limit. We only make some
preliminary observations on what questions one needs to address for
this approach to be indeed compatible with the action on links defined using the 
braid group endomorphisms and the minimal complexity braid representatives. 
Essentially, we are again asking a question similar to Question \ref{quescomplexity}
above, on whether certain operations, in this case cabling, preserve the minimal
complexity property. There are some indications to this effect, at least for what
concerns the first entry of the complexity function of \cite{BirMe1}, that is, the
braid index. In fact, it is shown in \cite{Will} (see also \cite{Kawa}) that the braid 
index of the $n$-cabling of a link $L$ is $n$-times the braid index of $L$, so 
the $n$-cabling of a minimal braid for $L$ is a minimal braid for the $n$-cabling of $L$. 
Also, by \cite{Kawa}, if the Jones conjecture holds for $L$, it also holds for its cabled
versions. The question of whether this property of the braid index 
extends to the rest of the complexity function of
\cite{BirMe1} depends on the Seifert surface of the cabled link and its braid foliations.
Since this will locally look like $n$ copies of the original Seifert surface and foliations,
it seems likely that the rest of the complexity function should also behave well under 
the cabling operations, but this requires a more careful investigation. If one can indeed
extend in this way the action from a link to all its $n$-cabled version, one can take 
the resulting action as being the one that induces the action on the 3-manifold
invariants, when the latter are constructed according to the procedure of \cite{Wenzl}.

One can also, in this approach based on braid representatives of links,
introduce the same kind of categorical constructions discussed in the
first part of this section, with categories of braid representatives, 1-morphisms
given by surgery equivalences and and 2-morphisms  by
Markov de/stabilizations and conjugations. For a categorical framework
using braid representatives of links see for instance \cite{Yetter}.

\medskip

At this stage it is not clear how far one can develop this point of view. It would be
interesting to interpret the universal WRT invariant as a morphism from a suitable
``endomotive of homology 3-spheres", which means a crossed product
$\cA\rtimes \N$ by an action of $\N$ on an algebra $\cA$ of admissible links and Hoste moves
implementing the surgery equivalence between them, to the endomotive constructed out of the
one-variable Habiro ring earlier in this paper. Such a construction would make it
possible to import the techniques of noncommutative
geometry developed in \cite{CCM} for the theory of endomotives to 
the context of 3-manifold invariants.

\end{document}